\documentclass[12pt]{article}
\usepackage[a4paper]{geometry}
\usepackage{amsthm}
\usepackage{amsmath}
\usepackage{amssymb}
\usepackage{amsfonts}
\usepackage{graphicx}
\usepackage{graphics}
\usepackage[bf,SL,BF]{subfigure}

 \numberwithin{equation}{section}
\def\eref#1{(\ref{#1})}
\def\N{\mathbb{N}}
\def\P{\mathbb{P}}
\def\R{\mathbb{R}}
\def\T{\mathcal{T}}
\def\No{\mathcal{N}}
\def\S{\mathcal{S}}
\def\Y{\mathcal{Y}}

\def\To{\mathbb{T}}

\def\Z{\mathcal{Z}}
\def\M{\mathcal{M}}

\def\B{\mathcal{B}}
\def\ED{\mathcal{E}}

\def\U{\mathcal{U}}
\def\I{\mathcal{I}}
\def\J{\mathcal{J}}

\def\<{\big\langle}
\def\>{\big\rangle}

\def\Vol{\operatorname{Vol}}
\def\diiv{\operatorname{div}}
\def\dist{\operatorname{dist}}

\def\Tr{\operatorname{Trace}}

\def\det{\operatorname{det}}

\def\esssup{{\operatorname{esssup}}}
\def\supp{{\operatorname{support}}}

\newtheorem{Lemma}{Lemma}[section]

\newtheorem{Theorem}{Theorem}[section]

\theoremstyle{remark}
\newtheorem{Remark}{Remark}[section]

\theoremstyle{definition}

\newtheorem{Definition}{Definition}[section]

 \theoremstyle{definition}
  \newtheorem{example}{Example}

\newcommand{\figbox}[1]{%
  \fbox{%
    \vbox to 1in{%
    \vfil
    \hbox to 2in{%
      \hfil
      #1%
      \hfil}%
    \vfil}}}

\newcommand{\goodgap}{%
  \hspace{\subfigcapskip}}

\begin{document}

\title{Metric based up-scaling\protect\footnotetext{AMS 1991 {\it{Subject Classification}}. Primary 34E13,35B27 ; secondary
68P30,  60F05, 35B05.} \protect}          

\author{Houman Owhadi\footnote{
 California Institute of Technology
Applied \& Computational Mathematics, Control \& Dynamical systems,
MC 217-50 Pasadena , CA 91125, owhadi@caltech.edu}  and Lei
Zhang\footnote{
 California Institute of Technology
Applied \& Computational Mathematics MC 217-50 Pasadena , CA 91125,
zhanglei@acm.caltech.edu}}
\date{November 15, 2005}
\maketitle \abstract{We consider divergence form elliptic operators
in dimension $n\geq 2$ with $L^\infty$ coefficients. Although
solutions of these operators are only H\"{o}lder continuous, we show
that they are differentiable ($C^{1,\alpha}$) with respect to
harmonic coordinates. It follows that numerical homogenization can
be extended to situations where the medium has no ergodicity at
small scales and is characterized by a continuum of scales by
transferring a new metric in addition to traditional averaged
(homogenized) quantities from subgrid scales into computational
scales and error bounds can be given. This numerical homogenization
method can also be used as a compression tool for differential
operators.
 }

\section{Introduction and main results}
Let $\Omega$ be a bounded and convex domain of class $C^2$. We
consider the following benchmark PDE
\begin{equation}\label{ghjh52}
\begin{cases}
-\diiv\big(a(x)\nabla u(x)\big)=g\quad \text{in}\quad \Omega \\
u=0 \quad \text{in}\quad \partial \Omega.
\end{cases}
\end{equation}
Where $g$ is a function in $L^{\infty}(\Omega)$. And $x\rightarrow
a(x)$ is a mapping from $\Omega$ into the space of positive definite
symmetric matrices. We assume $a$ to be  symmetric, uniformly
 elliptic  with entries in $L^\infty(\Omega)$.
\begin{itemize}
\item Is it possible to up-scale \eref{ghjh52}?
\end{itemize}
Homogenization theory (\cite{BeLiPa78}, \cite{JiKoOl91}) allows us
to do so by transferring bulk (averaged) information from sub-grid
scales into computational scales. This transfer from a numerical
homogenization point of view is justified under two fundamental
assumptions:
\begin{itemize}
\item Ergodicity at small scales, and
\item  scale separation.
\end{itemize}
Can we get rid of these assumptions?
\begin{itemize}
\item Can we numerically homogenize \eref{ghjh52} when $a$ is arbitrary? in particular
when $a$ is characterized by a continuum of scales with no
ergodicity at small scales?
\end{itemize}
What do we mean by numerical homogenization when $a$ is arbitrary?
It is important to recall that F. Murat and L. Tartar's theory of
H-convergence \cite{MR1493039} provides a mathematical framework for
analysis of composites in complete generality, without any need for
geometrical hypotheses such as periodicity or randomness. This
theory is based on a powerful tool called compensated compactness or
div-curl lemma introduced in the 70's by Murat and Tartar
\cite{MR506997}, \cite{MR584398} which has been characterized by a
wide range of applications and refinements \cite{MR1225511}. Here we
consider homogenization from a slightly different point of view: we
want to solve
 \eref{ghjh52} on a coarse mesh and we want to understand which information should be transferred
 from fine scales to coarse scales when the entries of $a$ are
 arbitrary. For that purpose we need a new form of compensation given in section \ref{jsjh61}.

It is important to observe that if one needs to solve \eref{ghjh52}
only one time (with one $g$) the method proposed here does not
reduce the number of numerical operations\footnote{ It will be shown
in \cite{OwZh05b} that for parabolic and reaction-diffusion
equations the situation is different:
 the methodology introduced in this paper can be used to reduce
the number of operations even when one needs to solve these
equations only once.}. Indeed we need to compute \eref{ghjh52}
$n$-times ($n$ being the space dimension) with $0$ in the right hand
side in \eref{ghjh52} and linear boundary conditions. However if one
needs to solve \eref{ghjh52} for a large number ($M$) of different
right hand sides $g$ ($M>>n$, which would happen if one tries to
optimize specific properties of $u$ with respect to $g$) then the
methods proposed here have a practical use since they basically say
that after solving \eref{ghjh52} $n$ times it is sufficient to solve
\eref{ghjh52} on a coarse mesh (with $N^{0.01}$ nodes instead of $N$
for instance) $M$ times.

Let us recall that fast methods based on hierarchical matrices are
available \cite{MR1993936, BeCh05, Beb05, Bebe05, Beben05} for
solving \eref{ghjh52} in $O\big(N (\ln N)^{n+3}\big)$ operations.

 Finally the point of view of this paper is to observe the
``redundancy'' of solutions of \eref{ghjh52} at small scales rather
than ``fast computation''. Indeed it is not obvious that
\eref{ghjh52} can be homogenized when the medium does not satisfy
the usual periodicity, ergodicity or scale-separation assumptions.
Moreover if \eref{ghjh52} can indeed be homogenized it is important
to understand what minimal quantity of information should be kept
from small scales to obtain an accurate homogenized operator? Once
the correct coarse parameters are identified (the bulk quantities
and the up-scaled metric), one can try to model, estimate or
simulate them but the first step is to identify them.

\subsection{A new form of compensation}\label{jsjh61}
To introduce the new form of compensation, we need to introduce the
so called $a$-harmonic coordinates associated to \eref{ghjh52}, i.e.
 the weak solution of the following boundary value problem
\begin{equation}\label{dgdgfsghsf62}
\begin{cases}
\diiv a \nabla F=0 \quad \text{in}\quad \Omega\\
F(x)=x \quad \text{on}\quad \partial \Omega.
\end{cases}
\end{equation}
By \eref{dgdgfsghsf62} we mean that $F$ is a $n$-dimensional vector
field $F(x)=\big(F_1(x),\ldots,F_n(x)\big)$ such that each of its
entries satisfies
\begin{equation}
\begin{cases}
\diiv a \nabla F_i=0 \quad \text{in}\quad \Omega\\
F_i(x)=x_i \quad \text{on}\quad \partial \Omega.
\end{cases}
\end{equation}
The new compensation phenomenon  is controlled by the following
object:
\begin{Definition}
We write
\begin{equation}
\sigma:= {^t \nabla F a\nabla F}.
\end{equation}
\end{Definition}
We write $\mu_{\sigma}$ the anisotropic distortion of $\sigma$
defined by
\begin{equation}
\mu_{\sigma}:=\esssup_{x\in \Omega}
\Big(\frac{\lambda_{\max}\big(\sigma(x)\big)}{\lambda_{\min}\big(\sigma(x)\big)}\Big).
\end{equation}
Where $\lambda_{\max}(M)$ ($\lambda_{\min}(M)$) denote the maximal
(minimal) eigenvalue of $M$.
\begin{Definition}
In dimension $n=2$, we say that $\sigma$ \it{is stable} if and only
if $\mu_{\sigma}<\infty$ and there exist a constant $\epsilon>0$
such that $\big(\Tr(\sigma)\big)^{-1-\epsilon}\in L^{1}(\Omega)$
\end{Definition}
\begin{Remark}
According to \cite{MR1892102} in dimension two if $a$ is smooth then
$\sigma$ is stable. According to \cite{MR2001070}, $F$ is always an
homeomorphism in dimension two even with $a_{i,j}\in
L^\infty(\Omega)$.
\end{Remark}

\begin{Theorem}\label{th2}
Assume  that $\sigma$ is stable and $n=2$. Then there exist
constants $\alpha>0$ and  $C>0$ such that $(\nabla F)^{-1} \nabla u
\in C^\alpha(\Omega)$ and
\begin{equation}\label{ksjsh5651}
\big\|(\nabla F)^{-1} \nabla u\big\|_{C^\alpha(\Omega)} \leq C
\|g\|_{L^\infty(\Omega)}.
\end{equation}
\end{Theorem}

\begin{Remark}
The constant $\alpha$ depends on
$\Omega,\lambda_{\max}(a)/\lambda_{\min}(a)$  and $\mu_{\sigma}$.
The constant $C$ depends on the constants above, $\lambda_{\min}(a)$
and
$\big\|\big(\Tr(\sigma)\big)^{-1-\epsilon}\big\|_{L^{1}(\Omega)}$.
We use the notation $\lambda_{\max}(a):=\sup_{x\in \Omega}
\sup_{|\xi|=1}{^t\xi a\xi}$ and \\$\lambda_{\max}(a):=\inf_{x\in
\Omega} \inf_{|\xi|=1}{^t\xi a\xi}$.
\end{Remark}
\begin{Remark}
If one considers a sequence $a_\epsilon$ such that
$\mu_{\sigma_{\epsilon}}$ and\\
$\big\|\big(\Tr(\sigma_\epsilon)\big)^{-1-\epsilon}\big\|_{L^{1}(\Omega)}$
are uniformly bounded away from $\infty$ and
$\lambda_{\min}(a_\epsilon)$ and $\lambda_{\max}(a_\epsilon)$ are
uniformly bounded away from $0$ and $\infty$ then \eref{ksjsh5651}
is uniformly true. If we consider a periodically oscillating
sequence $a_{\epsilon}(x)=a(\frac{x}{\epsilon})$,
$\mu_{\sigma_1}<\infty$ and
$\big\|\big(\Tr(\sigma_1)\big)^{-1-\epsilon}\big\|_{L^{1}(\Omega)}<\infty$
then \eref{ksjsh5651} is uniformly true.
\end{Remark}

\begin{Remark}
It is easy to check from the proof of \ref{th2} that if
$\big(\Tr(\sigma_\epsilon)\big)^{-1}\in L^{q}(\Omega)$ (with $q>2$)
then it is possible to replace $\|g\|_{L^\infty(\Omega)}$ by
$\|g\|_{L^p(\Omega)}$ in \eref{ksjsh5651} with $p$ depending on $q$.
More precisely if $\big(\Tr(\sigma_\epsilon)\big)^{-1}\in
L^{\infty}(\Omega)$ then it is possible to replace
$\|g\|_{L^\infty(\Omega)}$ by $\|g\|_{L^{2+\epsilon}(\Omega)}$ in
\eref{ksjsh5651} for any $\epsilon >0$.
\end{Remark}

\begin{Remark}
We don't need $a$ to be symmetric to obtain \ref{th2} but for the
clarity of the presentation we have chosen to restrict ourselves to
that case. We have not analyzed Neumann boundary conditions, this
will done in a further work.
\end{Remark}

\begin{figure}[htbp]
  \begin{center}
    \subfigure[$a$ in log scale.
\label{particulsslll}]
          {\includegraphics[width=0.3\textwidth,height= 0.3\textwidth]{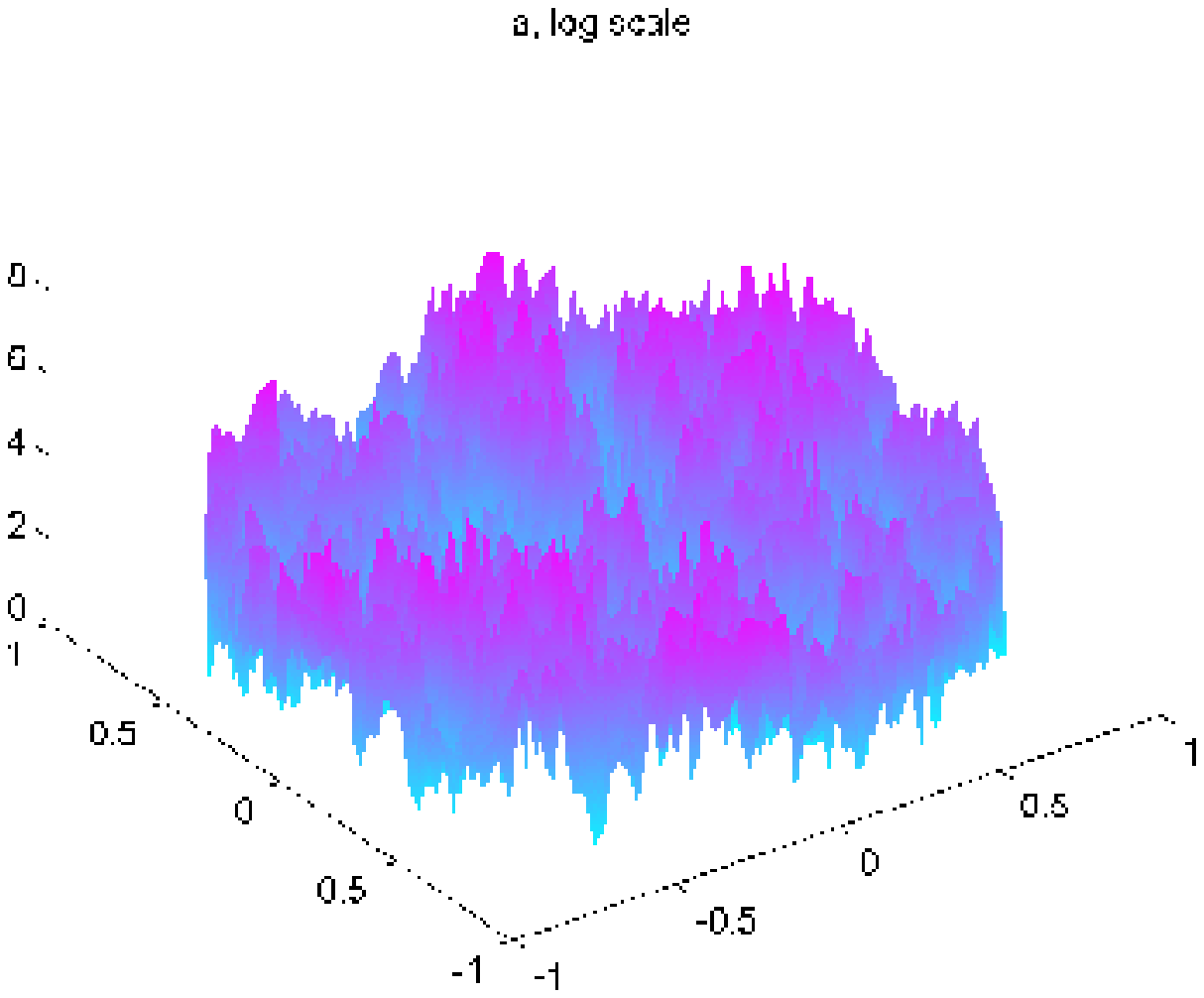}}
     \goodgap
    \subfigure[$\nabla u$.\label{figfffrrr}]
          {\includegraphics[width=0.3\textwidth,height= 0.3\textwidth]{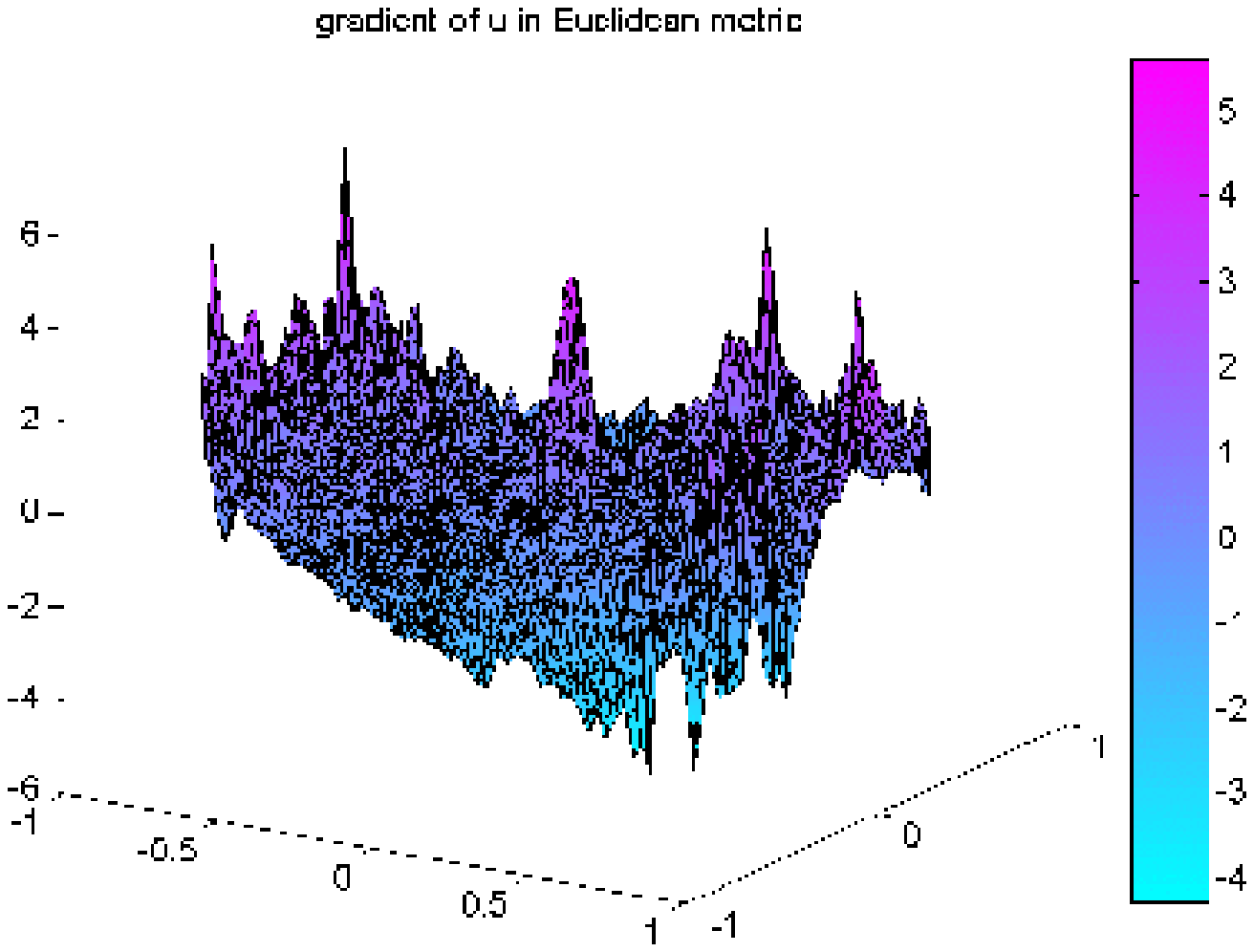}}
    \goodgap
   \subfigure[$(\nabla F)^{-1}\nabla u$.\label{figdjdueAfirstB}]
        {\includegraphics[width=0.3\textwidth,height= 0.3\textwidth]{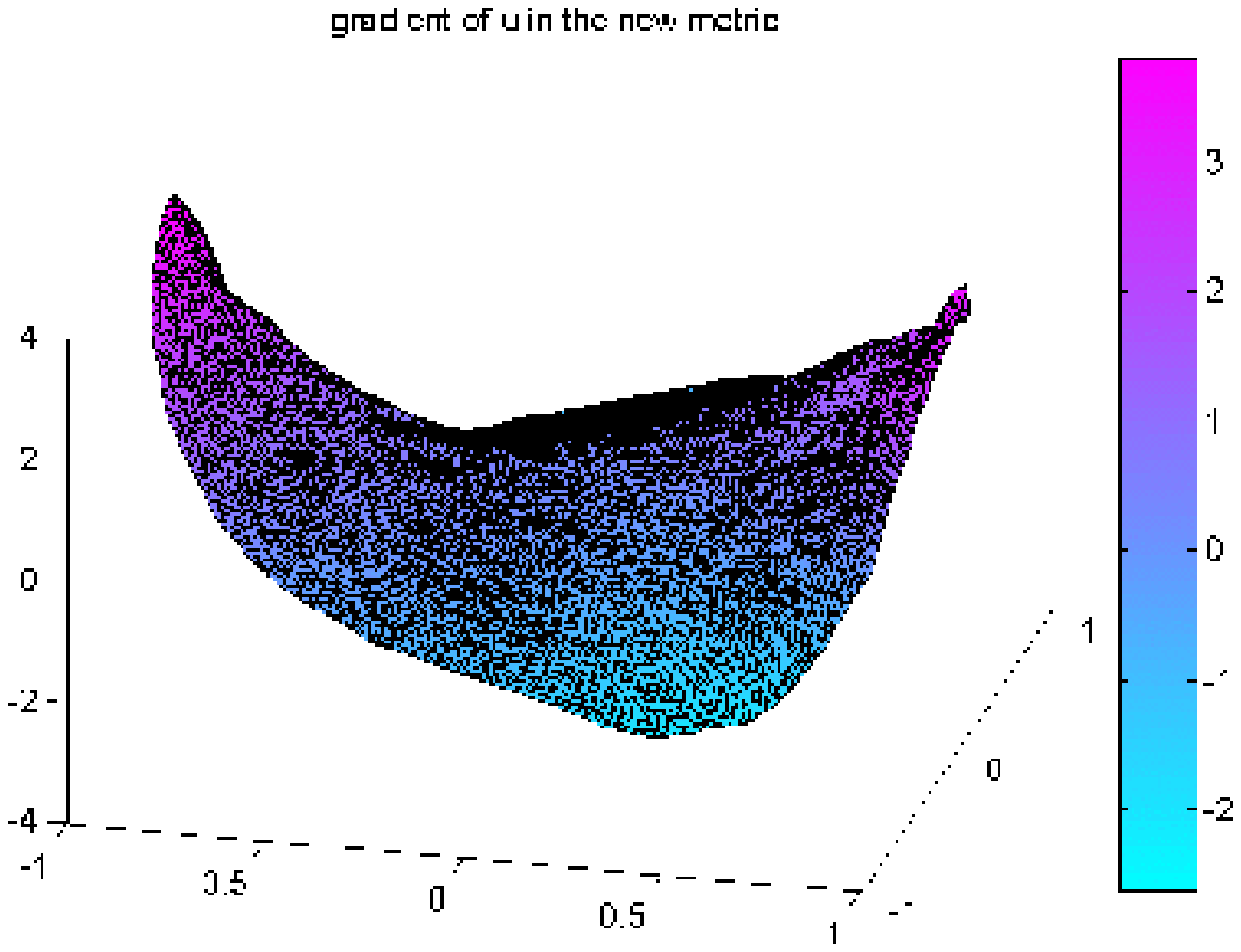}}\\[-10pt]
    \caption{Change of metric on the disk.}
    \label{figyyassAfirst}
\end{center}
\end{figure}

\begin{figure}[htbp]
  \begin{center}
    \subfigure[$a$ in log scale.
\label{partwial}]
          {\includegraphics[width=0.4\textwidth,height= 0.3\textwidth]{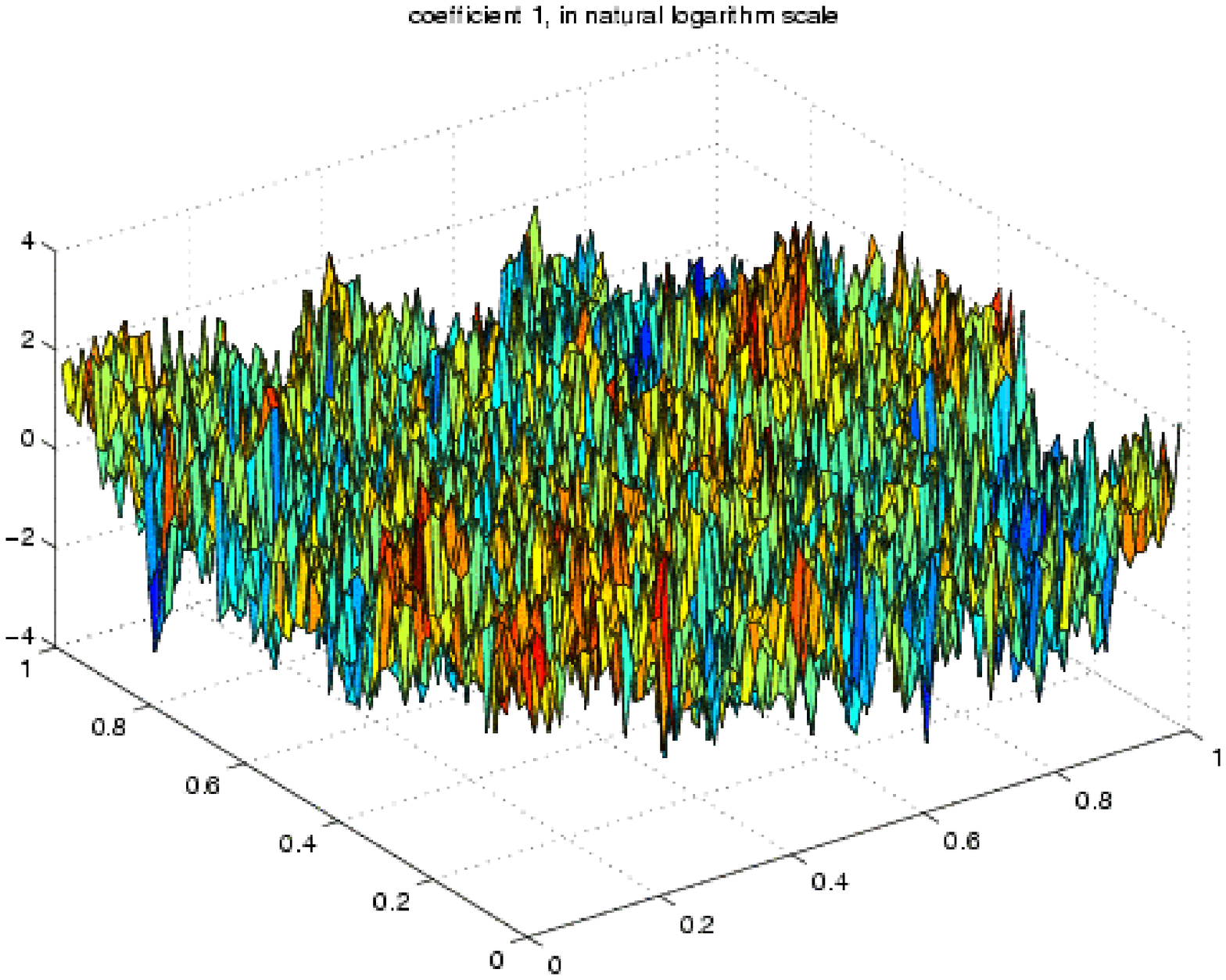}}
     \goodgap
    \subfigure[one of the entries $\nabla F$.\label{figafrwrr}]
          {\includegraphics[width=0.4\textwidth,height=
          0.3\textwidth]{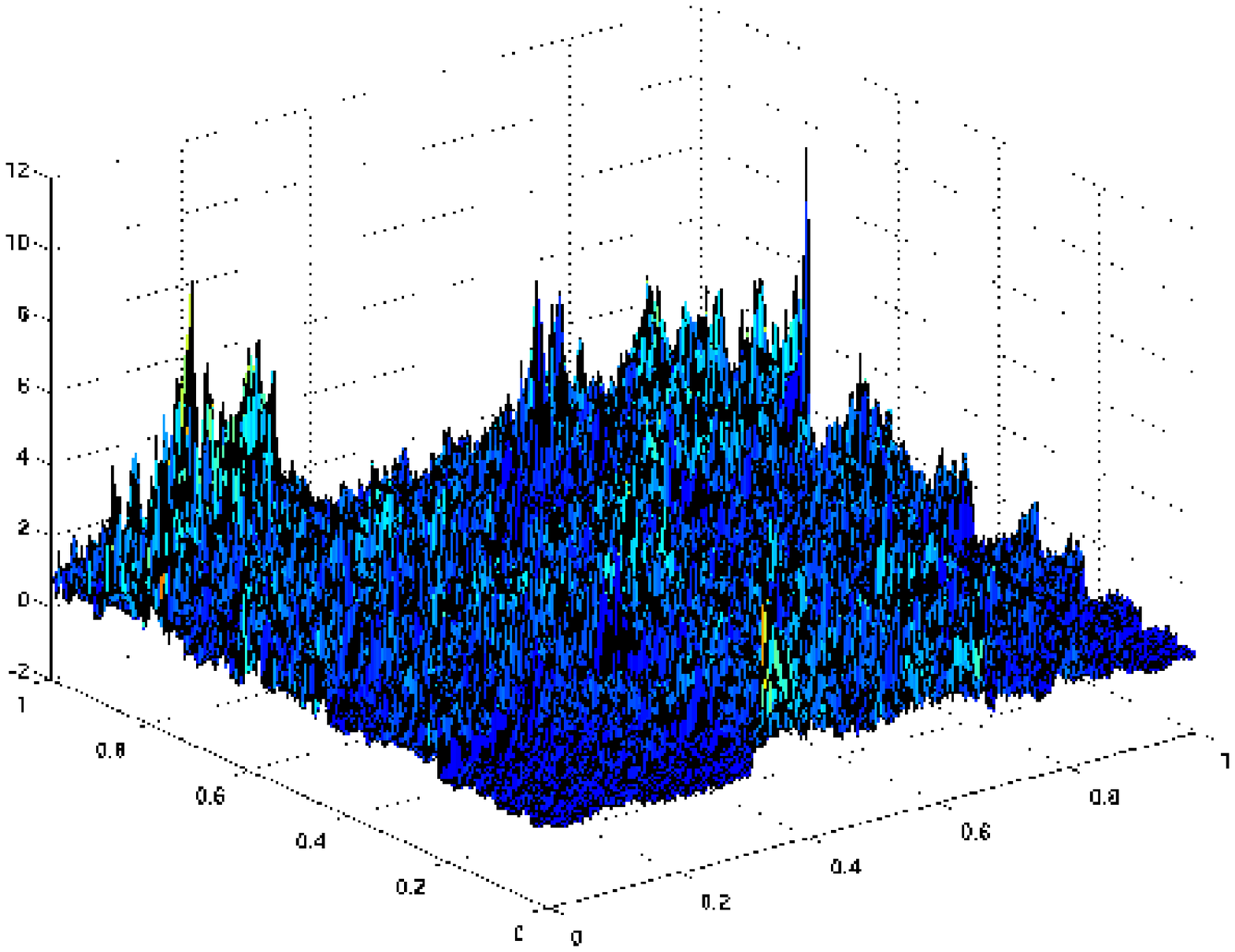}}\\
    \subfigure[$\nabla u$.
\label{partiawwl}]
          {\includegraphics[width=0.4\textwidth,height= 0.3\textwidth]{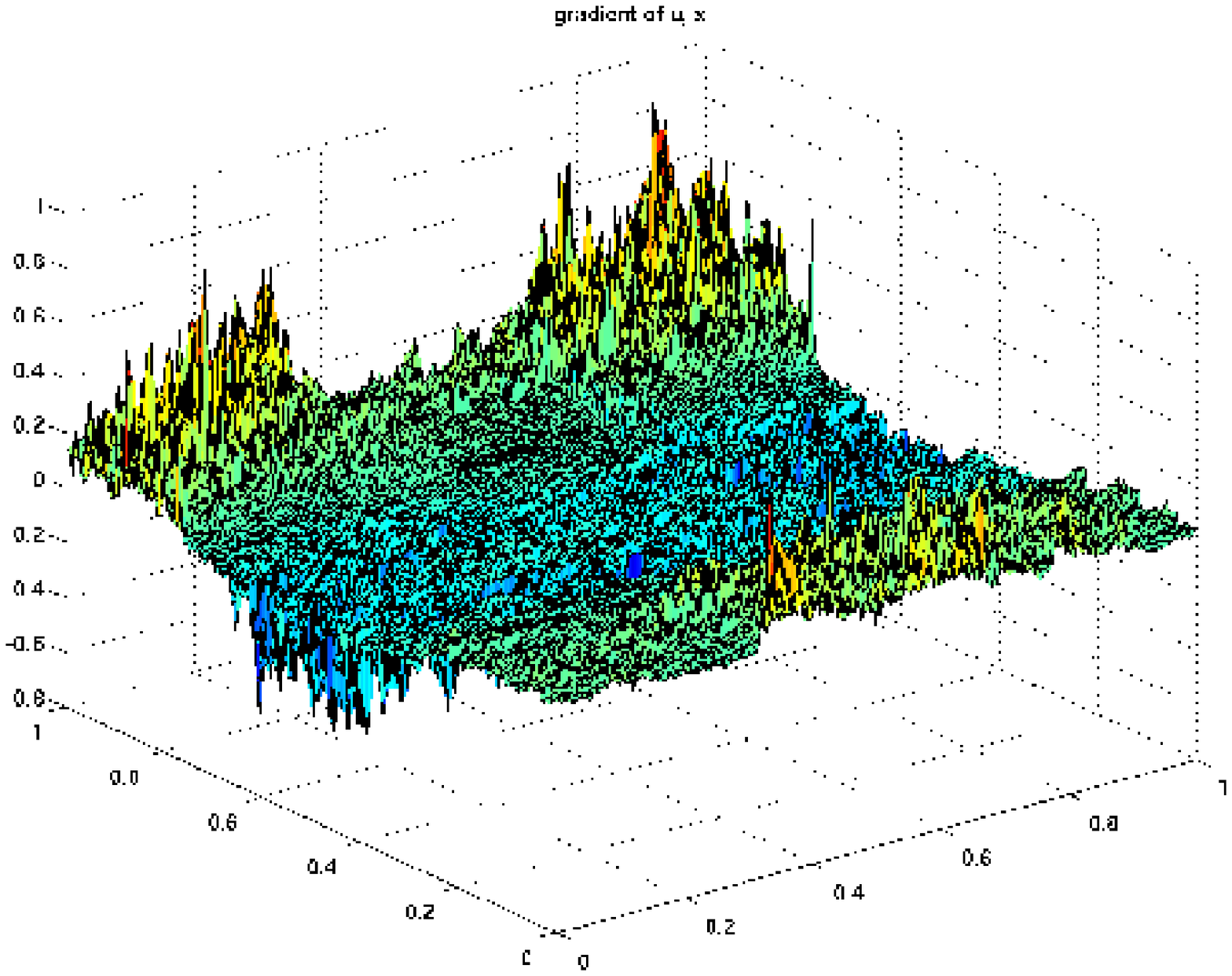}}
     \goodgap
    \subfigure[$(\nabla F)^{-1}\nabla u$.\label{figafrwwrr}]
          {\includegraphics[width=0.4\textwidth,height= 0.3\textwidth]{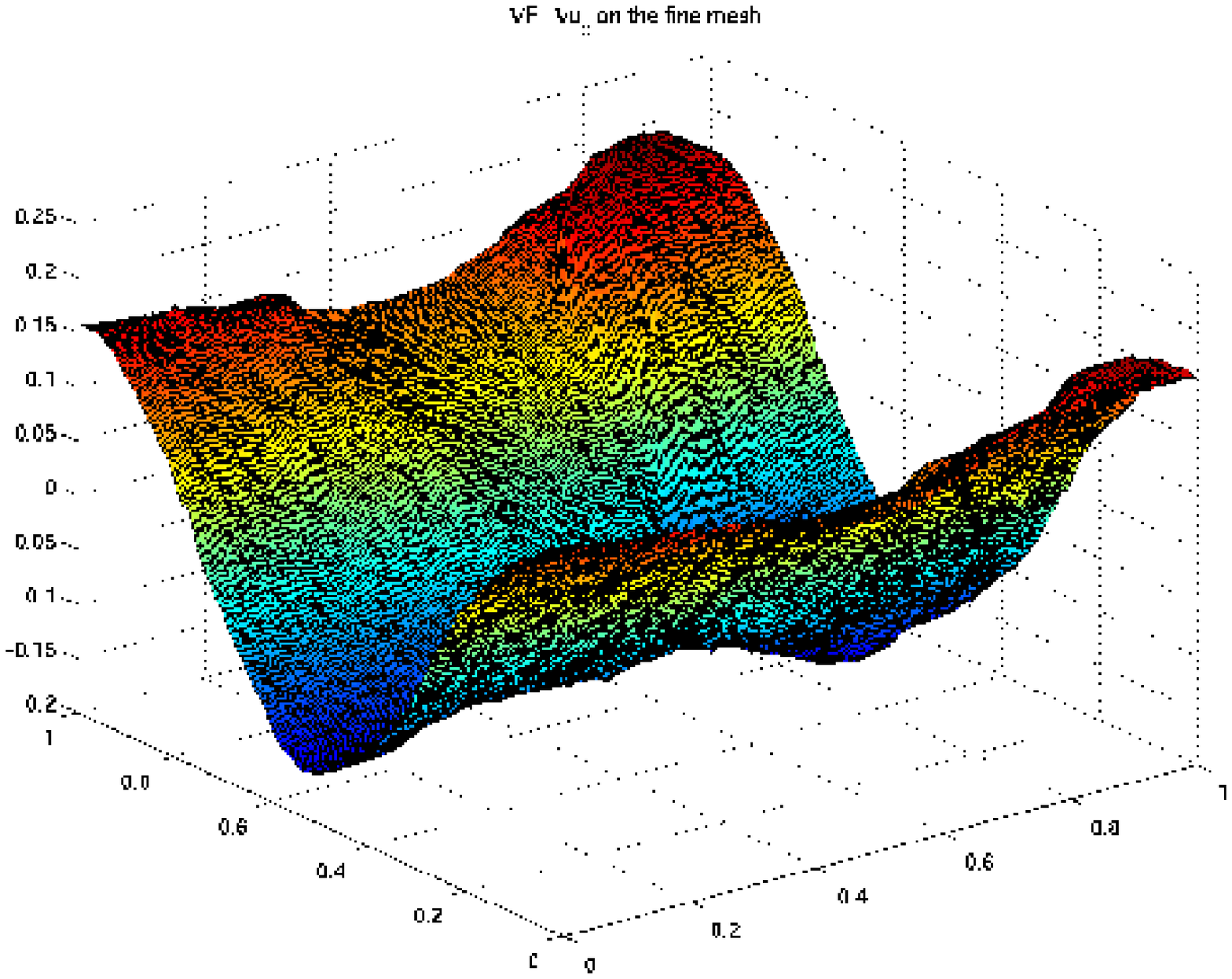}}
\\[-10pt]
    \caption{Change of metric on the torus.}
    \label{figyaqAfwwirst}
\end{center}
\end{figure}

\begin{Remark}
We will use the notation $\nabla_F u:= (\nabla F)^{-1}\nabla u$. In
dimension two, it is known (\cite{MR1892102}, \cite{MR2020365},
\cite{MR2001070}) that the determinant of $\nabla F$ is strictly
positive almost everywhere and the object $\nabla_F u$ is well
defined. In dimension three and higher $\nabla_F u$ is well defined
when $\sigma$ is stable.
\end{Remark}
This phenomenon can be observed numerically. In figures
\ref{particulsslll} and \ref{partwial} $a$ is given by a product of
random functions oscillating over a continuum of scales. The entries
of the matrix $\nabla F$ are in $L^{p}$ (figure \ref{figafrwrr}),
the entries of the gradient of $u$ in the Euclidean metric are in
$L^{p}$ (figures \ref{figfffrrr} and \ref{partiawwl}) yet $(\nabla
F)^{-1}\nabla u$ is H\"{o}lder continuous (figures
\ref{figdjdueAfirstB} and \ref{figafrwwrr}).


 Let
us now introduce the compensation phenomenon in dimension $n\geq 3$.
We call $\beta_\sigma$ the Cordes parameter associated to $\sigma$
defined by
\begin{equation}\label{sshgdswsd7641}
\beta_{\sigma}:=\esssup_{x\in \Omega}\Big(
n-\frac{\big(\Tr[\sigma(x)]\big)^2}{\Tr[{^t\sigma(x)\sigma(x)}]}\Big).
\end{equation}
Observe that since $\beta_{\sigma}$ is also given by
\begin{equation}\label{sshgdd7641}
\beta_{\sigma}=\esssup_{x\in \Omega}\Big( n-\frac{\big(\sum_{i=1}^n
\lambda_{i,\sigma(x)}\big)^2}{\sum_{i=1}^n
\lambda_{i,\sigma(x)}^2}\Big).
\end{equation}
where $(\lambda_{i,M})$ denotes the eigenvalues of $M$, it is a
measure of the anisotropy of $\sigma$.
\begin{Definition}
In dimension $n\geq 3$, we say that $\sigma$ is stable if and only
if, $\beta_{\sigma}<1$ and if $n\leq 4$ that there exist a constant
$\epsilon>0$ such that
$\big(\Tr(\sigma)\big)^{\frac{n}{2}-2-\epsilon}\in L^{1}(\Omega)$
\end{Definition}
\begin{Remark}
According to \cite{MR1892102} and \cite{MR2073507} in dimension
three and higher $\sigma$ can be unstable even if $a$ is smooth. We
refer to figure \ref{cap:computswsh} for an explicit example.
\end{Remark}
 Let
us write
\begin{equation}
\|v\|_{W^{2,p}_0 (\Omega)}:=\Big(\int_{\Omega} \big(\sum_{i,j=1}^n
|\partial_i\partial_j v|^2\big)^\frac{p}{2}\,dx\Big)^\frac{1}{p}
\end{equation}
\begin{Theorem}\label{th22}
Assume  that $\sigma$ is stable and $n\geq 3$. Then $F$ is an
automorphism on $\Omega$. Moreover there exist constants $p>2$ and
$C>0$ such that $u \circ F^{-1} \in W^{2,p}_0 (\Omega)$ and
\begin{equation}
\big\|u \circ F^{-1}\big\|_{W^{2,p}_0 (\Omega)} \leq C
\|g\|_{L^\infty(\Omega)}.
\end{equation}
\end{Theorem}
\begin{Remark}
The constant $p$ depends on $n,\Omega$ and $\beta_{\sigma}$. The
constant $C$ depends on the constants above,  $\lambda_{\min}(a)$
and if $n\leq 4$  on
$\big\|\big(\Tr(\sigma)\big)^{\frac{n}{2}-2-\epsilon}\big\|_{L^{1}(\Omega)}$.
\end{Remark}
In the following theorem we do not need  to assume $\Omega$ to be
convex.
\begin{Theorem}\label{ksjhgfts721}
Assume $n\geq 2$ and $(\Tr(\sigma))^{-1}\in L^\infty(\Omega)$. Let
$p>2$. There exist a constant $C^*=C^*(n,\partial \Omega)>0$ such
that if $\beta_{\sigma}<C^*$ then there exists a real number
$\gamma>0$ depending only on $n,\Omega$ and $p$ such that
\begin{equation}\label{hdhsxdzssssgxc7}
\begin{split}
\big\|(\nabla F)^{-1}\nabla u\big\|_{C^{\gamma}(\Omega)}^2\,dt \leq
C  \| g\|_{L^p(\Omega)}^2.
\end{split}
\end{equation}
\end{Theorem}
\begin{Remark}
The constant $C$ in \eref{hdhsxdzssssgxc7} depends on $n$, $\gamma$,
$\Omega$, $C^*$, $\lambda_{\min}(a)$, $\|a\|_{L^\infty(\Omega)}$,
$\mu_\sigma$ and
$\big\|(\Tr(\sigma))^{-1}\big\|_{L^\infty(\Omega_T)}$.
\end{Remark}

\subsection{Dimensionality reduction}\label{gfssg5}
Observe that \eref{ghjh52} is a priori an infinite dimensional
problem since $a$ and $g$ can be irregular at all scales. Yet
according to theorem \ref{th2} and \ref{th22}, whatever the choice
of $g$, at small scales,  solutions to \eref{ghjh52} are correlated
to $F$ which lives in a functional space of dimension $n$. More
precisely we will propose a rigorous justification of a variation of
the multi-scale finite element method\footnote{Let us recall that
the multi-scale finite element method is inspired from Tartar's
oscillating test functions \cite{MR557520}} introduced by Hou and Wu
\cite{MR1455261} in its form refined by Allaire and Brizzi
\cite{AlBr04} in situations where the medium is not assumed to be
periodic or ergodic (these methods are already rigorously justified
when the medium is periodic \cite{MR1455261}, \cite{AlBr04}).

Let $\T_h$ be a coarse conformal mesh on $\Omega$ composed of
$n$-simplices(triangles in dimension two and tetrahedra in dimension
three). Here $h$ is the usual resolution of the mesh defined as the
maximal length of the edges of the tessellation. Let us call
$\gamma(\T_h)$ the maximum over of the $n$-simplices $K$ of $\T_h$
of the ratio between the radius of the smallest ball containing $K$
and the largest ball inscribed in $K$. We assume $\gamma(\T_h)$  to
be uniformly bounded in $h$.

We write $V_h\subset H^1(\Omega)$  the set of piecewise linear
functions on the coarse mesh vanishing at the boundary of the
tessellation. We write $\No_h$ the set of interior nodes of the
tessellation and $\varphi_i$ ($i\in \No_h$) the usual nodal basis
function of $V_h$ satisfying
\begin{equation}
\varphi_i(y_j)=\delta_{ij}.
\end{equation}
We consider the elements $(\psi_i)_{i\in \No_h}$ defined by
\begin{equation}\label{gsfsgfd41}
\psi_i:=\varphi_i\circ F(x).
\end{equation}

\begin{figure}[htbp]
  \begin{center}
    \subfigure[$\varphi_i$
\label{partiawwwqsaswwl}]
          {\includegraphics[width=0.4\textwidth,height= 0.3\textwidth]{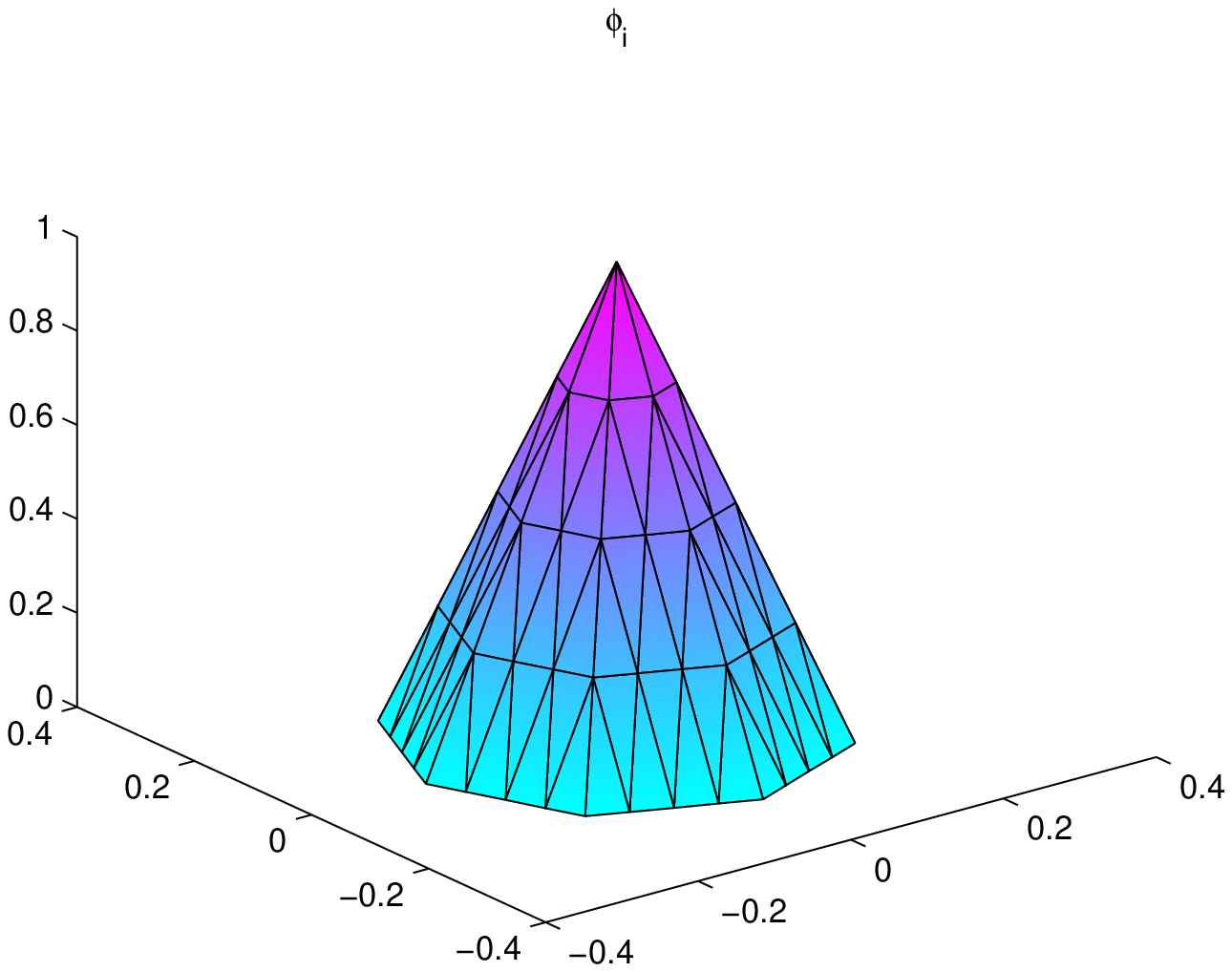}}
     \goodgap
    \subfigure[$\psi_i=\varphi_i \circ F$\label{figssafsarwwrr}]
          {\includegraphics[width=0.4\textwidth,height= 0.3\textwidth]{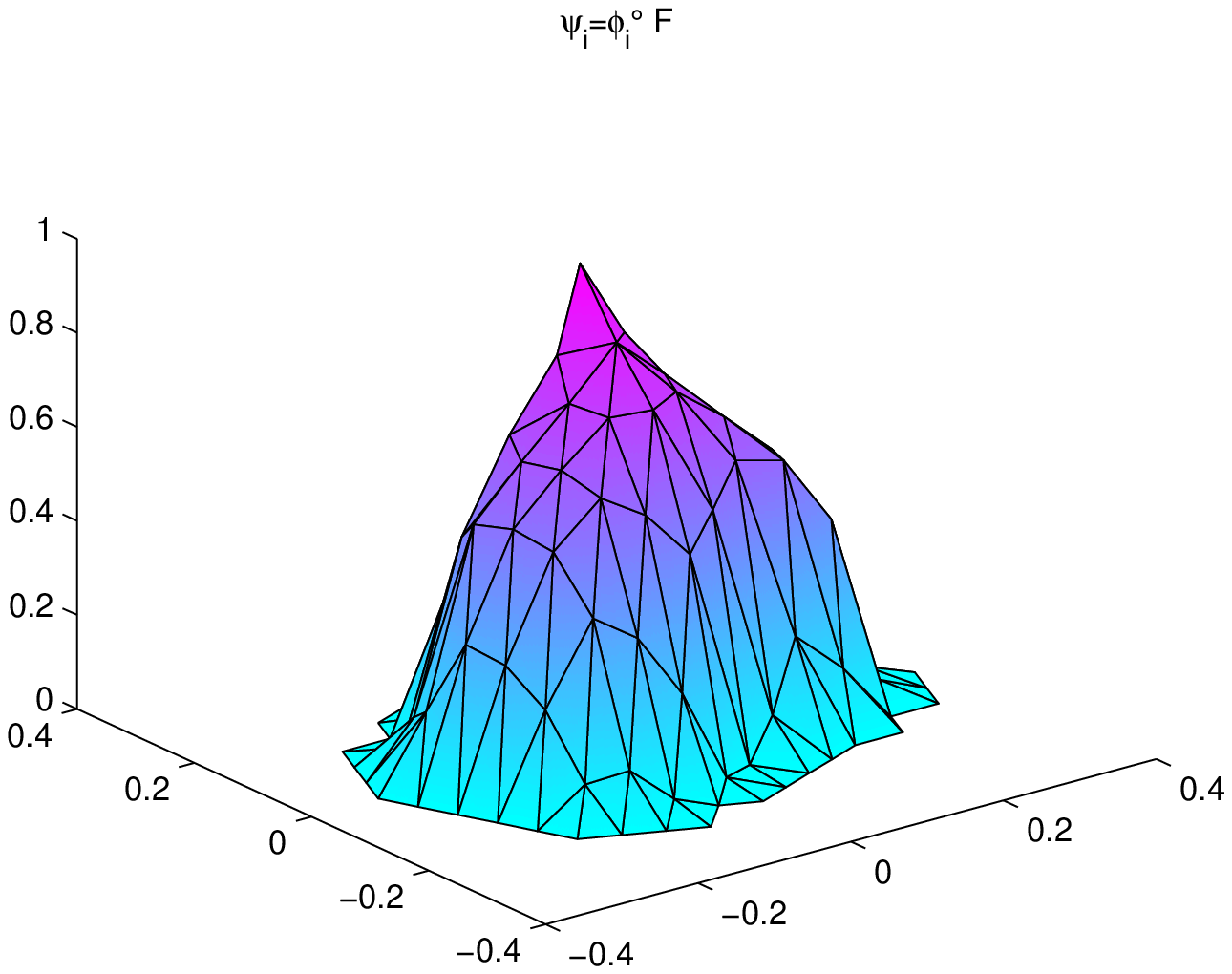}}
\\[-10pt]
    \caption{The Galerkin elements}
    \label{figyaqAfwwewirst}
\end{center}
\end{figure}
Let us write $u_h$ the solution of the Galerkin scheme associated to
\eref{dgdgfsghsf62} based on the elements $(\psi_i)_{i\in \No_h}$.
Observe that the number of elements is on the order of $h^{-n}$ and
we have the following theorem

\begin{Theorem}\label{sfsfssaaadd653}
Assume that $\sigma$ is stable and $n=2$. Then there exist constants
$\alpha,C>0$ such that
\begin{equation}\label{sddsdsaxssszsjsdsdhdgdgz7112}
\begin{split}
\|u-u_h\|_{H^1}\leq C  h^\alpha \|g\|_{L^\infty(\Omega)}.
\end{split}
\end{equation}
\end{Theorem}
\begin{Remark}
The constant $\alpha$ depends only on $n, \Omega$ and $\mu_\sigma$.
The constant $C$ depends on the objects mentioned above plus
$\lambda_{\min}(a),\lambda_{\max}(a)$, $\gamma(\T_h)$ and
$\big\|\big(\Tr(\sigma)\big)^{-1-\epsilon}\big\|_{L^{1}(\Omega)}$.
\end{Remark}
\begin{Remark}
Theorem \ref{sfsfssaaadd653} is also valid with $\alpha=1$ as in
theorem \ref{sfsfsswz653}. The only difference between these two
theorems lies in the constant $C$. In the proof theorem
\ref{sfsfssaaadd653} we use the property $u\circ F^{-1}\in
C^{1,\alpha}(\Omega)$ and in the proof of theorem \ref{sfsfsswz653}
we use the property $u\circ F^{-1}\in W^{2,2}(\Omega)$.
\end{Remark}
\begin{Remark}
Let us recall that $u_h$ is an element of the space $X_h$ spanned by
$(\psi_i)_{i\in \No_h}$ obtained as the solution of the linear
problem
\begin{equation}
a[\psi_i, u_h]=(\psi_i,g)_{L^2(\Omega)}.
\end{equation}
Where $a[,]$ is the bilinear form on $H^1_0(\Omega)$ defined by
\begin{equation}
a[v,w]:=\int_{\Omega}\nabla v a \nabla w.
\end{equation}
It follows from theorem \ref{sfsfssaaadd653} that solutions to
\eref{ghjh52} live in the $H^1$-norm neighborhood of a low
dimensional space.
\end{Remark}
\begin{Remark}
The proof of theorem \ref{sfsfssaaadd653} is done for the exact
 function $\psi_i$ and not for its discrete version on a fine mesh.
If $a$ is regular at a given small scale $h_0$ then it is easy to
check that theorem \ref{sfsfssaaadd653} remains valid as long as the
edges of the fine mesh are smaller than $h_0$. A more intriguing
case is when $a$ is discrete and discontinuous on a fine mesh.
Numerical experiments  show that theorems such as
\ref{sfsfssaaadd653} and \ref{th2} remain valid but to justify them
for the discrete version of the harmonic coordinates $F$ and
elements $\psi_i$ one would have to adapt our theorems to the
discrete setting. In order to remain concise we have chosen to not
include that adaptation in this paper.
\end{Remark}

\begin{Remark}
We keep the composition rule used in \cite{AlBr04}. The only
difference between the elements \eref{gsfsgfd41} and the ones
proposed by Hou, Wu, Allaire and Brizzi lies in the fact that we use
the global solution to \eref{dgdgfsghsf62} and not a local one
computed on each triangle of the coarse mesh through an
over-sampling technique. We refer to remark \ref{hgshsh512} for
further comments.
\end{Remark}
\begin{Remark}\label{gsfgsftwt541}
Write $S$ the stiffness matrix $a[\psi_i,\psi_j]$. $S^{-1}$ is in
general dense (characterized by $N^2$ entries where $N$ is the
number of nodes of the mesh). Yet surprisingly by combining theorem
\ref{sfsfssaaadd653} with theorem 5.4 of \cite{MR1993936} one can
obtain that $S^{-1}$ can be approximated (in $L^2$ norm) with
hierarchical matrix $M_H$ such that the matrix vector product by
$M_H$ requires only $O(N (\ln N)^{n+3})$ operations.
\end{Remark}

In dimension $n\geq 3$ we have the following estimate.
\begin{Theorem}\label{sfsfsswz653}
Assume that $\sigma$ is stable,   $n\geq 3$ and\\
$\big\|\big(\Tr(\sigma) \big)^{-1}\big\|_{L^\infty(\Omega)}<\infty$.
Then there exist constants $\alpha,C>0$ such that
\begin{equation}\label{sddsdssessjsdsdhdgdgz7112}
\begin{split}
\|u-u_h\|_{H^1}\leq C  h \|g\|_{L^\infty(\Omega)}.
\end{split}
\end{equation}
\end{Theorem}
\begin{Remark}
 The constant $C$ depends on $n$, $\gamma(\T_h)$, $\Omega$, $\beta_{\sigma}$,
 $\lambda_{\max}(a)$, $\lambda_{\min}(a)$ and $\big\|\big(\Tr(\sigma) \big)^{-1}\big\|_{L^\infty(\Omega)}$.
\end{Remark}

\subsection{Galerkin with localized elements}\label{hgdjhdg541}
For the clarity of the paper we will restrict ourselves from now on
to dimension two, the generalization of the statements to higher
dimensions is conditioned on the stability of $\sigma$ (and the
application of theorem \ref{ksjhgfts721}).

\begin{figure}[htbp]
  \begin{center}
\includegraphics[width=0.6\textwidth,height= 0.5\textwidth]{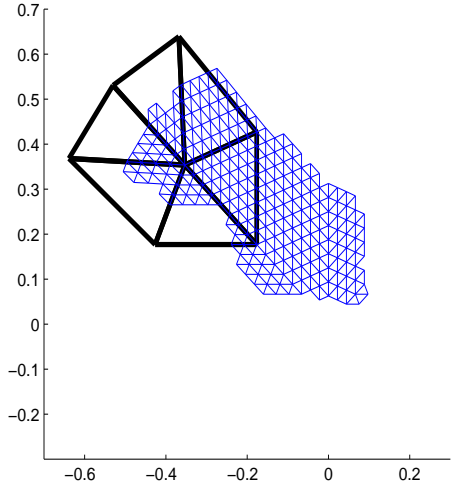}
    \caption{Support of the elements $\varphi_i$ and $\psi_i$.}
    \label{figyaqAfsswwirst}
\end{center}
\end{figure}
The elements \eref{gsfsgfd41} can be highly distorted and non local
(figure \ref{figyaqAfsswwirst}) since
\begin{equation}
\supp(\psi_i):=F^{-1}\big(\supp (\varphi_i)\big).
\end{equation}
It follows that the elements $\psi_i$ are piecewise linear on a fine
mesh different from the one on which $a$ is defined and $F$ has been
computed.  Is it possible to avoid that difficulty by solving
\eref{ghjh52} on a coarse mesh with localized elements? The answer
is yes but the price to pay for the localization will be the
discontinuity of the elements and the fact that the accuracy of the
method will depend on a weak aspect ratio of the triangles of the
tessellation in the metric induced by $F$. More precisely consider a
triangle $K$ of the tessellation, call $a,b,c$ the nodes of $K$ and
$\theta$ the interior angle of the triangle
$\big(F(a),F(b),F(c)\big)$ which is the closest to $\pi/2$.
 We call weak aspect ratio of the triangle $K$ in the
metric induced by $F$ the quantity
\begin{equation}
\eta^F_{\min}(K):=\frac{1}{\sin(\theta)}.
\end{equation}
So $\eta^F_{\min}(K)$ is large if the triangle
$\big(F(a),F(b),F(c)\big)$ is flat (all its interior angles are
close to $0$ or $\pi$). We define
\begin{equation}
\eta^*_{\min}:=\sup_{K\in \T_h}\eta^F_{\min}(K).
\end{equation}
Let us recall that  although the coefficients of the PDE
\eref{ghjh52} are irregular it is well known \cite{MR0352696} that
$F$ is H\"{o}lder continuous. Thus it makes sense to look at the
value of $F$ at a specific point. Now let $v$ be a function defined
on the nodes of the triangle $K\in \T_h$, let us write $a,b,c$ the
nodes of that triangle. It is usual to look at the coarse gradient
of $v$ evaluated at the nodes of the triangle $K$, i.e. the vector
defined by
\begin{equation}\label{54s21q}
\nabla v(K):=\left(\begin{array}{c}
b-a \\
c-a \end{array} \right)^{-1} \left( \begin{array}{c}
v(b)-v(a) \\
v(c)-v(a) \end{array} \right).
\end{equation}
If $\eta^F_{\min}(K)<\infty$ then the following object called the
gradient of $v$ evaluated on the coarse mesh with respect to the
metric induced by $F$ is well defined.
\begin{equation}\label{543swwwq21q}
\nabla_F v(K):=\left(\begin{array}{c}
F(b)-F(a) \\
F(c)-F(a) \end{array} \right)^{-1} \left( \begin{array}{c}
v(b)-v(a) \\
v(c)-v(a) \end{array} \right).
\end{equation}
\begin{Definition}\label{shdggd61}
We say that the tessellation $\T_h$ is not unadapted to $F$ if and
only if the determinant of $\nabla F(K)$ is strictly positive for
all $K\in \T_h$.
\end{Definition}
\begin{Remark}
Observe that if the tessellation $\T_h$ is not unadapted to $F$ then
$\eta^*_{\min}(K)<\infty$, the definition \ref{shdggd61} contains
the additional condition that there is no inversion in the images of
the triangles of $\T_h$ by $F$.
\end{Remark}
Now consider the nodal elements $(\xi_i)_{i\in \No_h}$, defined by
\begin{equation}\label{jshshsg51}
\begin{cases}
\xi_i(x_j)=\delta_{ij} \\
\nabla_F \xi(x)=\text{constant within each}\quad K\in \T_h.
\end{cases}
\end{equation}
\begin{figure}[htbp]
  \begin{center}
\includegraphics[width=0.7\textwidth,height= 0.6\textwidth]{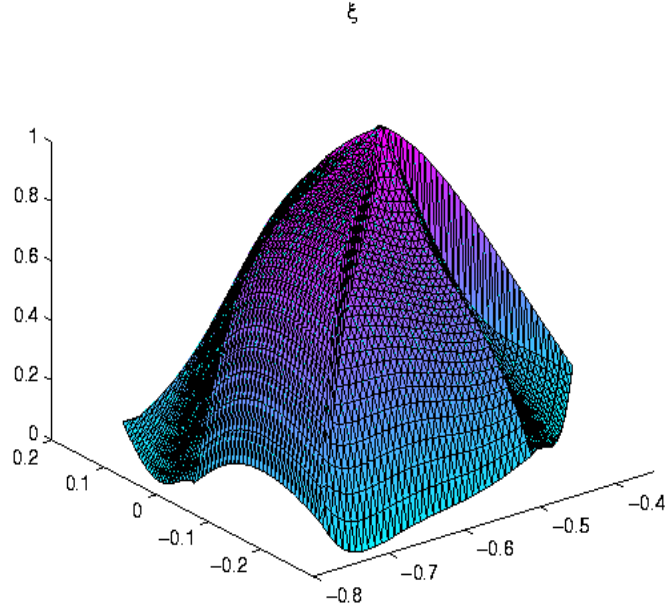}
    \caption{Localized Galerkin elements $\xi_i$}
    \label{figyaqAfwwisawrst}
\end{center}
\end{figure}
If the mesh is not unadapted to $F$ then the elements (figure
\ref{figyaqAfwwisawrst}) \eref{jshshsg51} are well defined and given
by
\begin{equation}
\begin{cases}
\xi_i(x)=1+\big(F(x)-F(x_i)\big) \nabla_F \varphi_i(K)\quad
\text{if}\quad i\sim
K\quad \text{and}\quad x\in K \\
\xi_i(0)=0 \quad \text{in other cases}.
 \end{cases}
\end{equation}
Where the notation $i\sim K$ means that $i$ is a node of $K$.
Observe that the elements $\xi_i$ are discontinuous at the
boundaries of the triangles of the coarse mesh however they are
easier to implement since they are localized in these triangles.
Write $Z_h$ the vector space spanned by the functions $\xi_i$. For
$K\in \T_h$ we write $a_K$ the bilinear form on $H^1(K)$ defined by
\begin{equation}
a_K[v, w]:=\int_{K} {^t\nabla v} a \nabla w.
\end{equation}
We will write $H^1(\T_h)$ the space of functions $v\in L^2(\Omega)$
such that the restriction of $v$ to each triangle $K$ belongs to
$H^1(K)$. We will write for $v,w \in H^1(\T_h)$
\begin{equation}
a^*[v,w]:=\sum_{K\in \T_h }a_K[v, w].
\end{equation}
The localized finite element method can be formulated in the
following way: look for $u^f \in Z_h$ such that for all $i\in
\No_h$,
\begin{equation}\label{gdhgdfctesstfc5}
a^*[\xi_i, u^f ]=(\xi_i,g)_{L^2(\Omega)}.
\end{equation}
We have the following estimate
\begin{Theorem}\label{sfsfssasadd653}
Assume that $\sigma$ is stable and that the mesh is not unadapted to
$F$. Then there exists a constant $\alpha>0$ such that
\begin{equation}\label{sddsdsaxsssszsjsdsdhdgdgz7112}
\begin{split}
\big(a^*[u-u^f]\big)^\frac{1}{2} \leq C \eta^*_{\min} h^\alpha
\|g\|_{L^\infty(\Omega)}.
\end{split}
\end{equation}
\end{Theorem}
\begin{Remark}
For a bilinear from $B[,]$ we write  $B[v]:=B[v,v]$.
\end{Remark}
\begin{Remark}
The constant $\alpha$ depends only on $n, \Omega, \epsilon$ and
$\mu_\sigma$. The constant $C$ depends on the  objects mentioned
above plus
$\big\|\big(\Tr(\sigma)\big)^{-1-\epsilon}\big\|_{L^{1}(\Omega)}$.
\end{Remark}
\begin{Remark}
The bilinear operator $a^*[,]$ on $Z_h$ is characterized by a
constant matrix within each triangle $K\in \T_h$ equal to
\begin{equation}
{^t\big(}\nabla F(K)\big)^{-1}\big< {^t\nabla F}a\nabla F\big>_K
\big(\nabla F(K)\big)
\end{equation}
where $<v>_K$ means the average of $v$ over $K$ with respect to the
Lebesgue measure ($<v>_K:=\frac{1}{\Vol(K)}\int_{K}v(x)\,dx$,
$\Vol(K)$ being the Lebesgue measure of volume of $K$).
\end{Remark}
The error bound given in theorem \ref{sfsfssasadd653} is given in
the norm induced by $a^*[.]$. We would like to obtain an error bound
with respect to the usual $H^1$ norm. Observe that $u^f$ is
discontinuous at the boundaries of the triangles of the coarse mesh
so we have to find an accurate way to interpolate $u^f$ in the whole
space using its values at the nodes of the coarse mesh.
 Let us write
$F(\No_h)$ the image of the nodes of $\T_h$ by $F$. Let us write
$\T^F$ the triangulation of  $F(\No_h)$. Let $\varphi_i^F$ be the
standard piecewise linear nodal basis of $\T^F$.  Let us write
$\J_h$ the interpolation operator from the space of functions
defined on the nodes of $\T_h$ into $H^1(\Omega)$ defined by
\begin{equation}\label{gdfd541}
\J_h v(x):= \sum_{i\in \No_h} v(x_i) \varphi_i^F \circ F(x).
\end{equation}
Observe that for $i\in \No_h$, $v(x_i)=\J_h v(x_i)$. We have the
following estimate
\begin{Theorem}\label{sfsasd653}
Assume that $\sigma$ is stable and that the mesh is not unadapted to
$F$. Then there exist constants $\alpha, C_f>0$ such that
\begin{equation}\label{sqwsszsjsddgz7112}
\begin{split}
\| u-\J_h u^f\|_{H^1(\Omega)} \leq C_f h^\alpha
\|g\|_{L^\infty(\Omega)}.
\end{split}
\end{equation}
\end{Theorem}
\begin{Remark}
The constant $\alpha$ depends only on $n, \Omega$ and $\mu_\sigma$.
The constant $C_f$ can be written
\begin{equation}
C_f:= C \eta^*_{\min}
\Big(\min\big(\eta^*_{\max}\eta_{\max}^3,\nu^*\big)\Big)^\frac{1}{2}
\end{equation}
where $C$ depends on the objects mentioned above plus
$\lambda_{\min}(a), \lambda_{\max}(a)$ and
$\big\|\big(\Tr(\sigma)\big)^{-1-\epsilon}\big\|_{L^{1}(\Omega)}$.
$\eta_{\max}$ is  defined by $\frac{1}{\sin \theta}$ where $\theta$
the interior angle of the triangles of $\T_h$ closest to $0$ or
$\pi$. $\eta_{\max}^*$ is  defined by $\frac{1}{\sin \gamma}$ where
$\gamma$ the interior angle of the triangles of $\T^F$ closest to
$0$ or $\pi$. Moreover
\begin{equation}
\nu^*:=\sup_{K\in \T_h} \frac{\Vol(K^F)}{\Vol(F(K))}
\end{equation}
where $K^F$ is the triangle whose nodes are the images of the nodes
of $K$ by $F$.
\end{Remark}

\subsection{Numerical homogenization from the information point of
view.}\label{gfgxfs541} The Galerkin scheme described in
\ref{gfssg5} and \ref{hgdjhdg541} are based on elements containing
the whole fine scale structure of $F$. This represents too much
information. We can wonder what minimal quantity of information
should be kept from the scales in order to up-scale \eref{ghjh52}?
We would like to keep an accurate version of \eref{ghjh52} with
minimal computer memory. This point touches the compression issue.
Images can be compressed. Can the same thing be done with operators?

This question has received an answer within the context of the fast
multiplication of vectors with fully populated special matrices
arising in various applications \cite{FeSt04}. Let us recall the
fast multipole method and the hierarchical multipole method designed
by L. Greengard and V. Rokhlin \cite{MR918448}. Wavelet based
methods have been designed by G. Beylkin, R. Coifman and V. Rokhlin
\cite{MR1201316, MR1276529, MR1085827}. The concept of Hierarchical
matrices has been developed by W. Hackbusch et al. \cite{MR1981528}.
More precisely we refer to \cite{MR1993936, BeCh05, Beb05, Bebe05,
Beben05}. The Hierarchical matrix method is based on a compression
of the inverse of the stiffness matrix (see remark
\ref{gsfgsftwt541}). Here we consider compression from the point of
view of numerical homogenization. We  look at the operator
\eref{ghjh52} as a bilinear form on $H^1_0(\Omega)$ and we will use
$V_h$ as space of test functions to {\it zoom at} the operator
associated to $a$ at a given arbitrary resolution.
\begin{equation}
a:
\begin{cases}
H^1_0(\Omega)\times  H^1_0(\Omega) &\rightarrow \R \\
(v,w) &\rightarrow \int_{\Omega}{^t\nabla v a \nabla w}.
\end{cases}
\end{equation}
The up-scaled or compressed operator, written $\U_h a$ will
naturally be a bilinear form on the space of piecewise linear
functions on the coarse mesh with Dirichlet boundary condition.
\begin{equation}
\U_h a:
\begin{cases}
V_h \times  V_h &\rightarrow \R \\
(v,w) &\rightarrow \U_ha[v,w].
\end{cases}
\end{equation}
The question is how to choose $\U_h a$? To answer that question we
can integrate \eref{ghjh52} against a test function $\phi$ in $V_h$,
then we obtain that
\begin{equation}
\int_{\Omega}\nabla \phi a \nabla u=\int_{\Omega}\phi g.
\end{equation}
We will use the test function $\phi$ to "look at" the operator
\eref{ghjh52} at the given resolution $h$. We can decompose the
first term in the integral above as a sum of integrals over the
triangles of the coarse mesh to obtain (we assume that $\sigma$ is
stable)
\begin{equation}\label{gfdfd541}
\int_{\Omega}\nabla \phi a \nabla u=\sum_{K\in \T_h}\int_{K}\nabla
\phi(x) a(x) \nabla F(x) \big(\nabla F(x)\big)^{-1}\nabla u(x)\,dx.
\end{equation}
Now $\nabla \phi$ is constant within each triangle $K\in \T_h$.
$\big(\nabla F(x)\big)^{-1}\nabla u(x)$ is H\"{o}lder continuous
thus almost a constant within each triangle $K$ and equal to the
gradient of $u$ evaluated on the coarse mesh with respect to the
metric induced by $F$, i.e. the following vector
\begin{equation}\label{54321q}
\nabla_F u(K):=\left(\begin{array}{c}
F(b)-F(a) \\
F(c)-F(a) \end{array} \right)^{-1} \left( \begin{array}{c}
u(b)-u(a) \\
u(c)-u(a) \end{array} \right).
\end{equation}
Where $a,b,c$ are the nodes of the triangle $K$. It follows that the
tensor $a\nabla F$ can be averaged over each triangle of the coarse
mesh and we will write $<a\nabla F>_K$ its average. In conclusion a
good candidate for the up-scaled operator $\U_h a$
 is the bilinear form given by the following formula: for $v,w\in V_h$
\begin{equation}\label{jhgeehg41}
\U_ha [v,w]:=\sum_{K\in \T_h} \int_{K} {^t \nabla v}\big< a \nabla
F\big>_K \big(\nabla F(K)\big)^{-1} \nabla w.
\end{equation}
Observe that the only information kept from the small scales in the
compressed operator \eref{jhgeehg41} are the  bulk quantities $\big<
a \nabla F\big>_K$ and the non averaged quantities $F(b)-F(a)$ where
$a$ and $b$ are nodes of the triangles of the coarse mesh. The
latter quantity can be interpreted as a deformation of the coarse
mesh induced by the small scales (or a new distance defining coarse
gradients). In the particular case where $a=M(\frac{x}{\epsilon})$
and $M$ is ergodic then as $\epsilon \downarrow 0$  $<a\nabla F>_K$
converges to the usual  effective conductivity obtained from
homogenization theory and $\nabla F(K)$ converges to the identity
matrix. It follows that the object \eref{jhgeehg41} recovers the
formulae obtained from homogenization theory when the medium is
ergodic and characterized by scale separation. Let us now show that
formula \eref{jhgeehg41} can be accurate beyond these assumptions.

\begin{figure}[httb]
\begin{center}
\includegraphics[%
  scale=0.3]{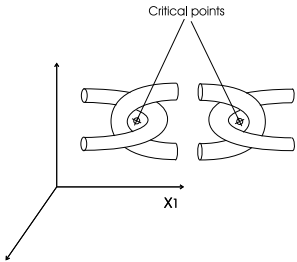}
\caption{\label{cap:computswsh}$a$ in dimension three}
\end{center}
\end{figure}

To estimate the compression accuracy we have to use the up-scaled
operator $\U_h a$ to obtain an approximation of the linear
interpolation of $u$ on the coarse mesh. We look for $u^m \in V_h$
such that for all $i\in \No_h$,
\begin{equation}\label{gdhgdfcssstetfc5}
\U_h a[\varphi_i, u^m ]=(\varphi_i,g)_{L^2(\Omega)}.
\end{equation}
The price to pay for the loss of information on the small scales is
the loss of ellipticity. This loss can be caused by two correlated
factors:
\begin{itemize}
\item The new metric can generate flat triangles.
\item The up-scaled operator can become singular.
\end{itemize}
The first factor is due to the localization of the scheme. The
second factor does not appear with Galerkin schemes. It is not
observed in dimension two but it can't be  avoided in dimension
higher or equal to three in the sense that the up-scaled operator
has no reason to remain elliptic and local. Indeed consider a box of
dimension three, and set in that box empty tubes of low boundary
conductivity as shown in figure \ref{cap:computswsh}. Set the left
side of the box to temperature $0^0c$ and the right side to
temperature $100^0c$. Then an inversion in the temperature profile
is produced around the critical points shown in figure
\ref{cap:computswsh} (see \cite{MR1892102} and \cite{MR2073507},
instead of increasing from left to right in these region temperature
decreases). Now as the operator is up-scaled, the information on the
geometry of the tubes is lost but the inversion phenomenon remains
in the loss of ellipticity and locality of the operator. We will
address  this issue further in a forthcoming paper.

\begin{figure}[htbp]
  \begin{center}
    \subfigure[Refined $F$.
\label{particulsslsll}.]
          {\includegraphics[width=0.3\textwidth,height= 0.3\textwidth]{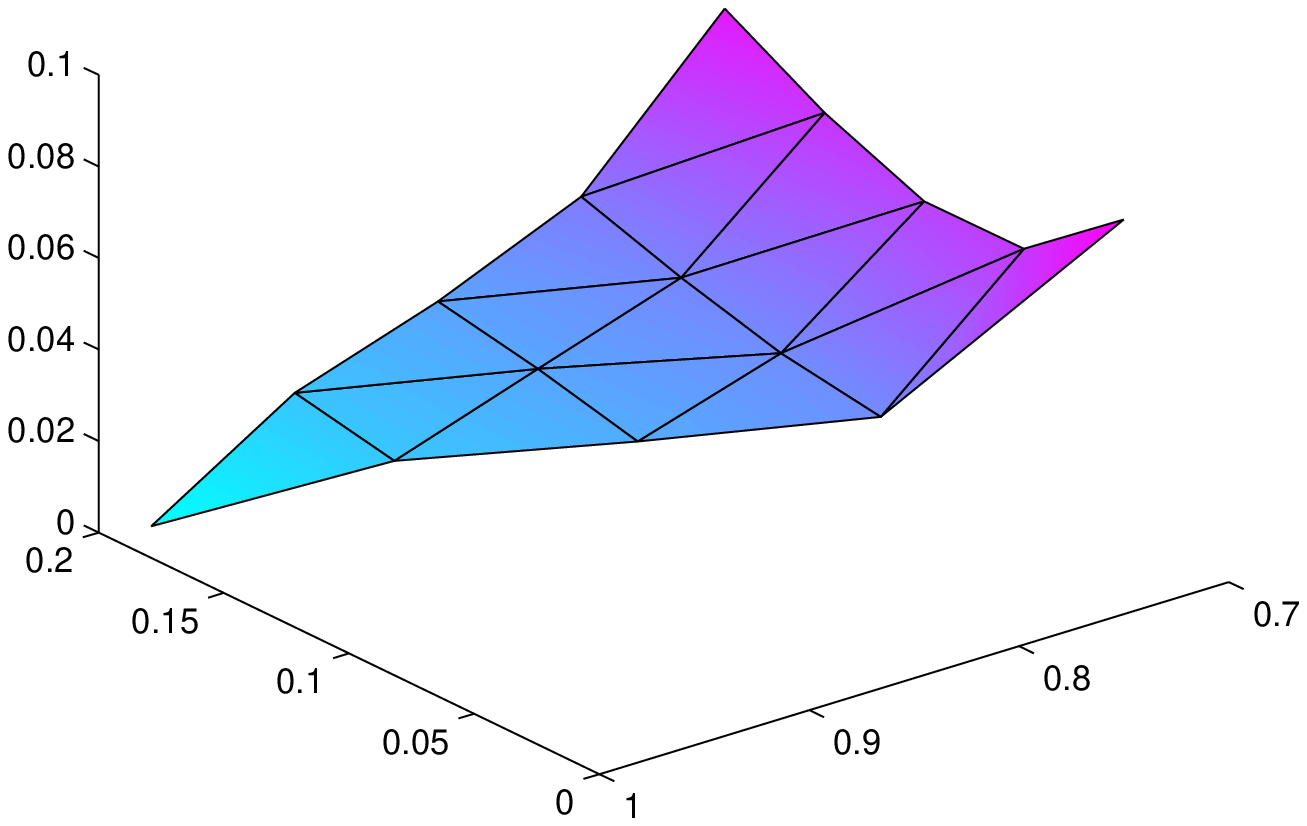}}
     \goodgap
    \subfigure[Refined $F$.\label{fisgfffrrr}]
          {\includegraphics[width=0.3\textwidth,height= 0.3\textwidth]{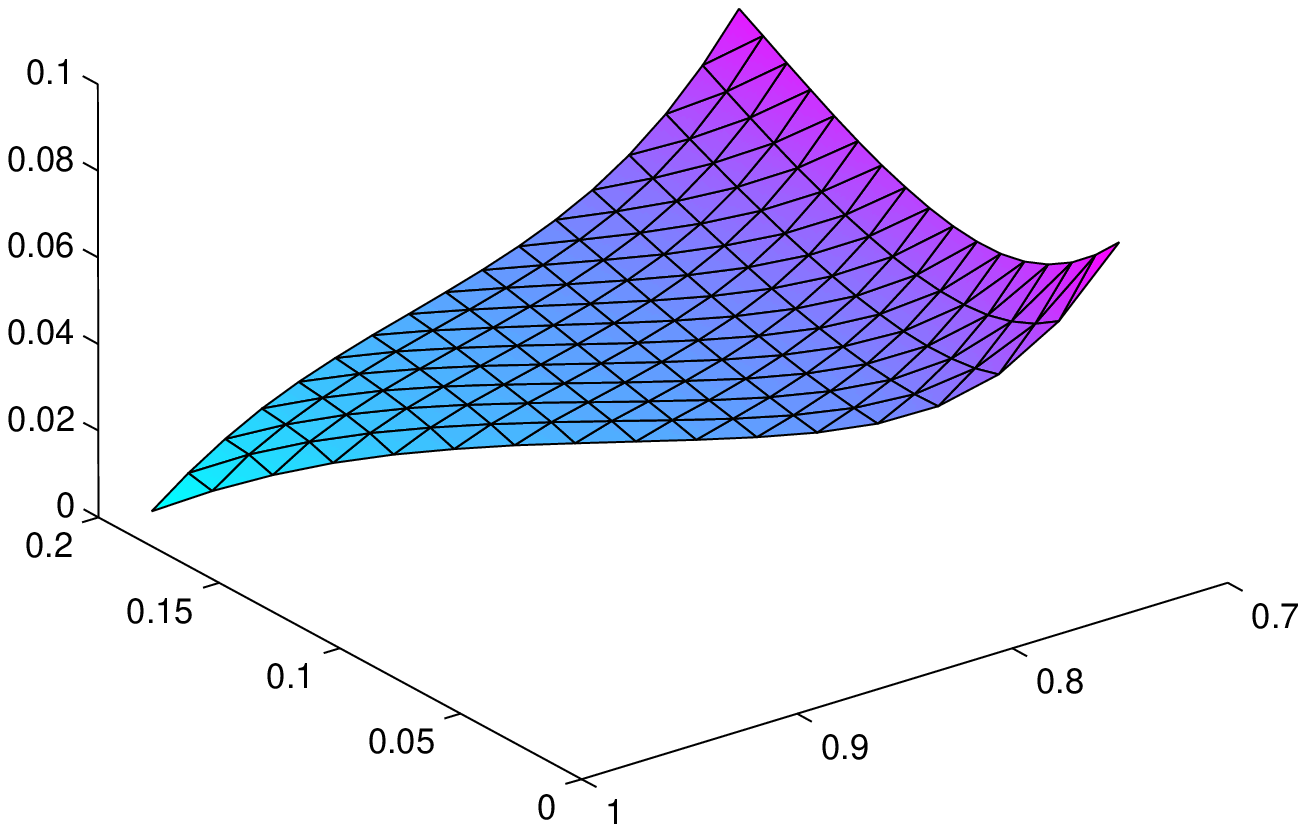}}
    \goodgap
   \subfigure[$\J_h u$.\label{figdsjdsueAfirstB}]
        {\includegraphics[width=0.3\textwidth,height= 0.3\textwidth]{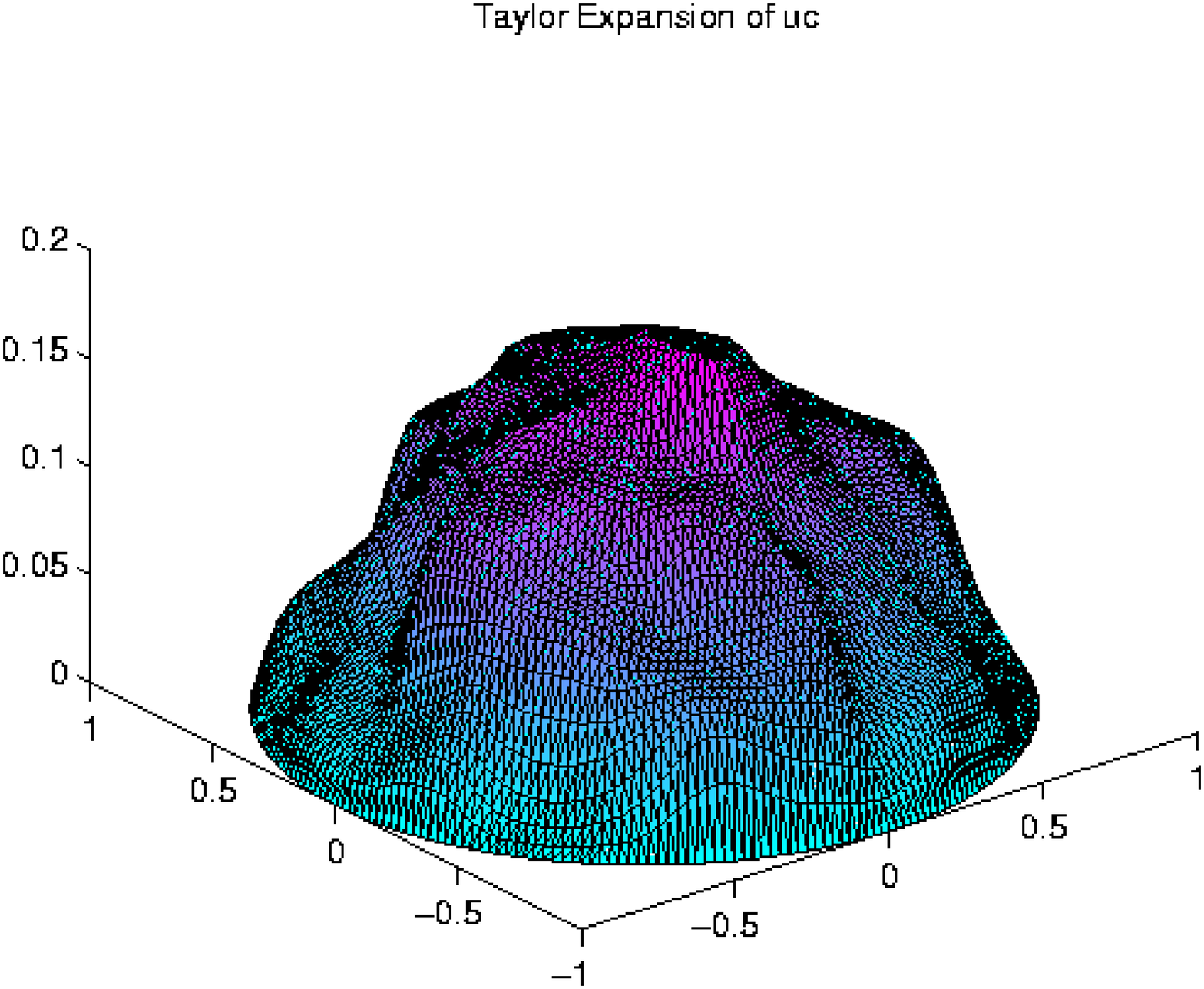}}\\[-10pt]
    \caption{Taylor expansion with respect to the new metric.}
    \label{figyyswassAfirst}
\end{center}
\end{figure}

Nevertheless it is possible to prove that once stability  is
achieved then the method is accurate (if $\sigma$ is stable). More
precisely, for a nodal function $v$ let us define the homogeneous
Dirichlet form on the graph induced by $\T_h$:
\begin{equation}
\ED_h[v]:=\sum_{i\sim j} |v_i-v_j|^2
\end{equation}
We write $i\sim j$ when those nodes share an edge on the coarse
mesh. Let us define  the following stability parameter of the scheme
\begin{equation}\label{gsgfd763}
\S^m:=\inf_{w \in V_h} \sup_{v\in V_h} \frac{\U_h
a[v,w]}{\big(\ED_h[v]\big)^\frac{1}{2}
\big(\ED_h[w]\big)^\frac{1}{2}}.
\end{equation}
Observe that $\S^m$ depends only on the up-scaled parameters so we
have a control on the stability.
\begin{Definition}
We say that the scheme is stable if and only if $\S^m>0$.
\end{Definition}
\begin{Remark}
Let us recall that for $v\in V_h$, $\ED_h[v]$ can be bounded from
below and above by the $L^2$-norm of the gradient of $v$. More
precisely
\begin{equation}\label{ddff351ssad32435}
\frac{1}{4 \eta_{\max}}  \ED_h[v] \leq \|\nabla v\|_{L^2(\Omega)}^2
\leq \eta_{\max} \ED_h[v].
\end{equation}
$\eta_{\max}=1/\sin(\theta)$ where $\theta$ is the closest interior
angle of the triangles of $\T_h$ to $0$ or $\pi$.
\end{Remark}
\begin{Remark}
In practice in dimension two the condition number of the scheme
associated to the up-scaled operator is as good as the one obtained
from a Galerkin scheme by solving a local cell problem.
\end{Remark}
Let us write $\I_h u$ the linear interpolation of $u$ over $\T_h$:
\begin{equation}
\I_h u:=\sum_{i\in \No_h} u(x_i) \varphi_i (x).
\end{equation}
 We have the following estimate
\begin{Theorem}\label{sfsfxadd653}
Assume that $\sigma$ and the scheme are stable and that the mesh is
not unadapted to $F$. Then there exist constants $\alpha,C_m>0$ such
that
\begin{equation}\label{sqxsssszsjsdsdhdgdgz7112}
\begin{split}
\|\I_h u-u^m\|_{H^1(\Omega)} \leq C_m h^\alpha
\|g\|_{L^\infty(\Omega)}.
\end{split}
\end{equation}
\end{Theorem}
\begin{Remark}\label{sgsff41}
The constant $\alpha$ depends only on $n, \Omega$ and $\mu_\sigma$.
The constant $C_m$ can be written
\begin{equation}
C_m:= C \frac{\eta^*_{\min} \eta_{\max}}{S^m}.
\end{equation}
where $C$ depends on the objects mentioned above plus
$\lambda_{\min}(a), \lambda_{\max}(a)$ and
$\big\|\big(\Tr(\sigma)\big)^{-1-\epsilon}\big\|_{L^{1}(\Omega)}$.
\end{Remark}
\begin{figure}[htbp]
  \begin{center}
    \subfigure[$u$
\label{particulssaslll}]
          {\includegraphics[width=0.3\textwidth,height= 0.3\textwidth]{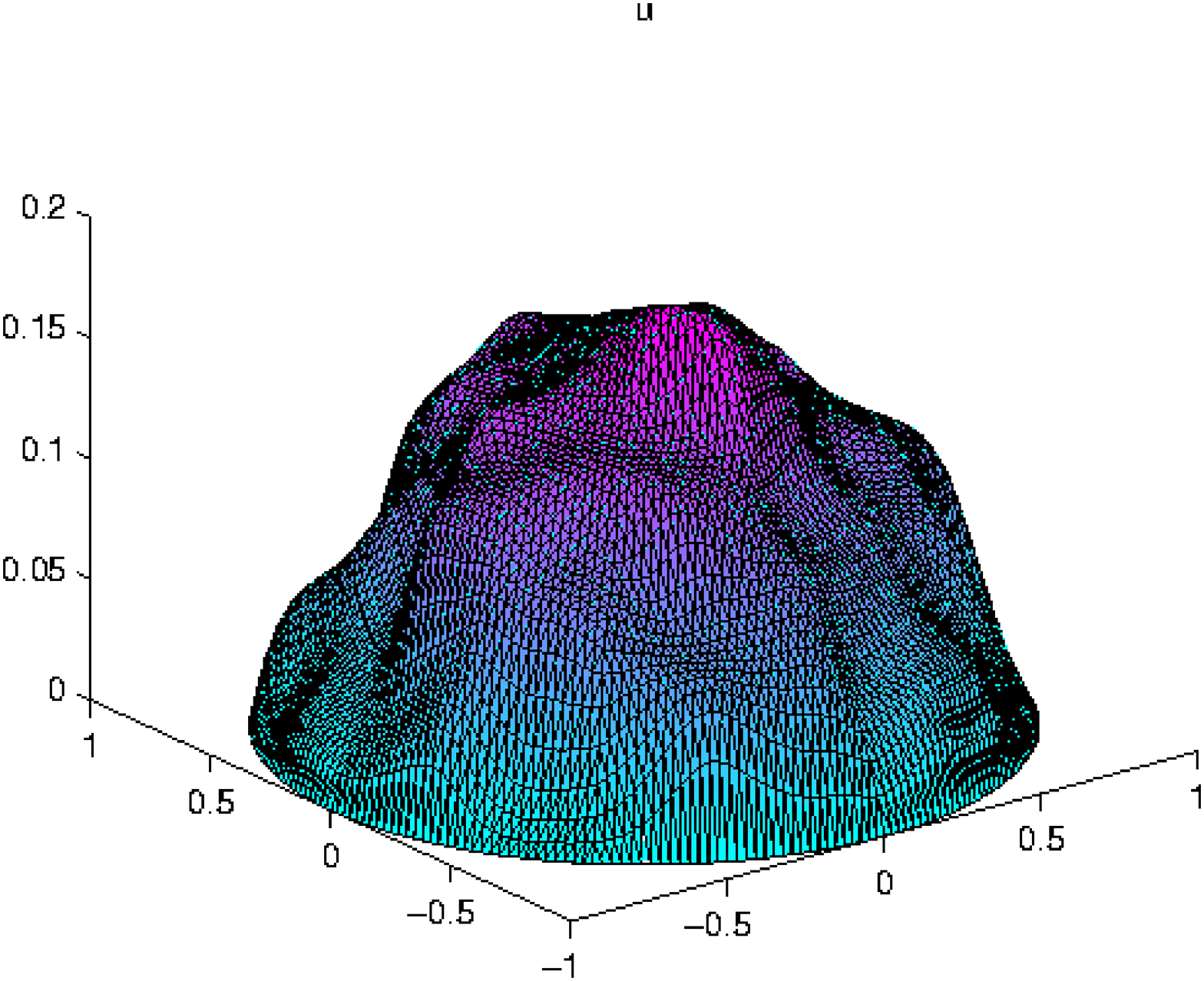}}
     \goodgap
    \subfigure[$u^m$\label{figfffsarrr}]
          {\includegraphics[width=0.3\textwidth,height= 0.3\textwidth]{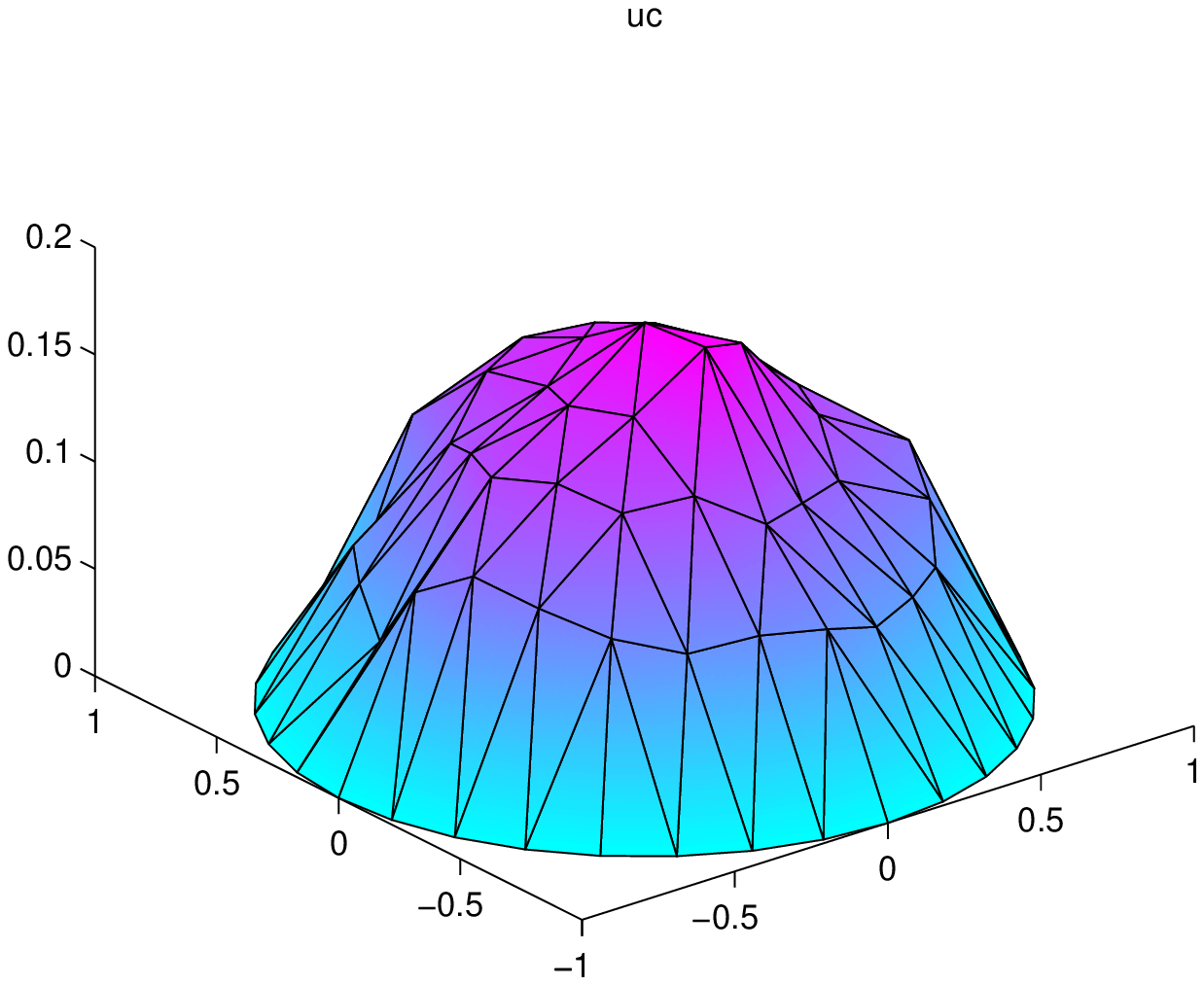}}
    \goodgap
   \subfigure[Triangle of the coarse mesh.\label{figdjdueAfirstBs}]
        {\includegraphics[width=0.3\textwidth,height= 0.3\textwidth]{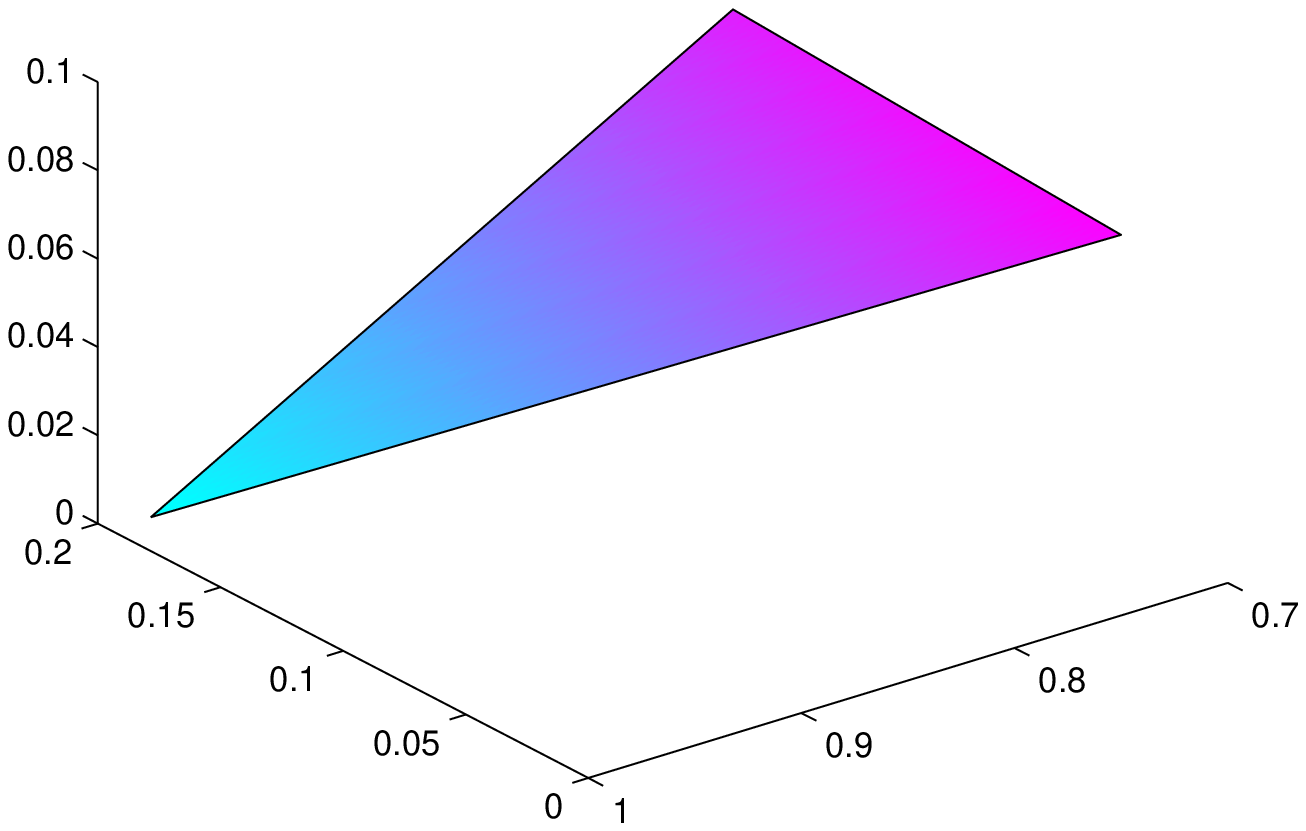}}\\[-10pt]
    \caption{$u$ estimated with the up-scaled operator.}
    \label{figyyassAfirstsss}
\end{center}
\end{figure}
The compressed operator allows us to capture the solution of
\eref{ghjh52}  on a coarse mesh (figure \ref{figyyassAfirstsss}).
What information should be added to the compressed operator in order
to obtain fine resolution approximation of $u$? The answer is a
finer resolution of $F$ (figure \ref{figyyswassAfirst}). Indeed let
$\J_h$ be the interpolation operator introduced in \eref{gdfd541},
we then have the following estimate
\begin{Theorem}\label{sfsfsdd653}
Assume that $\sigma$ and the scheme are stable and that the mesh is
not unadapted to $F$. Then there exist  constants $\alpha,C_m>0$
such that
\begin{equation}\label{sqwsszsjsdsdhdgdgz7112}
\begin{split}
\| u-\J_h u^m\|_{H^1(\Omega)} \leq C_m h^\alpha
\|g\|_{L^\infty(\Omega)}.
\end{split}
\end{equation}
\end{Theorem}
\begin{Remark}\label{sqwff41}
The constant $\alpha$ depends only on $n, \Omega$ and $\mu_\sigma$.
The constant $C_m$ can be written
\begin{equation}
C_m:= C \Big(\frac{\eta^*_{\max} \eta_{\max}}{S^m}\Big)^\frac{1}{2}
\end{equation}
where $C$ depends on the objects mentioned above plus
$\lambda_{\min}(a), \lambda_{\max}(a)$ and
$\big\|\big(\Tr(\sigma)\big)^{-1-\epsilon}\big\|_{L^{1}(\Omega)}$.
\end{Remark}

\subsubsection{Coherent multi-resolution.}\label{gfsgfsfs541}
We want to compress a physical system from a fine scale description
(F) to a coarse scale description (C) there are two ways of doing
so:
\begin{itemize}
\item Either we up-scale directly from (F) to (C)
\item Either we do so in two steps: from (F) to an intermediate
scale (I), then from (I) to (C).
\end{itemize}
Now if the scales (F), (I) and (C) are not completely separated a
technique solely based on averaging (up-scaling of bulk quantities)
would produce two different results (depending on the presence of an
intermediate step or not). Thus it is important to check the
consistency of  the numerical homogenization method if the metric
information $F(x_i)-F(x_j)$ is up-scaled in addition to traditional
bulk quantities $<a\nabla F>_K$.

Let $\T^0, \ldots , \T^n$ be a multi-resolution tessellation of
$\Omega$. Each
 $\T^i$ is a regular conformal tessellation of $\Omega$. Moreover
 $\T^{i+1}$ is a refinement of $\T^i$.
Let us write $V^i$ the space of piecewise linear functions on
$\T^i$. We write $\B^i$ the space of bilinear operators on $\T^i$.
We want to compress (up-scale) the bilinear operator $a[,]$ on the
multi-grid $\T^0, \ldots , \T^n$. We assume that the smallest scale
$n$ is fine enough to capture the irregularities of $a$, in that
case we define $a^n$ such that for all $v,u \in V^n$
\begin{equation}\label{gsffs544131}
a^n[v,u]:= a[v,u].
\end{equation}
Since the gradient of an element of $V^n$ is constant within each
triangle of $\T^n$, $a^n[,]$ can be defined by a mapping from $\T^n$
onto $\M_n$ the space of $n\times n$ constant matrices. We will
write
 $a^n(K)$ the constant matrix associated\footnote{the bilinear form \eref{gsffs544131} can be written as a sum of integrals over $K\in
 \T^n$.
 $\nabla v$ and $\nabla u$  are constant over these triangles, thus $a^n_K$ as a bilinear form over $V^n$ is determined by a matrix.} to  $K\in \T^n$. Similarly
 each bilinear operator of $\B^i$ can be defined by mapping from $\T^i$ onto
$\M_n$. We define for $k\leq  p$, $\U^{k,p}$ the up-scaling operator
mapping $\B^{p}$ onto $\B^{k}$ in the following way (we assume that
the tessellations $\T_i$ are not unadapted to $F$ and that the
respective schemes are stable). Let $B \in \B^{p}$.
 \begin{itemize}
 \item Let $F \in V^p$ be the solution of
 \begin{equation}
 \begin{cases}
 B[v,F]=0 \quad \text{for all}\quad v \in V^p\\
 F=x \quad \text{on}\quad \Gamma^i.
 \end{cases}
 \end{equation}
where we have written $\Gamma^i$ the boundary of $\T^i$.
\item The bilinear form $\U^{k,p} B$ is defined by its matrices $\U^{k,p}
B(K)$ for $K\in \T^k$.
\begin{equation}
\U^{k,p} B(K):= \big< B \nabla F\big>_K \big(\nabla F(K)\big)^{-1}.
\end{equation}
$\big<.\big>_K$ stands for the averaging operator
\begin{equation}
\big< B \nabla F\big>_K:= \frac{1}{\Vol(K)}\sum_{T\in \T^p, T\subset
K} \Vol(T) B(T)\nabla F(T).
\end{equation}
 \end{itemize}
For $k \leq p \leq q$, $\U$ satisfies the semi-group property
\begin{equation}\label{jhdgdd63}
\U^{k,q}=\U^{k,p} \U^{p,q}.
\end{equation}
Note also that $\U^{q,q}=I_d$. In particular if we define for $k\in
\{1,\ldots,n\}$
\begin{equation}
a^k:=\U^{k,n} a^n
\end{equation}
as the up-scaled operator, the following semi-group property is
satisfied: for $k\leq p \leq n$
\begin{equation}\label{gfdgdfd63}
a^k:=\U^{k,p} a^p.
\end{equation}
Semi-group properties \eref{jhdgdd63}, \eref{gfdgdfd63} are
essential to the consistency and coherence of the numerical
homogenization method.

\subsection{Numerical homogenization from the transport point of
view.}\label{hsghsg6531} The elliptic operator appearing in
\eref{ghjh52} can be seen as the generator of a stochastic
differential equation. This stochastic differential equation can
reflect the transport process of a pollutant in a highly
heterogeneous medium such as soil. The following operator $\Delta
-\nabla V\nabla$ whose numerical homogenization is similar to the
one of \eref{ghjh52} can represent a physical system evolving in a
highly irregular energy landscape $V$. The simple fact that this
 evolution taking place in a
continuous domain can be captured by a Markov chain evolving on a
 graph is far from being obvious \cite{MR1473568}. Our point of view here is to
 accurately simulate a Markov chain living on a fine graph by an ``up-scaled'' Markov chain living on a coarse
 graph. The main question is how to choose the
the jump rate $\gamma_{ij}$ of the random walk between the nodes of
the coarse graph? The answer to that question is given by a finite
volume method.

Let us write $\T^*_h$ the dual mesh associated to $\T_h$. $\T^*_h$
can be obtained by drawing segments from the midpoints of the edges
of the triangles of $\T_h$ to an interior point in these triangles
(the circumcenter to obtain a Vorono\"{i} tessellation but one can
also choose the barycenter).

Let us write $V_i$ the control volume associated to the node $i$ of
the primal mesh and $\chi_i$ the characteristic function of $V_i$.
The finite volume method can be expressed in the following way: look
for $u^v\in \Z_h$ ($\Z_h$ being the space spanned by the elements
$\xi_i$ introduced in \eref{jshshsg51}) such that for all $i\in
\No_h$
\begin{equation}\label{hsgs161}
a^*[\chi_i,u^v]=(\chi_i,g)_{L^2}.
\end{equation}
Again, it follows from equation \eref{hsgs161} that the only
information kept from the scales are the usual bulk quantities
(effective conductivities at the edges of the dual mesh) plus the
metric information $F(b)-F(a)$ where $a$ and $b$ are nodes of the
triangles of the primal mesh. Observe also that it is possible to
generate with this finite volume method a coherent multi-resolution
compression similar to the one introduced in subsection
\ref{gfsgfsfs541}. According to \eref{hsgs161} the good choice for
the jump rates of the random walk should be
\begin{equation}\label{hdgdh651}
\gamma_{ij}=a^*[\chi_i,\xi_j]\quad \text{if}\quad i\sim j \quad
\text{and}\quad i\not=j.
\end{equation}
To properly describe the transport process one should look at a
parabolic operator instead of the elliptic one. This issue will be
addressed in \cite{OwZh05b}, we will restrict ourselves to the
elliptic case characterizing the equilibrium properties of the
random walk.  Let us write $\S^v$ the stability of the up-scaled
finite volume operator. It is defined by
\begin{equation}\label{jjhsjhsk726}
\S^v:=\inf_{w\in \Z_h} \sup_{v \in \Y_h}
\frac{a^*[v,w]}{(\ED_h[v])^\frac{1}{2} (\ED_h[w])^\frac{1}{2}}.
\end{equation}

\begin{Theorem}\label{sfsfsdd653sss}
Assume that $\sigma$ and the scheme are stable ($S^v>0$) and that
the mesh is not unadapted to $F$. Then there exist constants
$\alpha,C_v>0$ such that
\begin{equation}\label{sqhdgdgz7112sss}
\begin{split}
\|\I_h u-\I_h u^v\|_{H^1(\Omega)} \leq C_{v,1} h^\alpha
\|g\|_{L^\infty(\Omega)}
\end{split}
\end{equation}
and
\begin{equation}\label{sqhdgdgz7sss112}
\begin{split}
\| u-\J_h u^v\|_{H^1(\Omega)} \leq C_{v,2} h^\alpha
\|g\|_{L^\infty(\Omega)}.
\end{split}
\end{equation}
\end{Theorem}
\begin{Remark}\label{sqwfaswsaf41}
The constant $\alpha$ depends only on $n, \Omega$ and $\mu_\sigma$.
The constant $C_{v,1}$ can be written
\begin{equation}\label{gsfgdfd6521}
C_{v,1}:= C \frac{\eta^*_{\min} \eta_{\max}}{S^m}
\big(\lambda_{\max}(\sigma)\big)^\frac{1}{2}.
\end{equation}
The constant $C_{v,2}$ can be written
\begin{equation}\label{jhdgdd6331}
C_{v,2}:= C \Big(\frac{\lambda_{\max}(\sigma)\eta^*_{\max}
\eta_{\max}}{S^v}\Big)^\frac{1}{2}.
\end{equation}
 $C$ depends on the objects mentioned above plus
$\lambda_{\min}(a), \lambda_{\max}(a)$ and\\
$\big\|\big(\Tr(\sigma)\big)^{-1-\epsilon}\big\|_{L^{1}(\Omega)}$.
\end{Remark}
\begin{Remark}
Observe that we need the additional condition
$\lambda_{\max}(\sigma)<\infty$ to prove the convergence of the
method. Numerical experiments show that although the finite volume
method keeps very few information from small scales it is more
stable and accurate than the method presented in \ref{gfgxfs541} (it
is also more stable and almost as accurate as Galerkin method in
which the whole fine scale scale structure of $F$ is up-scaled).
That is why we believe that the constants \eref{gsfgdfd6521} and
\eref{jhdgdd6331} are not optimal.
\end{Remark}

\subsection{Explicit formulae in laminar cases.}
The harmonic coordinates can be explicitly computed in dimension
one. In this subsection we will analyze a toy model to understand
the effect of the new metric when the coefficients of the partial
differential equations are characterized by an infinite number of
overlapping scales. Our point is to show that new metric becomes
multi-fractal.

 Let $\Omega:=(0,1)$. Let $V\in L^\infty(\Omega)$ and write
$w(V)$ the weak solution of the following Dirichlet problem:
\begin{equation}\label{dgvdtdvt1}
\begin{cases} -\frac{1}{2}e^{2 V}\diiv\big(e^{-2 V } \nabla w\big)=f\\
w=0 \quad \text{on}\quad \partial\Omega.
\end{cases}
\end{equation}
With $f\in L^\infty(\Omega)$. We write $w(V)$ the solution of
\eref{dgvdtdvt1}. Write $\mu_+$ and $\mu_-$  the probability
measures defined on the Borelian subset of $(0,1)$ by
\begin{equation}\label{dgvdttadvt12}
\begin{split}
\mu_+[0,x]:=\frac{\int_0^x e^{2V(z)} dz}{\int_0^1 e^{2V(z)} dz}
\quad \text{and} \quad \mu_-[0,x]:=\frac{\int_0^x e^{-2V(z)}
dz}{\int_0^1 e^{-2V(z)} dz}.
\end{split}
\end{equation}
Define $D(V):=\big(\int_0^1 e^{2V(z)} dz\int_0^1 e^{-2V(z)}
dz\big)^{-1}$. $D(V)$ is to put into correspondence with the bulk
quantity $<a\nabla F>$ which is a constant in dimension one. Observe
that \eref{dgvdtdvt1} can be explicitly rewritten
\begin{equation}\label{dgvssddqttazdvt12}
\begin{split}
-\frac{1}{2}\frac{d}{d\mu_-} D(V) \frac{d}{d\mu_+}u=f.
\end{split}
\end{equation}
In that sense the metric based numerical homogenization is exact in
dimension one. Equations such as \eref{dgvssddqttazdvt12} have been
studied in dimension one in \cite{MR2017701} and \cite{MR1885763} in
order to introduce a measure-theoretic way of
defining differential operators on fractal sets on the real line.\\
 Let $V_n$ be a sequence in $L^\infty(\Omega)$. We call
a probability measure non degenerate if and only if it is atom-less
and the mass of any non void open subset of $\Omega$ is strictly
positive. We have the following theorem
\begin{Theorem}\label{sfgsfd6637}
If $\mu_+^{V_n}$ and $\mu_-^{V_n}$ weakly converge to non degenerate
probability measures $\mu_+$ and $\mu_-$ then $D(V_n)w(V_n)$
converges pointwise to the  unique solution of the following
differential equation with Dirichlet boundary condition
\begin{equation}\label{dgssvssddqttazdvt12}
\begin{split}
-\frac{1}{2}\frac{d}{d\mu_-} \frac{d}{d\mu_+}\psi=f.
\end{split}
\end{equation}
\end{Theorem}
In the case of classical homogenization observe that $\mu_+$ and
$\mu_-$ are simple Lebesgue probability measures on $[0,1]$.

Write $\To$ the torus of dimension one and side of length one and
let $U\in C^1(\To)$. Let  $\rho \in \N /\{0,1\}$. Write $T_{\rho}$
the scaling operator denied on the space of functions by
$T_{\rho}U(x):=U(\rho x)$. Write
\begin{equation}\label{hsghsghgqq55677}
S_n U:=\sum_{p=0}^{n-1} T_{\rho^p} U.
\end{equation}
Take $V_n=S_n U$ in theorem \ref{sfgsfd6637}. Then by
Perron-Frobenius-Ruelle theorem \cite{MR2000h:37063}, $\mu_+^{V_n}$
and $\mu_-^{V_n}$ weakly converge to some probability measure
$\mu_+$ and $\mu_-$ (eigenvectors of the Ruelle transfer operator)
and it is easy to check that they are non degenerate. Let us notice
that similarly it is possible to show that $-\frac{1}{n} \ln D(V_n)$
converges to the sum of topological pressures of $U$ and $-U$ with
respect to the shift induced by the multiplication by $\rho$ on the
space of $\rho$-adic decompositions \cite{Ow00a}, \cite{BeOw00b}.
Theorem \ref{sfgsfd6637} is telling us that the regularity of $\psi$
corresponds to the regularity of $\mu_+[0,x]$ thus it is natural to
wonder what is the regularity of that harmonic measure. To answer
that question we will consider the paradigm of binomial measures. We
refer to \cite{Ri97} for a  detailed introduction to this subject,
for the sake of completeness we will recall its main lines below in
our framework. We take $U(x)\in L^\infty(\To)$ with ($a\not= b$)
\begin{equation}
2 U(x)=
\begin{cases} a \quad \text{for}\quad 0\leq x<\frac{1}{2} \\
b \quad \text{for}\quad \frac{1}{2}\leq x<1.
\end{cases}
\end{equation}
Let us write
\begin{equation}
m_0=\frac{e^a}{e^a+e^b} \quad \text{and}\quad
m_1=\frac{e^b}{e^a+e^b}.
\end{equation}
Then $\mu_+^{V_n}$ weakly converges to $\mu_+$ and that measure is
self-similar in the sense that it satisfies
\begin{equation}\label{dgvsswddqqqzdsvt12}
\begin{split}
\mu_+([x,y])=m_0 \mu_+([2x,2y])+m_1\mu_+([2x-1,2y-1]).
\end{split}
\end{equation}
For $x\in (0,1)$ we write $x_1x_2\ldots$ its $2$-adic decomposition
in the sense that $x=\sum_{n=0}^\infty x_n 2^{n}$. We write
$I_{y_1,\ldots,y_p}$ the cylinder of $x\in (0,1)$ such that
$x_1\ldots x_p=y_1\ldots y_p$. $\mu_+$ is atom-less but singular
with respect to the Lebesgue measure. Moreover it is easy to check
\cite{Ri97} that $\mu_+$ and therefore $\psi$ are Holder continuous
and their Holder continuity exponent is not a constant. More
precisely let us write
\begin{equation}
\alpha_n(x):=-\frac{\ln \big(\mu_+[I_{x_1\ldots x_n}]\big)}{n \ln 2}
\end{equation}
and
\begin{equation}
\alpha(x):=\lim_{n\rightarrow \infty}\alpha_n(x).
\end{equation}
Whenever this limit exists. This limit exists for almost all $x$
(with respect to Lebesgue measure) and its value depends on the
dyadic expansion of $x$. Writing $l_n(x)$ the number of one
appearing in the first $n$ digits of $x$ we have
\begin{equation}
\alpha(x)=\lim_{n\rightarrow \infty}-(1-\frac{l_n(x)}{\ln
n})\log_2(m_0)-\frac{l_n(x)}{\ln n} \log_2(m_1).
\end{equation}
Thus $\alpha(x)$ can take all the values between $-\log_2(m_0)$ and
$-\log_2(m_1)$. However for almost all $x$ with respect to Lebesgue
measure
\begin{equation}
\alpha(x)=-\frac{1}{2}\log_2(m_0 m_1)=-\frac{1}{\ln 2}
\big(\frac{a+b}{2}-\ln (e^a+e^b)\big).
\end{equation}
Now it is possible to obtain from  large deviation theory
(\cite{El85} theorem $2$ and \cite{Ri97}) that
\begin{equation}
\P\Big(\alpha_n(x)\in (\alpha-\epsilon,
\alpha+\epsilon)\Big)\rightarrow 1+c^*(\alpha)
\end{equation}
$\P$ being the uniform probability measure on $(0,1)$,
$c(q)=1-\log_2(m_0^q +m_1^q)$ and
  $c^*$ denotes the Legendre transform of $c$, i.e.
$c^*(\alpha)=\inf_q\big(q\alpha -c(q) \big)$. Thus the metric
associated to our up-scaling method is multi-fractal \cite{FrPa85},
\cite{Mel02}, \cite{MR1244423}. Let us recall that multi-fractal
formalism was originally introduced to describe the regularity of
the velocity field in Turbulence \cite{FrPa85} and explaining
intermittency.

\subsection{Literature.}
The issue of numerical homogenization partial differential equations
with heterogeneous coefficients has received a great deal of
attention and many  methods have been proposed\footnote{Justified in
dimension one,  in the case of periodic or ergodic media with scale
separation or in the case of partial differential equations with
sufficiently smooth coefficients}. Let us mention a few of them.
\begin{itemize}
\item Multi-scale finite element methods \cite{MR1740386}, \cite{MR1613757}, \cite{MR1455261}, \cite{MR2123115}, \cite{MR1232956}, \cite{MR1194543}, \cite{AlBr04}.
\item Multi-scale finite volume methods \cite{JLP03}.
\item Heterogeneous Multi-scale Methods \cite{EV04}.
\item Wavelet based homogenization \cite{MR1492791}, \cite{MR1618846}, \cite{LC04}, \cite{MR1614457}, \cite{MR1614980}, \cite{MR1354913}.
\item Residual free bubbles methods \cite{MR2006324}.
\item Discontinuous enrichment methods \cite{MR2007030}, \cite{MR1870426}.
\item Partition of Unity Methods \cite{FY05}.
\item Energy Minimizing Multi-grid Methods \cite{MR1756048}.
\end{itemize}
The methods mentioned above are part of a larger quest aimed at
capturing   high dimensional problems with a few coarse parameters
\cite{MR1675219}, \cite{MR2050644}, \cite{Kev05}, \cite{MR619904} .
Paraphrasing the outcome of a recent DOE workshop \cite{DOE04}, we
may understand the physics of multi-scale structures at each
individual scale nevertheless \emph{without the capability to
"bridge the scales," a significant number of important scientific
and engineering problems will remain out of reach}.

\section{Proofs.}

\subsection{Compensation.}
Let us  prove theorem  \ref{th2}. We need a variation of Campanato's
result \cite{CM5} on non-divergence form elliptic operators. Let us
write for a symmetric matrix $M$,
\begin{equation}
\nu_M:=\frac{\sum_{i=1}^n \lambda_{i,M}}{\sum_{i=1}^n
\lambda_{i,M}^2}.
\end{equation}
We consider the following Dirichlet problem:
\begin{equation}\label{dcaqssaslkwaq21}
L_M v=f
\end{equation}
with $L_M:=\sum_{i,j=1}^n M_{ij}(x) \partial_i \partial_j$. We
assume $M$ to be elliptic and symmetric.
\begin{Theorem}\label{hdgjhdgd7}
Assume that $\beta_M<1$. If $\Omega$ is convex, then, there exists a
real number $p>2$ depending only on $n,\Omega$ and $\beta_M$ such
that if $f\in L^p(\Omega)$ the Dirichlet problem
\eref{dcaqssaslkwaq21} has a unique solution satisfying
\begin{equation}\label{hdhdgc7}
\|v\|_{W^{2,p}_0(\Omega)}\leq \frac{C }{1-\beta_M^\frac{1}{2}}
\|\nu_M f\|_{L^p(\Omega)}.
\end{equation}
\end{Theorem}
\begin{Remark}
$\beta_M$ is the Cordes parameter \eref{sshgdd7641} associated to
$M$.
\end{Remark}
\begin{proof}
 Theorem \ref{hdgjhdgd7} is a straightforward adaptation of theorem 1.2.1 of
 \cite{MPG00}, for the sake of  completeness we will give the main lines of ideas leading to
  estimate \eref{hdhdgc7}. Let us recall the Miranda-Talenti estimate \cite{MPG00}.
\begin{Lemma}\label{dcaqssawappmbq21}
Let $\Omega\subset \R^n$ be a bounded and convex domain of class
$C^2$. Then for each $v\in W^{2,2}_0(\Omega)$ it results
\begin{equation}\label{sgfdfd651}
\int_{\Omega} \sum_{i,j=1}^n (\partial_i \partial_j v)^2\,dx \leq
\int_{\Omega} (\Delta v)^2\,dx.
\end{equation}
\end{Lemma}
The Laplacian $\Delta: W^{2,p}_0(\Omega)\rightarrow L^p(\Omega)$ is
an isomorphism for each $p>1$. Let $\Delta^{-1}(p)$ be the inverse
operator $\Delta^{-1}: L^p(\Omega)\rightarrow W^{2,p}_0$. It is
clear from \eref{sgfdfd651} that $\|\Delta^{-1}(2)\|\leq 1$. Let
$r\in (2,\infty)$, by the convexity of the norms we have
\begin{equation}
\|\Delta^{-1}(p)\|\leq C(p)
\end{equation}
with
\begin{equation}
 C(p):=\|\Delta^{-1}(r)\|^{\frac{r(p-2)}{p(r-2)}}.
\end{equation}
 Let $v$ be a solution of \eref{dcaqssaslkwaq21} (we refer to \cite{MPG00} for the existence of $v$ which is obtained  from a fix point theorem), we have
\begin{equation}\label{dcssdaqsaq21}
\|v\|_{W^{2,p}_0(\Omega)} \leq \|\Delta^{-1}(p)\| \|\Delta
v\|_{L^p(\Omega)}.
\end{equation}
Observing that $\Delta v=\nu_M f +\Delta v - \nu_M L_M v$ one can
obtain
\begin{equation}\label{dcssdasdqqsaq21}
\|\Delta v\|_{L^p(\Omega)}\leq \|\nu_M f\|_{L^p(\Omega)}+ \|\Delta v
- \nu_M L_M v\|_{L^p(\Omega)}.
\end{equation}
Then following the proof of theorem theorem 1.2.1 of \cite{MPG00} we
have
\begin{equation}\label{dcaqsaszaq21}
\|\Delta v - \nu_M L_M v\|_{L^p(\Omega)}^p \leq \int_{\Omega}
\beta_M^{p/2}\big(\sum_{i,j=1}^n (\partial_i \partial_j
v)^p\big)\,dx.
\end{equation}
Let us choose $p>2$ such that $1-C(p) \beta_M^{1/2} \geq
(1-\beta_M^{1/2})/2$. Combining \eref{dcssdaqsaq21},
\eref{dcssdasdqqsaq21} and \eref{dcaqsaszaq21} we obtain that
\begin{equation}\label{dcaqsswaszaq21}
\Big( \int_{\Omega}  \big(\sum_{i,j=1}^n (\partial_i \partial_j
v)^p\big)\,dx\Big)^\frac{1}{p} \leq \frac{2
C(p)}{1-\beta_M^{1/2}}\|\nu_M f\|_{L^p(\Omega)}.
\end{equation}
Which leads to estimate \eref{hdhdgc7}.
\end{proof}
Remembering the Sobolev embedding inequality
\begin{equation}\label{hdhsadgc7}
\|\nabla v\|_{C^{1-\frac{n}{p}}(\bar{\Omega})}\leq C
\|v\|_{W^{2,p}_0(\Omega)}
\end{equation}
theorem \ref{hdgjhdgd7} implies the H\"{o}lder continuity of $v$ in
dimension $n=2$.

We assume that $\sigma$ is stable. We write $F^{-1}$ the inverse of
$F$ (which is well defined if $\sigma$ is stable). Let us write $Q$
the symmetric positive matrix given by the following equation
\begin{equation}\label{sddasjdsdhdgdgz7112}
Q(y):= \Big(\frac{\big({^t\nabla F} a \nabla F\big)}{|\det
\big(\nabla F\big)|}\Big)\circ F^{-1}(y).
\end{equation}
Let us write $w$ the solution of the following equation: for all
$\hat{\varphi}\in C_0^\infty$,
\begin{equation}\label{sghsaddjh52}
 \sum_{i,j=1}^n Q_{ij}
\partial_i\partial_j w=-
\frac{\hat{g}}{|\det \big(\nabla F\big)|\circ F^{-1}}.
\end{equation}
Let us now prove the following theorem,
\begin{Theorem}\label{th3}
Assume that $\sigma$ is stable and that $\Omega$ is convex. Then
there exists
 constants $p>2$, $C>0$  such that
the solution of \eref{sghsaddjh52} belongs to $W^{2,p}_0(\Omega)$
and satisfies
\begin{equation}
\|\nabla w\|_{W^{2,p}_0(\Omega)}\leq
\frac{C}{1-\beta_\sigma^\frac{1}{2}} \|g\|_{L^\infty(\Omega)}.
\end{equation}
\end{Theorem}
\begin{Remark}
$\alpha$ depends on $\Omega$ and $\beta_\sigma$.  $C$  depends on
$\lambda_{\min}(a)$ and if $n\leq 4$ on
$\big\|\big(\Tr(\sigma)\big)^{\frac{n}{2}-2-\epsilon}\big\|_{L^{1}(\Omega)}$.
\end{Remark}
\begin{proof}
 Now
let us observe that
\begin{equation}\label{jkjh8}
\frac{\nu_Q}{\det(\nabla F)\circ
F^{-1}}=\frac{\Tr(\sigma)}{\Tr(\sigma^2)}\circ F^{-1}.
\end{equation}
Using the change of variables $y=F(x)$ and choosing $1/q'+1/q=1$ we
obtain that
\begin{equation}
\|\frac{\nu_Q}{\det(\nabla F)\circ F^{-1}}
\hat{g}\|_{L^p(\Omega)}\leq \|g\|_{L^{pq'}(\Omega)}
\big\|\frac{\Tr(\sigma)}{\Tr(\sigma^2)} (\det(\nabla
F))^\frac{1}{pq}\big\|_{L^{pq}(\Omega)}.
\end{equation}
It is easy to check that
\begin{equation}\label{hgsfddhgf651}
\big\|\frac{\Tr(\sigma)}{\Tr(\sigma^2)} (\det(\nabla
F))^\frac{1}{pq}\big\|_{L^{pq}(\Omega)}^{pq} \leq
\frac{C_{qp,n}}{(\lambda_{\min}(a))^\frac{n}{2}} \int_{\Omega}
\big(\Tr(\sigma)\big)^{\frac{n}{2}-qp}.
\end{equation}
For $2\leq n\leq 4$  we choose $q=1$ in \eref{hgsfddhgf651} for
$n\geq 5$ we choose $q=\frac{n}{2p}$. Then a direct application of
theorem \ref{hdgjhdgd7} and estimate \eref{jkjh8} to equation
\eref{sghsaddjh52} implies the theorem (observe that
$\beta_Q=\beta_{\sigma}$).
\end{proof}
Let $\varphi \in C_0^\infty(\Omega)$. Write $\hat{\varphi}:=\varphi
\circ F^{-1}$. Using theorem \ref{th3} we obtain that
\begin{equation}\label{sghsddxddjh52}
\Big( \hat{\varphi}, \sum_{i,j=1}^n Q_{ij}
\partial_i\partial_j w\Big)_{L^2(\Omega)}=-\Big( \hat{\varphi},
\frac{\hat{g}}{|\det \big(\nabla F\big)|\circ F^{-1}}
\Big)_{L^2(\Omega)}.
\end{equation}
Using the change of variable $y=F(x)$ we deduce that
\begin{equation}\label{sghsddwexsddjh52}
\Big( \varphi, \sum_{i,j=1}^n \sigma_{ij} (\partial_i\partial_j
w)\circ F \Big)_{L^2(\Omega)}=-\Big( \varphi, g \Big)_{L^2(\Omega)}.
\end{equation}
Let us observe that
\begin{equation}\label{sghsddwexwssddjh52}
\sum_{i,j=1}^n \sigma_{ij} (\partial_i\partial_j w)\circ F =\diiv
\Big(a\nabla F \big((\nabla w)\circ F\big)\Big).
\end{equation}
It follows after an integration by parts that (and observing that
$\nabla F (\nabla w)\circ F=\nabla (w\circ F)$)
\begin{equation}\label{sghazexwssddjh52}
a[\varphi, w\circ F]=(\varphi,g)_{L^2(\Omega)}
\end{equation}
It follows from the uniqueness of the solution of the Dirichlet
problem \eref{sghazexwssddjh52} that $w\circ F=u$. Theorem
\ref{th22} is then a straightforward consequence of theorem
\ref{th3} and the equality $u\circ F^{-1}=w$.

\begin{Remark}
 Using the change of variables $y=F(x)$ we obtain for
all $\varphi \in C^{\infty}_0(\Omega))$
\begin{equation}\label{ghsxdsdjsh52}
a[\varphi,u]= Q[\hat{\varphi},\hat{u}]
\end{equation}
The comparison between \eref{sghsddxddjh52} and \eref{ghsxdsdjsh52}
indicates that $Q$ is a divergence free matrix. We have not used
that property of $Q$ explicitly in our proof above but it is present
implicitly in the deduction of \eref{sghazexwssddjh52} from
\eref{sghsddwexwssddjh52}.
\end{Remark}

\begin{Remark}
The only place where we use the convexity of $\Omega$ is for the
validity of lemma \ref{dcaqssawappmbq21} (we refer to \cite{MPG00}).
\end{Remark}

The following lemma is a well known result obtained from the De
Giorgi-Moser-Nash theory (\cite{MR0093649}, \cite{MR0159138},
\cite{MR0100158}) of divergence form elliptic operators with
discontinuous coefficients (more precisely we refer to
\cite{MR0352696} for the Global H\"{o}lder regularity)

\begin{Lemma}\label{lem1}
There exists $C, \alpha'>0$ depending on $\Omega$ and
$\lambda_{\max}(a)/\lambda_{\min}(a)$ such that $F$ is $\alpha'$
H\"{o}lder continuous and
\begin{equation}
\|F\|_{C^{\alpha'}}\leq C.
\end{equation}
\end{Lemma}

Theorem  \ref{th2} is a straightforward consequence of the Sobolev
embedding inequality \eref{hdhsadgc7}, theorem \ref{th3}, lemma
\ref{hdgjhdgd7} and the fact that $\nabla_F u=\nabla \hat{u}\circ
F$. Let us observe that in dimension two, we have
\begin{equation}
\frac{1}{1-\beta_\sigma}=\frac{1}{2}(\mu_\sigma+\frac{1}{\mu_\sigma}).
\end{equation}
And the condition $\beta_\sigma<1$ is equivalent to
$\mu_{\sigma}<\infty$.

\subsubsection{H\"{o}lder continuity for $n\geq 3$ or non-convexity
of $\Omega$.} In this subsection we will not assume $\Omega$ to be
convex.  Let $N^{p,\lambda}(\Omega)$ $(1<p<\infty,\, 0<\lambda<n$)
be the weighted Morrey space formed by  functions
$v:\Omega\rightarrow \R$ such that
$\|v\|_{N^{p,\lambda}(\Omega)}<\infty$ with
\begin{equation}
\|v\|_{N^{p,\lambda}(\Omega)}=\sup_{x_0\in
\Omega}\Big(\int_{\Omega}|x-x_0|^{-\lambda} |v(x)|^p
\Big)^\frac{1}{p}.
\end{equation}
To obtain the H\"{o}lder continuity of $u\circ F^{-1}$ in dimension
$n\geq 3$ we will use corollary 4.1 of \cite{MR1903306}. We will
give  the result of S. Leonardi below in a form adapted to our
context. Consider the Dirichlet problem \eref{dcaqssaslkwaq21}. We
do not assume $\Omega$ to be bounded. We write
$W^{2,p,\lambda}(\Omega)$ the functions in $W^{2,p}(\Omega)$ such
that their second order derivatives are in $N^{p,\lambda}(\omega)$.
\begin{Theorem}\label{ksjhs721}
There exist a constant $C^*=C^*(n,p,\lambda,\partial \Omega)>0$ such
that if $\beta_M <C^*$ and $f\in N^{p,\lambda}(\Omega)$ then the
Dirichlet problem \eref{dcaqssaslkwaq21} has a unique solution in
$W^{2,p,\lambda}\cap W^{1,p}_0(\Omega)$. Moreover, if $0<\lambda
<n<p$ then $\nabla v \in C^{\alpha}(\Omega)$ with $\alpha=1-n/p$ and
\begin{equation}
\|\nabla v\|_{C^\alpha(\Omega)}\leq
\frac{C}{\lambda_{\min}(M)}\|f\|_{N^{p,\lambda}(\Omega)}
\end{equation}
where $C=C(n,p,\lambda,\partial \Omega)$.
\end{Theorem}
Theorem \ref{ksjhgfts721} is a straightforward application of
theorem \ref{ksjhs721}.

\subsection{Dimensionality reduction.}\label{dcaqsaupq21}
Let us prove theorem \ref{sfsfssaaadd653}. We write $X_h$ the linear
space spanned by the elements $\psi_i$. The solution of the Galerkin
scheme satisfies $a[u-u_h,v]=0$ for all $v\in X_h$. Thus
\begin{equation}
a[u-u_h]=\inf_{v\in X_h}a[u-u_h,u-v].
\end{equation}
It follows by Cauchy-Schwartz inequality that
\begin{equation}\label{sddsjdsdhdgdgz7112}
\begin{split}
a[u-u_h]\leq \inf_{v\in X_h}a[u-v].
\end{split}
\end{equation}
Now writing $\hat{v}:=v\circ F^{-1}$ and using the change of
variable  $y=F(x)$ we obtain
\begin{equation}\label{sddadsjdsdhdgdddgz7112}
a[u-v]=Q[\hat{u}-\hat{v}].
\end{equation}
It follows in dimension $n=2$ that
\begin{equation}\label{sddsjsdsdhdgdgz7112}
\begin{split}
\|u-u_h\|_{H^1}^2\leq \frac{D}{\lambda_{\min}(a)}\inf_{w \in V_h}
\|\nabla \hat{u}-\nabla w\|_{L^\infty(\Omega)}^2
\end{split}
\end{equation}
with
\begin{equation}\label{sddsjsdsdsshdgdgz7112}
\begin{split}
D:=\Tr\Big[\int_{\Omega} {^t \nabla F} a\nabla F\Big].
\end{split}
\end{equation}
Thus using the following standard approximation properties of the
elements $\varphi_i$ (see for instance \cite{ErGu04}),
\begin{equation}\label{sswdgddcz7112}
\begin{split}
\inf_{w \in V_h} \|\nabla \hat{u}-\nabla w\|_{L^2(\Omega)}\leq C
\gamma(\T_h) h^\alpha \|\hat{u}\|_{C^{1,\alpha}_0(\Omega)}
\end{split}
\end{equation}
we obtain that

\begin{equation}\label{sddsdsxsjsdsdhdgdgz7112}
\begin{split}
\|u-u_h\|_{H^1}\leq \gamma(\T_h)
\big(\frac{D}{\lambda_{\min}(a)}\big)^\frac{1}{2}\|\nabla
\hat{u}\|_{C^\alpha} h^\alpha.
\end{split}
\end{equation}
We conclude by observing that for $l\in \R^n$
\begin{equation}\label{sddsjsxxdsdsshdgdgz7112}
\begin{split}
\int_{\Omega} {^tl}{^t \nabla F} a\nabla Fl=\inf_{f\in
C^\infty_0(\Omega)} \int_{\Omega} {^t (l+\nabla f)} a (l+\nabla f).
\end{split}
\end{equation}
Theorem \ref{sfsfssaaadd653} becomes a direct consequence of
\eref{sddsjsdsdhdgdgz7112} and theorem \ref{th3}. \\
In dimension $n\geq 3$ we obtain from \eref{sddadsjdsdhdgdddgz7112}
that
\begin{equation}\label{sddsjsswdgz7112}
\begin{split}
\|u-u_h\|_{H^1}^2\leq
\frac{\lambda_{\max}(Q)}{\lambda_{\min}(a)}\inf_{w \in V_h} \|\nabla
\hat{u}-\nabla w\|_{L^2(\Omega)}^2.
\end{split}
\end{equation}
It is easy to obtain that
\begin{equation}\label{eqforQ}
\begin{split}
\lambda_{\max}(Q)\leq \big(\det(a)\big)^\frac{1}{2}
\mu_{\sigma}^\frac{n}{2}\big(\Tr(\sigma)\big)^{1-\frac{n}{2}}.
\end{split}
\end{equation}
We conclude by observing that $\mu_{\sigma}< C(\beta_\sigma)$ and
using the following standard approximation properties of the
elements $\varphi_i$ (see for instance \cite{ErGu04}).
\begin{equation}\label{sswdgz7112}
\begin{split}
\inf_{w \in V_h} \|\nabla \hat{u}-\nabla w\|_{L^2(\Omega)}\leq C
\gamma(\T_h) h \|\hat{u}\|_{W^{2,2}_0(\Omega)}.
\end{split}
\end{equation}

\subsection{Galerkin with localized elements.}\label{dcsswaqsaq21}
Let us prove theorem \ref{sfsfssasadd653}. We assume that the coarse
mesh is not unadapted to $F$. Let $K$ be a triangle of $\T_h$ and
let $a$ be a node of $K$ such that $\eta^F_{\min}(K)=\frac{1}{\sin
\theta}$ where $\theta$ is the interior angle between $(F(a),F(b))$
and $(F(a),F(c))$; ($b,c)$ being the other nodes of $K$. Let us
prove the following lemma
\begin{Lemma}\label{dgdgqwq536f}
\begin{equation}
\begin{split}
 |\nabla_F u(K)-\nabla_F u(a)|\leq  3 \eta^F_{\min}(K) \|\nabla \hat{u}\|_{C^\alpha}
 \|F\|_{C^{\alpha'}}^\alpha h^{\alpha \alpha'}.
\end{split}
\end{equation}
\end{Lemma}
\begin{proof}
It is easy to check that
\begin{equation}
\begin{split}
u(b)-u(a)=\big(F(b)-F(a)\big)\nabla \hat{u}\circ
F(a)+(F(b)-F(a))\cdot q_{ba}.
\end{split}
\end{equation}
Where the vector $q_{ba}$ is defined by
\begin{equation}\label{jhshs51}
\begin{split}
q_{ba}:=  \int_{0}^1 \Big[\nabla \hat{u}\big[F(a)+s (F(b)-F(a))\big
]-\nabla \hat{u}\big[F(a)\big]\Big]\,ds.
\end{split}
\end{equation}
We will use the notation $f_{ba}:=(F(b)-F(a))/|F(b)-F(a)|$. We will
write $f_{ba}^{\perp}$ the unit vector obtained by  a $90^o$
rotation of $f_{ba}$ towards $f_{ca}$. Defining  $q_{ca}$ as in
\eref{jhshs51} we obtain that
\begin{equation}
\begin{split}
 \nabla_F u(K)=\nabla_F u(a) + k
\end{split}
\end{equation}
with
\begin{equation}
k=q_{ba}-\lambda f_{ba}^{\perp}
\end{equation}
with
\begin{equation}
\lambda:=\frac{f_{ca}.(q_{ba}-q_{ca})}{f_{ca}.f_{{ba}^{\perp}}}.
\end{equation}
Which leads us to
\begin{equation}
\begin{split}
 |\nabla_F u(K)-\nabla_F u(a)|\leq  \frac{3}{f_{ca}.f_{{ba}^{\perp}}} \|\nabla \hat{u}\|_{C^\alpha}
 \|F\|_{C^{\alpha'}}^\alpha h^{\alpha \alpha'}.
\end{split}
\end{equation}
\end{proof}
The following lemma is a direct consequence of lemma
\ref{dgdgqwq536f}

\begin{Lemma}\label{dgdssgqaswq536f}
Let $K\in \T_h$ and let $x\in \Omega$  then
\begin{equation}
\begin{split}
 |\nabla_F u(K)-\nabla_F u(x)|\leq   3 \eta^*_{\min} \|\nabla \hat{u}\|_{C^\alpha}
 (1+\|F\|_{C^{\alpha'}}^\alpha) \big(h+\dist(x,K)\big)^{\alpha
 \alpha'}.
\end{split}
\end{equation}
\end{Lemma}
Let us write $\Z_h u$ the interpolation of $u$ over the space $Z_h$:
\begin{equation}\label{dcdcssssaq221}
\Z_h u(x)=\sum_{i\in \No_h} u(x_i) \xi_i(x).
\end{equation}
\begin{Lemma}\label{lem3}
We have
\begin{equation}
a^*[u-\Z_h u]\leq C \eta^*_{\min} \|\nabla \hat{u}\|_{C^\alpha}
 \|F\|_{C^{\alpha'}}^\alpha h^{\alpha \alpha'} D.
\end{equation}
\end{Lemma}
\begin{proof}
We have
\begin{equation}
a_K[u-\Z_h u]=\int_K {^t\big(} \nabla_F u(x) -\nabla_F u(K)\big)
\sigma(x) \big( \nabla_F u(x) -\nabla_F u(K)\big)\,dx
\end{equation}
with $\sigma(x):={^t\nabla F a\nabla F}$. Using the change of
variables $F(x)=y$ we obtain that
\begin{equation}
a_K[u-\Z_h u]=\int_{F(K)} {^t\big(} \nabla \hat{u}(y) -\nabla_F
u(K)\big) Q(y) \big( \nabla \hat{u}(y) -\nabla_F u(K)\big)\,dy
\end{equation}
from which we deduce that
\begin{equation}
a_K[u-\Z_h u]\leq \big( 3 \eta^*_{\min} \|\nabla
\hat{u}\|_{C^\alpha}
 \|F\|_{C^{\alpha'}}^\alpha h^{\alpha \alpha'}\big)^2
 \int_{F(K)}\sup_{|e|=1}{^teQe}.
\end{equation}
Thus
\begin{equation}
a^*[u-\Z_h u]\leq C \big(\eta^*_{\min} \|\nabla \hat{u}\|_{C^\alpha}
 \|F\|_{C^{\alpha'}}^\alpha h^{\alpha \alpha'}\big)^2 D
\end{equation}
where $D$ has been defined by \eref{sddsjsdsdsshdgdgz7112}.
\end{proof}
Theorem \ref{sfsfssasadd653} is implied by lemma \ref{lem3}, theorem
\ref{th3}, lemma \ref{lem1} and the  following inequality
\begin{equation}
a^*[u-u^f]\leq a^*[u-\Z_h u].
\end{equation}
Let us now prove theorem \ref{sfsasd653}. By the triangle inequality
\begin{equation}\label{hsgs6510}
a\big[u- \J_h u^f\big]\leq a\big[u- \J_h u\big]+a\big[\J_h u- \J_h
u^f\big].
\end{equation}
We write $\hat{\J}_h u:= (\J_h u)\circ F^{-1}$. $\hat{\J}_h u$ is a
linear interpolation of $\hat{u}$ on the tessellation $\T^F$. Now
using the identity
\begin{equation}
a[u- \J_h u]= Q[\hat{u}- \hat{\J}_h u]
\end{equation}
we obtain that
\begin{equation}\label{hsgs541}
a[u- \J_h u]\leq \|\hat{u}\|_{C^\alpha} \hat{h}^\alpha D.
\end{equation}
Where we have written $\hat{h}$ the maximal length of the edges of
$\T^F$. Observe that $\hat{h}\leq h^{\alpha'}\|F\|_{C^{\alpha'}}$.
Theorem \ref{sfsasd653} is a consequence of inequalities
\eref{hsgs6510}, \eref{hsgs541}, lemmas \ref{hsjhgs51} and
\ref{lem3}, theorem \ref{th3} and  the following inequality
$$a^*\big[\Z_h u-
 u^f\big]\leq  2 a^*\big[ u-
 \Z_h u\big].$$

\begin{Lemma}\label{hsjhgs51}
\begin{equation}\label{jshsgshg413}
a\big[\J_h u- \J_h u^f\big]\leq 2^{6}\eta^*_{\max}
\mu_\sigma^\frac{1}{2}\eta_{\max}^3
\big(\frac{\lambda_{\max}(a)}{\lambda_{\min}(a)}\big)^\frac{1}{2}
 a^*\big[\Z_h u-
 u^f\big]
\end{equation}
and
\begin{equation}\label{jaqgshg413}
a\big[\J_h u- \J_h u^f\big]\leq \mu_{\sigma}\nu^* a^*\big[\Z_h u-
 u^f\big]
\end{equation}
\end{Lemma}
\begin{proof}
Let us write $w:=\J_h u- \J_h u^f$. We need to bound $a[w]$.
 We have
\begin{equation}
a[w]=\frac{a^*[\J_h w]}{a^*[\Z_h w]} a^*[\Z_h w].
\end{equation}
 $a^*[\Z_h w]=a^*[\Z_h u - u^f]$  has already been estimated
in lemma \ref{lem3}. Observe that
\begin{equation}\label{eqq1}
a^*[\J_h w]=Q[\hat{w}].
\end{equation}
$\hat{w}$ is piecewise linear on $\T^F$. Using property
\eref{ddff351ssad32435} we obtain that
\begin{equation}\label{eqq2}
 Q[\hat{w}]
\leq \eta^*_{\max} \lambda_{\max}(Q) \ED_h[w].
\end{equation}
Moreover, observing that
\begin{equation}\label{gsfssw41}
Q\circ
F=\frac{\sigma}{\big(\det(\sigma)\big)^\frac{1}{2}}\big(\det(a)\big)^\frac{1}{2}
\end{equation}
we obtain that
\begin{equation}\label{sjhsggs6}
\lambda_{\max}(Q)\leq \mu_{\sigma}^\frac{1}{2}
\big(\det(a)\big)^\frac{1}{2}.
\end{equation}
Equation \eref{sjhsggs6} is valid in dimension $2$, in dimension
higher than $2$ we would use the inequality \eref{eqforQ}
\begin{equation}
\lambda_{\max}(Q)\leq \mu_{\sigma}^\frac{1}{2}
\big(\lambda_{\min}(\sigma)\big)^{\frac{n-2}{2}}
\big(\det(a)\big)^\frac{1}{2}.
\end{equation}
We now need to bound from below $a^*[\Z_h w]/\ED_h[w]$. For $K\in
\T_h$
 let us write $H(K)$ the matrix
\begin{equation}
H(K):=\int_{K}{^t (}\nabla F(K))^{-1}{^t \nabla F(x)}a\nabla F(x)
(\nabla F(K))^{-1}\,dx.
\end{equation}
We need to estimate  $\inf_{|l|=1} {^t l} H(K)l$. Let us write
$a,b,c$ the nodes of $K$ and
\begin{equation}
f(x):=\big(F(x)-F(a)\big) \big(\nabla F(K)\big)^{-1}\cdot l.
\end{equation}
Let us observe that
\begin{equation}
{^t l} H(K)l=\int_{K}{^t \nabla f}a\nabla f.
\end{equation}
Moreover $f(a)=0$, $f(b)=(b-a).l$ and $f(c)=(c-a).l$. Let us assume
without loss of generality that $|f(b)|/|b-a|\geq |f(c)|/|c-a|$.
Then
\begin{equation}\label{hgdfd6541}
{^t l} H(K)l \geq \inf_{w\in C^\infty(\Omega),\, w(a)=0,\,
w(b)=f(b)} \int_{K}{^t \nabla w}a\nabla w.
\end{equation}
The quantity appearing in \eref{hgdfd6541} is the resistance metric
distance between the points $a$ and $b$ and it is easy to check that
(see for instance lemma 1.1 of \cite{BCK05})
\begin{equation}
\inf_{w\in C^\infty(\Omega),\, w(a)=0,\, w(b)=f(b)} \int_{K}{^t
\nabla w}a\nabla w  \geq \lambda_{\min}(a) \Vol(K)
\big(\frac{|f(b)|}{|b-a|}\big)^2.
\end{equation}
Thus
\begin{equation}
{^t l} H(K)l\geq \lambda_{\min}(a) \Vol(K)
\big(\frac{|(b-a).l|}{|b-a|}\big)^2.
\end{equation}
Let us observe that
\begin{equation}
\frac{|(b-a).l|}{|b-a|} \geq \frac{1}{4\eta_{\max}}.
\end{equation}
It follows that
\begin{equation}
a^{*}[\Z_h w]\geq \|\nabla \I_h w\|_{L^2}^2 \lambda_{\min}(a)
\frac{1}{16\eta_{\max}^2}.
\end{equation}
Thus
\begin{equation}\label{kdjdhj6612}
\frac{a^{*}[\Z_h w]}{\ED_h[w]}\geq  \lambda_{\min}(a)
\frac{1}{2^{6}\eta_{\max}^3}.
\end{equation}
which leads us to equation \eref{jshsgshg413}.

\begin{Remark}
One of the methods employed with non-conformal elements to ensure
the stability and convergence of the scheme is the so called
\emph{patch test}. In our proof the stability condition and
convergence is ensured by \eref{kdjdhj6612} and an uniform lower
bound on $\eta_{\max}$.
\end{Remark}

To obtain \eref{jaqgshg413} let us observe that
\begin{equation}
a[w]=Q[\hat{w}].
\end{equation}
Thus
\begin{equation}
a[w]=\sum_{K\in \T_h}\int_{K^F} {^t\nabla \hat{w}(K^F)}Q(y) \nabla
\hat{w}(K^F)\,dy.
\end{equation}
It follows that
\begin{equation}
a[w]\leq  \nu^* \sum_{K\in \T_h}\int_{F(K)}
\frac{\lambda_{\max}(Q)}{\lambda_{\min}(Q)} {^t\nabla
\hat{w}(K^F)}Q(y) \nabla \hat{w}(K^F)\,dy.
\end{equation}
from equation \eref{gsfssw41} we obtain that
\begin{equation}\label{sjhsggs6w}
\frac{\lambda_{\max}(Q)}{\lambda_{\min}(Q)} \leq \mu_{\sigma}.
\end{equation}
Next, observing that
\begin{equation}
\sum_{K\in \T_h}\int_{F(K)} {^t\nabla \hat{w}(K^F)}Q(y) \nabla
\hat{w}(K^F)\,dy= a^*[\Z_h w].
\end{equation}
we obtain \eref{jaqgshg413}.
\end{proof}

\subsection{Numerical homogenization from the information point of view}
In this subsection we will prove theorem \ref{sfsfxadd653} and
theorem \ref{sfsfsdd653}. The method introduced in subsection
\ref{gfgxfs541} can be formulated in the following way: look for
$u^m \in V_h$ such that for all $i\in \No_h$,
\begin{equation}\label{gdhgdfctetfc5}
a^*[\varphi_i, \Z_h u^m ]=(\varphi_i,g)_{L^2{\Omega}}
\end{equation}
which implies the following finite volume orthogonality property for
all $i\in \No_h$,
\begin{equation}
a^*[\varphi_i, \Z_h u^m-u]=0.
\end{equation}
Let us write $w= u- \Z_h u^m$. By equation \eref{gsgfd763} we obtain
that
\begin{equation}
\big(\ED_h[w]\big)^\frac{1}{2} \leq \frac{1}{S^m} \sup_{v\in Z_h}
\frac{a^*[v,\Z_h w]}{\big(\ED_h[v]\big)^\frac{1}{2}}.
\end{equation}
By the orthogonality property we have
\begin{equation}
a^*[v,\Z_h w]=a^*[v,\Z_h u-u].
\end{equation}
Thus
\begin{equation}
a^*[v, \Z_h w] \leq \big(\lambda_{\max}(a)\big)^\frac{1}{2} \|\nabla
v\|_{L^2(\Omega)} \big(a^*[\Z_h u-u]\big)^\frac{1}{2}.
\end{equation}
Using the inequality
\begin{equation}\label{ddffsa3aq5132435}
\|\nabla v\|_{L^2(\Omega)}^2 \leq \eta_{\max} \ED_h[v].
\end{equation}
We deduce that
\begin{equation}\label{sghfgd98}
\ED_h[w]^\frac{1}{2} \leq \frac{1}{S^m} \big(\lambda_{\max}(a)
\eta_{\max}\big)^\frac{1}{2} \big(a^*[\Z_h u-u]\big)^\frac{1}{2}.
\end{equation}
It follows from \eref{ddff351ssad32435} that
\begin{equation}\label{sgshsg652}
 \|\nabla \I_h u- \nabla  u^m\|_{L^2(\Omega)}
\leq \frac{\eta_{\max}}{S^m} \big(\lambda_{\max}(a)
\big)^\frac{1}{2} \big(a^*[\Z_h u-u]\big)^\frac{1}{2}.
\end{equation}
And we deduce from Poincar\'{e} inequality that
\begin{equation}\label{hsgs51}
 \|\I_h u-  u^m\|_{L^2(\Omega)} \leq C_{\Omega} \frac{\eta_{\max}}{S^m} \big(\lambda_{\max}(a)
\big)^\frac{1}{2} \big(a^*[\Z_h u-u]\big)^\frac{1}{2}.
\end{equation}
We obtain theorem \ref{sfsfxadd653} by from equations
\eref{sgshsg652}, \eref{hsgs51}, lemma \ref{lem3} and theorem
\ref{th3}. Let us now prove theorem \ref{sfsfsdd653}. Using triangle
inequality we obtain
\begin{equation}
a\big[u- \J_h u^m\big]\leq a\big[u- \J_h u\big]+a\big[\J_h u- \J_h
u^m\big].
\end{equation}
The object $a\big[u- \J_h u\big]$ has already been bounded from
above by \eref{hsgs541}. Writing $w:=\J_h u- \J_h u^m$ we have
\begin{equation}
a[w]=\frac{a[\J_h w]}{\ED_h[w]} \ED_h[w].
\end{equation}
But $\ED_h[w]$  has already been estimated in equation
\eref{sghfgd98}. It remains to notice that
\begin{equation}
a[\J_h w]=Q[\hat{w}].
\end{equation}
From this point the arguments are similar to the ones employed in
subsection \ref{dcsswaqsaq21}, indeed
\begin{equation}
Q[\hat{w}] \leq \lambda_{\max}(Q) \|\nabla
\hat{w}\|_{L^2(\Omega)}\leq \eta^*_{\max} \lambda_{\max}(Q)
\ED_h[w].
\end{equation}

\subsection{Numerical homogenization from a transport point of view.} We assume
the mesh to be regular in the following sense: the nodes of the
 Vorno\"{i} diagram of $\T_h$ belong elements of the primal mesh. In
 dimension $2$ this means that for each triangle $K\in \T_h$ the
 intersection of the median of $K$ (the circumcenter) belongs to the interior of
 $K$.
 Let us write $Y_h$ the vector space spanned
by the functions $\chi_i$. For $v\in Z_h$ we define $\Y_h v$ by
\begin{equation}
\Y_h v:=\sum_{i\in \No_h} v_i \chi_i.
\end{equation}
The metric numerical homogenization method can be formulated in the
following way: look for $u^v \in Z_h$ (the space spanned by the
elements $\xi_i$) such that for all $i\in \No_h$,
\begin{equation}\label{gdhgdfctetfc5ee}
a^*[\chi_i, u^v]=(\chi_i,g)_{L^2{\Omega}}
\end{equation}
which implies the following finite volume orthogonality property for
all $i\in \No_h$,
\begin{equation}
a^*[\chi_i, u^v-u]=0.
\end{equation}
Equation \eref{gdhgdfctetfc5ee} can be written
\begin{equation}
\sum_{j\sim i} u^v_j \int_{\partial V_i} n.a.\nabla \xi_i
=\int_{V_i}g.
\end{equation}

Let us write $w:=\Z_h u-u^v$. From the equation \eref{jjhsjhsk726}
we obtain that
\begin{equation}\label{jjhsjdssdhsk726}
(\ED_h[w])^\frac{1}{2} \leq \frac{1}{\S^v} \sup_{v \in \Y_h}
\frac{a^*[v,w]}{(\ED_h[v])^\frac{1}{2}}.
\end{equation}
 Using the orthogonality property of
the finite volume method we obtain that for $v\in \Y_h$
\begin{equation}
a^*[v,w]=-\sum_{i\in \No_h} v_i \int_{\partial V_i} n.a \big(\nabla
\Z_h u- \nabla u\big).
\end{equation}
Writing $\ED^*_h$ the edges of the dual tessellation (edges of the
control volumes), we obtain that
\begin{equation}
a^*[v,w]=\sum_{e_{ij}\in \ED_h^*} (v_j-v_i) \int_{ e_{ij}} n_{ij}.a
\big(\nabla \Z_h u- \nabla u\big).
\end{equation}
Where $e_{ij}$ is the edge separating the control volume $V_i$ from
the control volume $V_j$ and  $n_{ij}$ is the unit vector orthogonal
to $e_{ij}$ pointing outside of $V_i$. It follows that
\begin{equation}
\begin{split}
a^*[v,w]\leq & \big(\ED_h[v]\big)^\frac{1}{2} \|\nabla_F \Z_h u
-\nabla_F u\|_{L^\infty(\ED_h^*)}\\& \Big(\sum_{e_{ij}\in \ED_h^*}
|e_{ij}|^2 \lambda_{\max}(a)
\lambda_{\max}(\sigma)\Big)^\frac{1}{2}.
\end{split}
\end{equation}
Now let us observe that
\begin{equation}
\sum_{e_{ij}\in \ED_h^*} |e_{ij}|^2 \leq  3 \eta_{\max}
\Vol(\Omega).
\end{equation}
We then have from equation \eref{jjhsjdssdhsk726}
\begin{equation}\label{dggcfgfdc53}
\begin{split}
\big(\ED_h[w]\big)^\frac{1}{2} \leq & \frac{3}{\S^v}\|\nabla_F \Z_h
u -\nabla_F u\|_{L^\infty(\ED_h^*)} \\& \eta_{\max} \Vol(\Omega)
\big(\lambda_{\max}(a) \lambda_{\max}(\sigma)\big)^\frac{1}{2}.
\end{split}
\end{equation}
Equation \eref{sqhdgdgz7112sss} of theorem \ref{sfsfsdd653sss} is
then a straightforward consequence of equation
\eref{ddff351ssad32435} and lemma \ref{dgdssgqaswq536f}. Let us now
prove equation \eref{sqhdgdgz7sss112} of theorem
\ref{sfsfsdd653sss}. By triangle inequality
\begin{equation}
a\big[u- \J_h u^v\big]\leq a\big[u- \J_h u\big]+a\big[\J_h u- \J_h
u^v\big].
\end{equation}
$a\big[u- \J_h u\big]$ has already been estimated in equation
\eref{hsgs541}. Writing $w:=\J_h u- \J_h u^v$ we have
\begin{equation}
a[w]=\frac{a^*[\J_h w]}{\ED_h[w]} \ED_h[w].
\end{equation}
But $\ED_h[w]$  has already been estimated in equation
\ref{dggcfgfdc53}. It remains to notice that $\frac{a^*[\J_h
w]}{\ED_h[w]}$ has already been estimated in equations \eref{eqq1}
and \eref{eqq2} to conclude the proof.

\section{Numerical Experiments}
\begin{figure}[httb]
\begin{center}
\includegraphics[%
  scale=0.3]{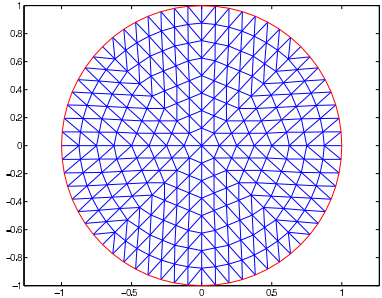}
\caption{\label{cap:computational-mesh}Coarse Grid}
\end{center}
\end{figure}

Let us now give illustrations of the implementation of this method.
The domain is the unit disk in dimension two. Equation \eref{ghjh52}
is solved on a fine tessellation characterized by $66049$ nodes and
$131072$ triangles. The coarse tessellation has $289$ nodes and
$512$ triangles (figure \ref{cap:computational-mesh}). It is
important to recall that since our methods involve the computation
of global harmonic coordinates, the memory requirement and CPU time
are not improved if one needs to solve \eref{ghjh52} only one time,
whereas localized methods such as the one of Hou and Wu or E. and
Engquist do improve the memory requirement or the CPU time.

The elliptic operator associated to equation \eref{ghjh52} has been
up-scaled to an operator defined on the coarse mesh (compression by
a factor of $300$) using $5$ different methods:
\begin{itemize}
\item The Galerkin scheme described in subsection \ref{gfssg5} using the multiscale finite element
$\psi_i$ noted FEM\_$\psi$.
\item The Galerkin scheme described in subsection \ref{hgdjhdg541} using the localized elements
$\xi_i$ noted FEM\_$\xi$.
\item The metric based compression scheme described in subsection
\ref{gfgxfs541} noted MBFEM.
\item The finite volume method described in subsection
\ref{hsghsg6531} noted FVM.
\item A multi-scale finite element method noted LFEM where $F$ is
computed locally\footnote{instead of globally} on each triangle $K$
of the coarse mesh as the solution of a cell problem with boundary
condition $F(x)=x$ on $\partial K$. This method has been implemented
in order to understand the effect of the removal of global
information in the structure of the metric induced by $F$.
\end{itemize}

\begin{figure}[!]
  \begin{center}
    \subfigure[$a$.]
    {\includegraphics[width=0.35\textwidth,height= 0.3\textwidth]{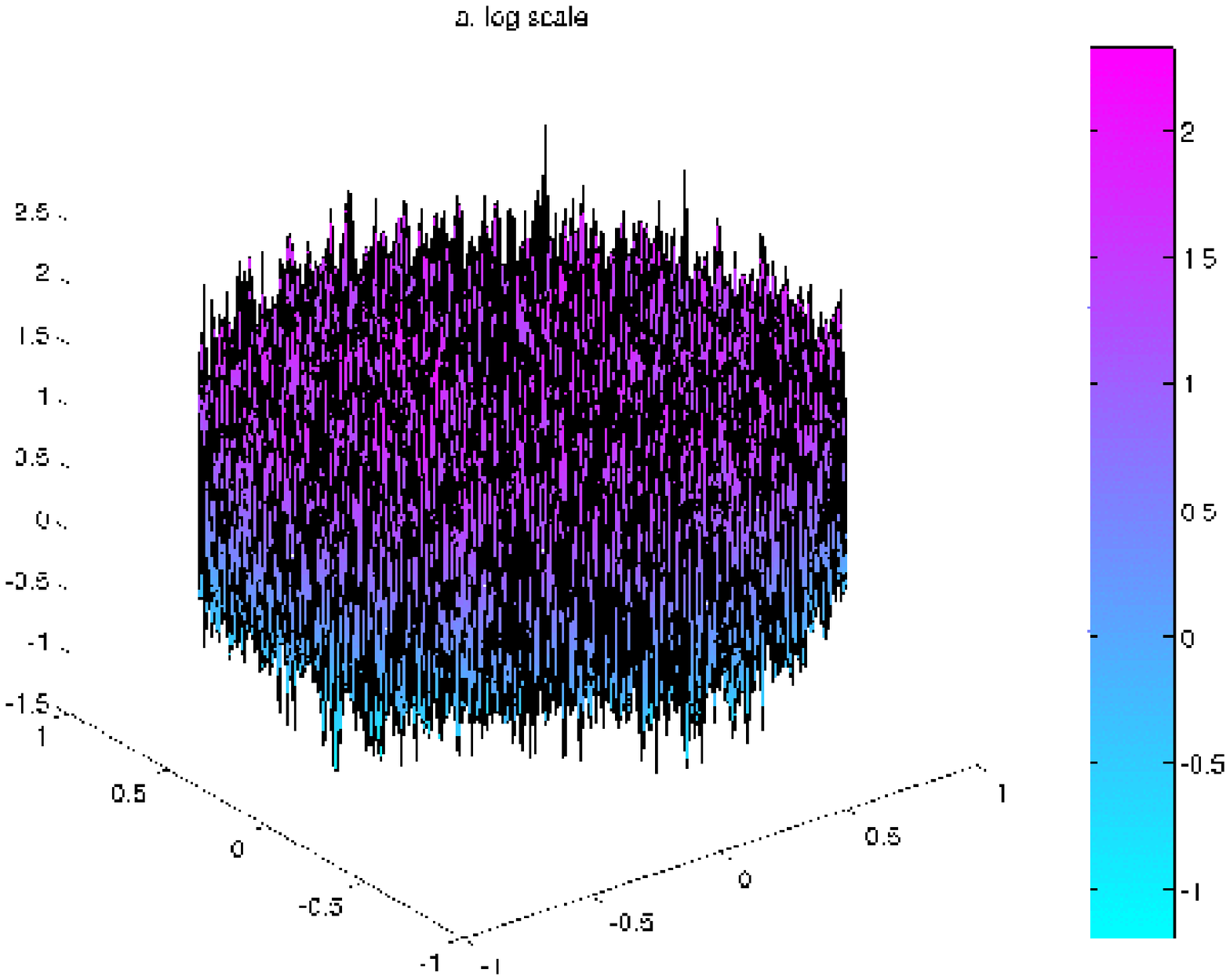}}
    \goodgap
    \subfigure[$\T^F$]
    {\includegraphics[width=0.35\textwidth,height= 0.3\textwidth]{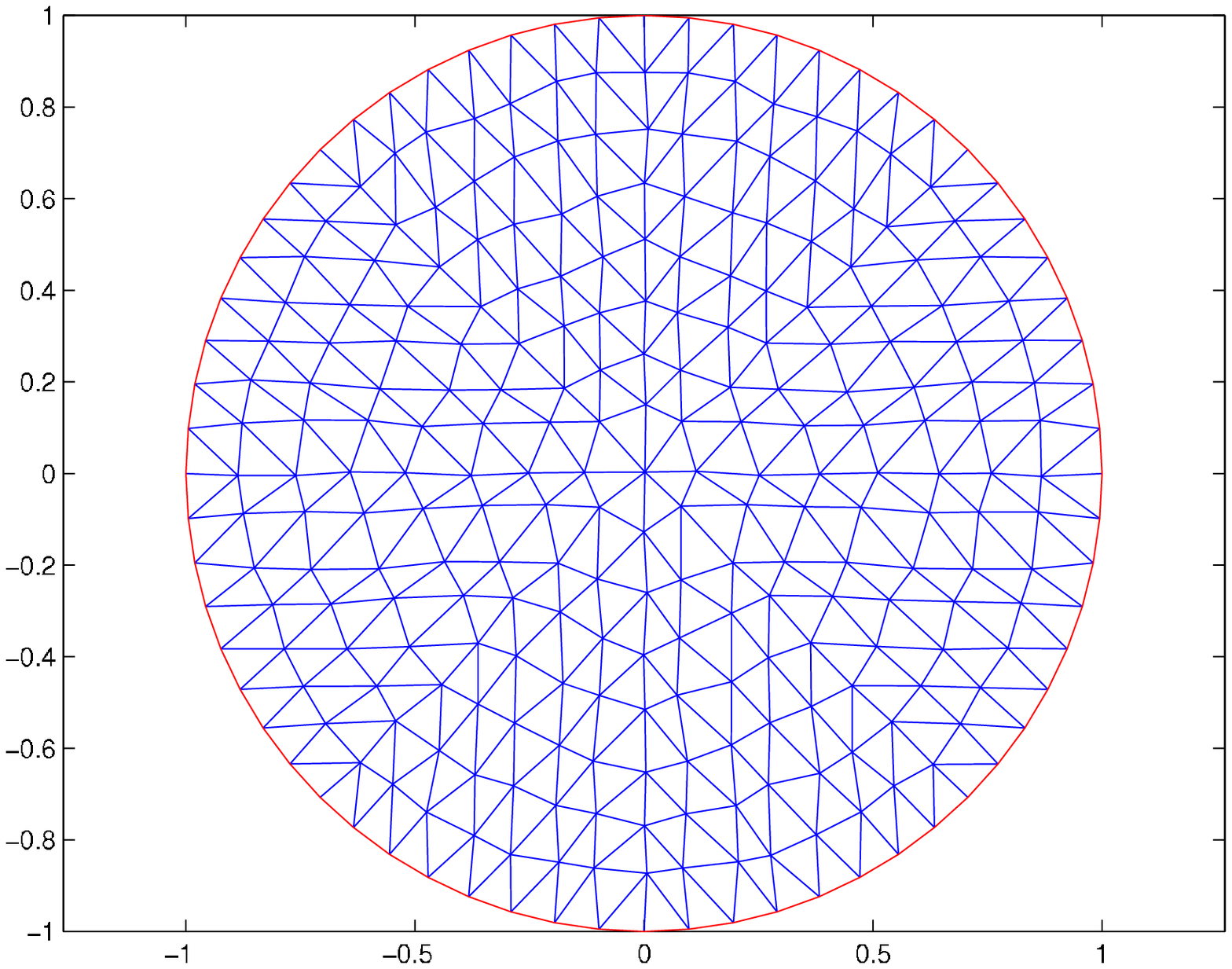}}\\
    \caption{\label{cap:Example Trigonometric}Example \ref{exa:trigonometric}, Trigonometric multi-scale.}
\end{center}
\end{figure}
\begin{figure}[httb]
  \begin{center}
    \subfigure[Condition Number.\label{condnumb1}]
    {\includegraphics[width=0.35\textwidth,height= 0.3\textwidth]{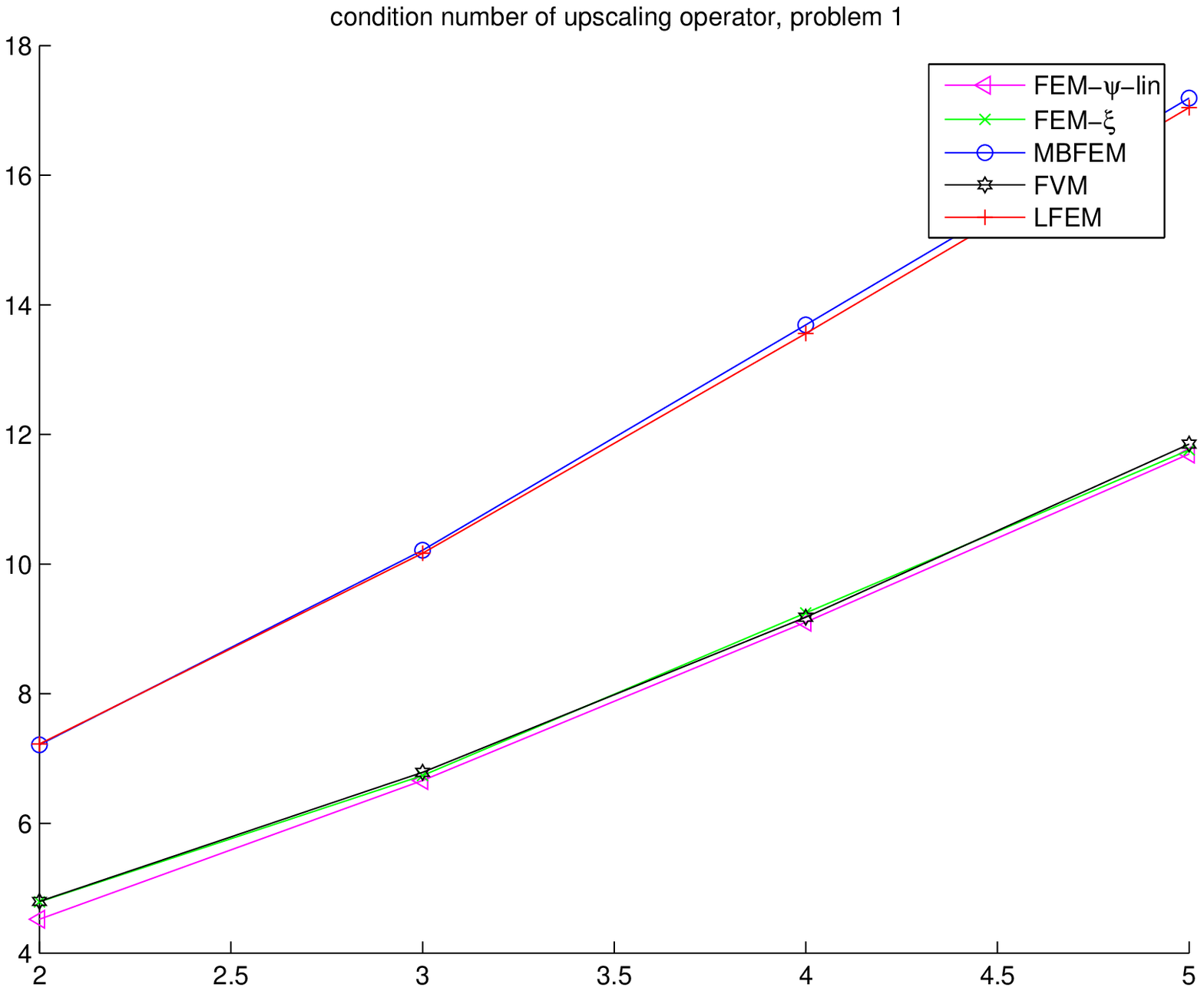}}
    \goodgap
    \subfigure[Coarse Mesh $L^1$ error.\label{l1err1}]
    {\includegraphics[width=0.35\textwidth,height= 0.3\textwidth]{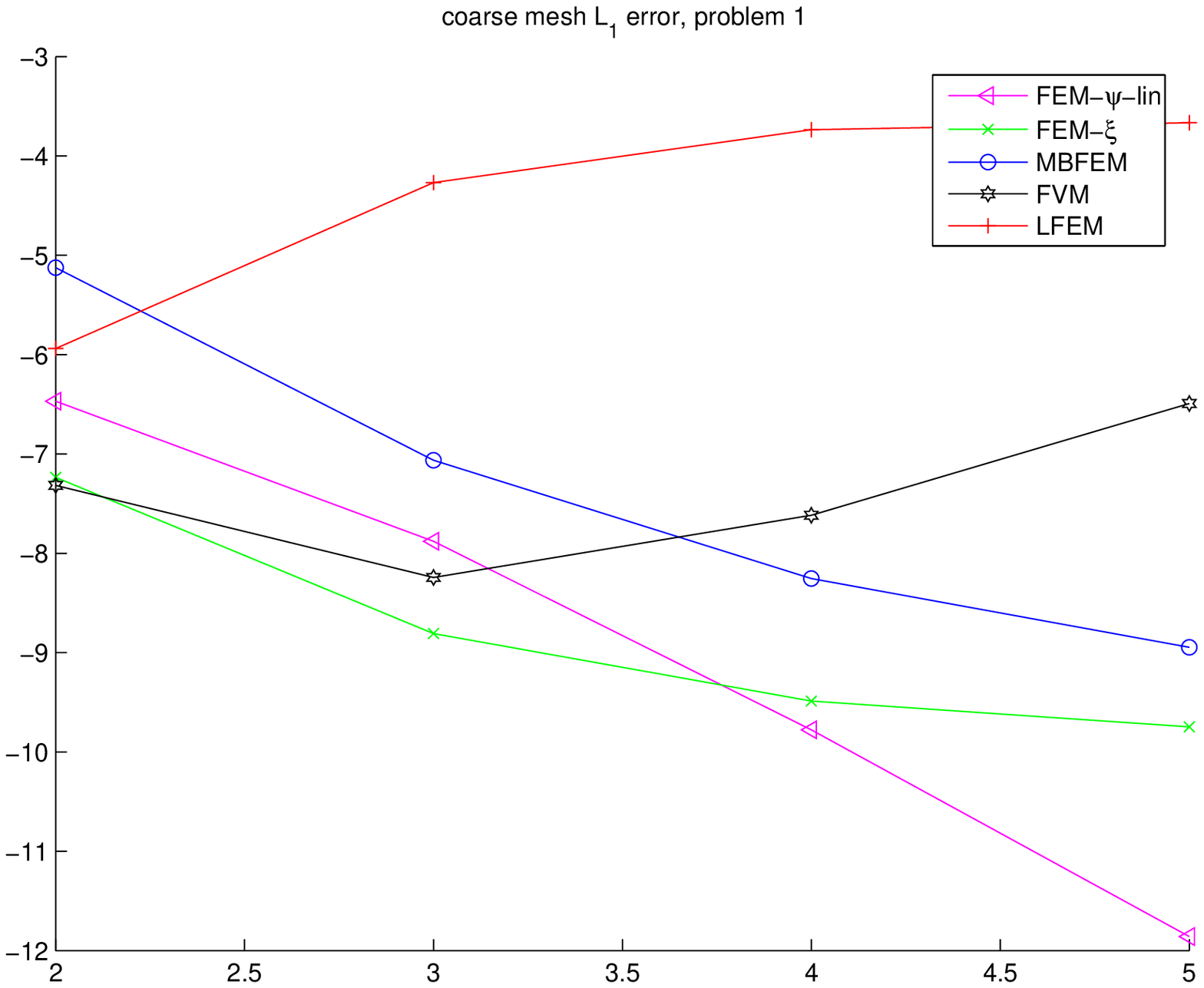}}\\
    \caption{\label{cap:Example Trigonometric2}Example \ref{exa:trigonometric}, Trigonometric Multiscale }
\end{center}
\end{figure}

\begin{example}
\label{exa:trigonometric}Trigonometric multiscale
\end{example}
$a(x)=\frac{1}{6}(\frac{1.1+\sin(2\pi x/\epsilon_{1})}{1.1+\sin(2\pi
y/\epsilon_{1})}+\frac{1.1+\sin(2\pi y/\epsilon_{2})}{1.1+\cos(2\pi
x/\epsilon_{2})}+\frac{1.1+\cos(2\pi x/\epsilon_{3})}{1.1+\sin(2\pi
y/\epsilon_{3})}+\frac{1.1+\sin(2\pi y/\epsilon_{4})}{1.1+\cos(2\pi
x/\epsilon_{4})}+\frac{1.1+\cos(2\pi x/\epsilon_{5})}{1.1+\sin(2\pi
y/\epsilon_{5})}+\sin(4x^{2}y^{2})+1)$,
 where
$\epsilon_{1}=\frac{1}{5}$,$\epsilon_{2}=\frac{1}{13}$,$\epsilon_{3}=\frac{1}{17}$,$\epsilon_{4}=\frac{1}{31}$,$\epsilon_{5}=\frac{1}{65}$.

Figure \ref{cap:Example Trigonometric} is an illustration of $\T^F$
the deformation of the coarse mesh (figure
\ref{cap:computational-mesh}) under the metric induced by $F$. The
deformation is small since the medium is quasi-periodic. The weak
aspect ratio for triangles the coarse mesh in the metric induced by
$F$ is $\eta^{*}_{\min}=1.1252$. Table \ref{caswError1}
 gives the relative error estimated on the nodes of the coarse mesh between the solution $u$ of the initial PDE
 \eref{ghjh52} and an approximation obtained from the up-scaled
 operator on the nodes of the coarse mesh. Table \ref{caswErrswaor1}
 gives
 the relative error estimated on the nodes of the fine mesh between  $u$ and the  $\J_h$-interpolation of the previous approximations with respect
 to $F$ on a fine resolution. Figure \ref{condnumb1} gives the condition number of the
stiffness matrix associated to the up-scaled operator versus
$-\log_2 h$ (logarithm of the resolution). Figure \ref{l1err1} gives
the relative $L_1$-distance between $u$ and its approximation on the
coarse mesh in log scale versus $-\log_2 h$ (logarithm of the
resolution). Observe that for the method LFEM this error increase
with the resolution this is an effect of the so called cell
resonance observed in \cite{MR2119937} and \cite{AlBr04}. This cell
resonance does not occur with the methods proposed in this paper.
The finite volume method is characterized by the the best stability
and one of the best accuracy at a coarse resolution. The increase in
the error observed for this method as the resolution is decreased is
a numerical artifact created by the fine mesh: one has to divide the
coarse tessellation into coarse control volumes. These coarse
control volumes are unions of the control volumes defined on a fine
mesh and when the refinement between the coarse and the fine mesh is
small and the  triangulation irregular it is not possible to divide
the coarse tessellation into control volumes intersecting the edges
of the primal mesh close to the midpoints of those edges and the
other control volumes close to the barycenters of the coarse
triangles.
\begin{table}[!]
\begin{center}
\begin{tabular}{|c||c|c|c|c|c|}

\hline Coarse Mesh Error&  FEM\_$\psi$& FEM\_$\xi$& MBFEM& FVM&
LFEM\tabularnewline \hline \hline
\begin{tabular}{c}
$L^{1}$\tabularnewline \hline $L^{2}$\tabularnewline \hline
$L^{\infty}$\tabularnewline \hline $H^{1}$\tabularnewline
\end{tabular}&
\begin{tabular}{c}
0.0042\tabularnewline \hline 0.0039\tabularnewline \hline
0.0059\tabularnewline \hline 0.0060\tabularnewline
\end{tabular}&
\begin{tabular}{c}
0.0022\tabularnewline \hline 0.0024\tabularnewline \hline
0.0090\tabularnewline \hline 0.0262\tabularnewline
\end{tabular}&
\begin{tabular}{c}
0.0075\tabularnewline \hline 0.0074\tabularnewline \hline
0.0154\tabularnewline \hline 0.0568\tabularnewline
\end{tabular}&
\begin{tabular}{c}
0.0032\tabularnewline \hline 0.0040\tabularnewline \hline
0.0117\tabularnewline \hline 0.0203\tabularnewline
\end{tabular}&
\begin{tabular}{c}
0.0411\tabularnewline \hline 0.0441\tabularnewline \hline
0.0496\tabularnewline \hline 0.0763\tabularnewline
\end{tabular}\tabularnewline
\hline
\end{tabular}
\caption{\label{caswError1}Example \ref{exa:trigonometric},
Trigonometric multi-scale.}
\end{center}
\end{table}

\begin{table}[!]
\begin{center}
\begin{tabular}{|c||c|c|c|c|c|}
\hline Fine mesh Error& FEM\_$\psi$& FEM\_$\xi$& MBFEM& FVM&
LFEM\tabularnewline \hline \hline
\begin{tabular}{c}
$L^{1}$\tabularnewline \hline $L^{2}$\tabularnewline \hline
$L^{\infty}$\tabularnewline \hline $H^{1}$\tabularnewline
\end{tabular}&
\begin{tabular}{c}
0.0042\tabularnewline \hline 0.0043\tabularnewline \hline
0.0063\tabularnewline \hline 0.0581\tabularnewline
\end{tabular}&
\begin{tabular}{c}
0.0085\tabularnewline \hline 0.0082\tabularnewline \hline
0.0112\tabularnewline \hline 0.0540\tabularnewline
\end{tabular}&
\begin{tabular}{c}
0.0053\tabularnewline \hline 0.0061\tabularnewline \hline
0.0154\tabularnewline \hline 0.0778\tabularnewline
\end{tabular}&
\begin{tabular}{c}
0.0080\tabularnewline \hline 0.0078\tabularnewline \hline
0.0141\tabularnewline \hline 0.0601\tabularnewline
\end{tabular}&
\begin{tabular}{c}
0.0593\tabularnewline \hline 0.0591\tabularnewline \hline
0.0597\tabularnewline \hline 0.0943\tabularnewline
\end{tabular}\tabularnewline
\hline
\end{tabular}
\caption{\label{caswErrswaor1}Example \ref{exa:trigonometric},
Trigononmetric Multiscale }
\end{center}
\end{table}

\clearpage

\begin{figure}[httb]
  \begin{center}
    \subfigure[Condition Number. \label{condnumb2}]
    {\includegraphics[width=0.35\textwidth,height= 0.3\textwidth]{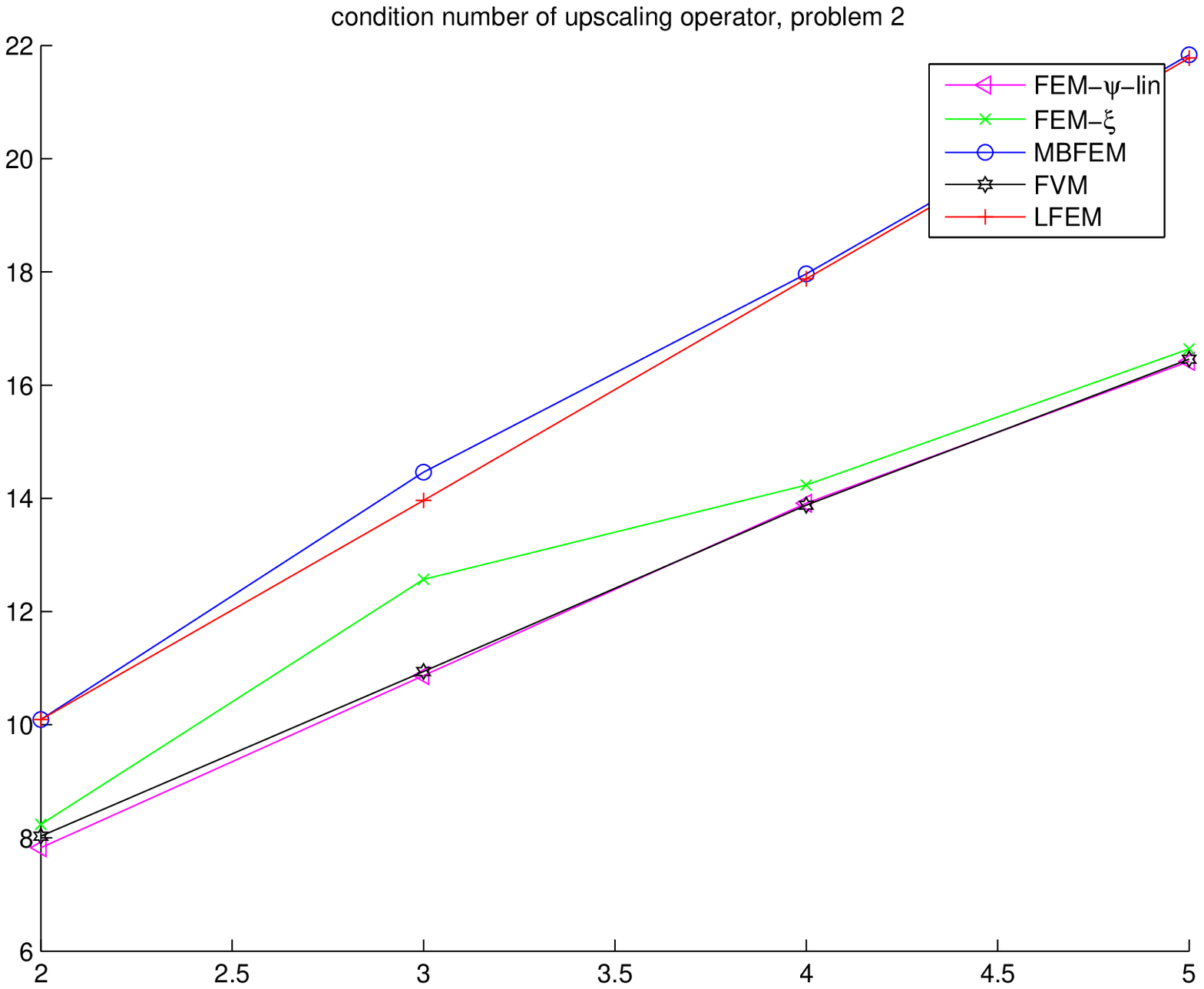}}
    \goodgap
    \subfigure[Coarse Mesh $L^1$ error. \label{l1err2}]
    {\includegraphics[width=0.35\textwidth,height= 0.3\textwidth]{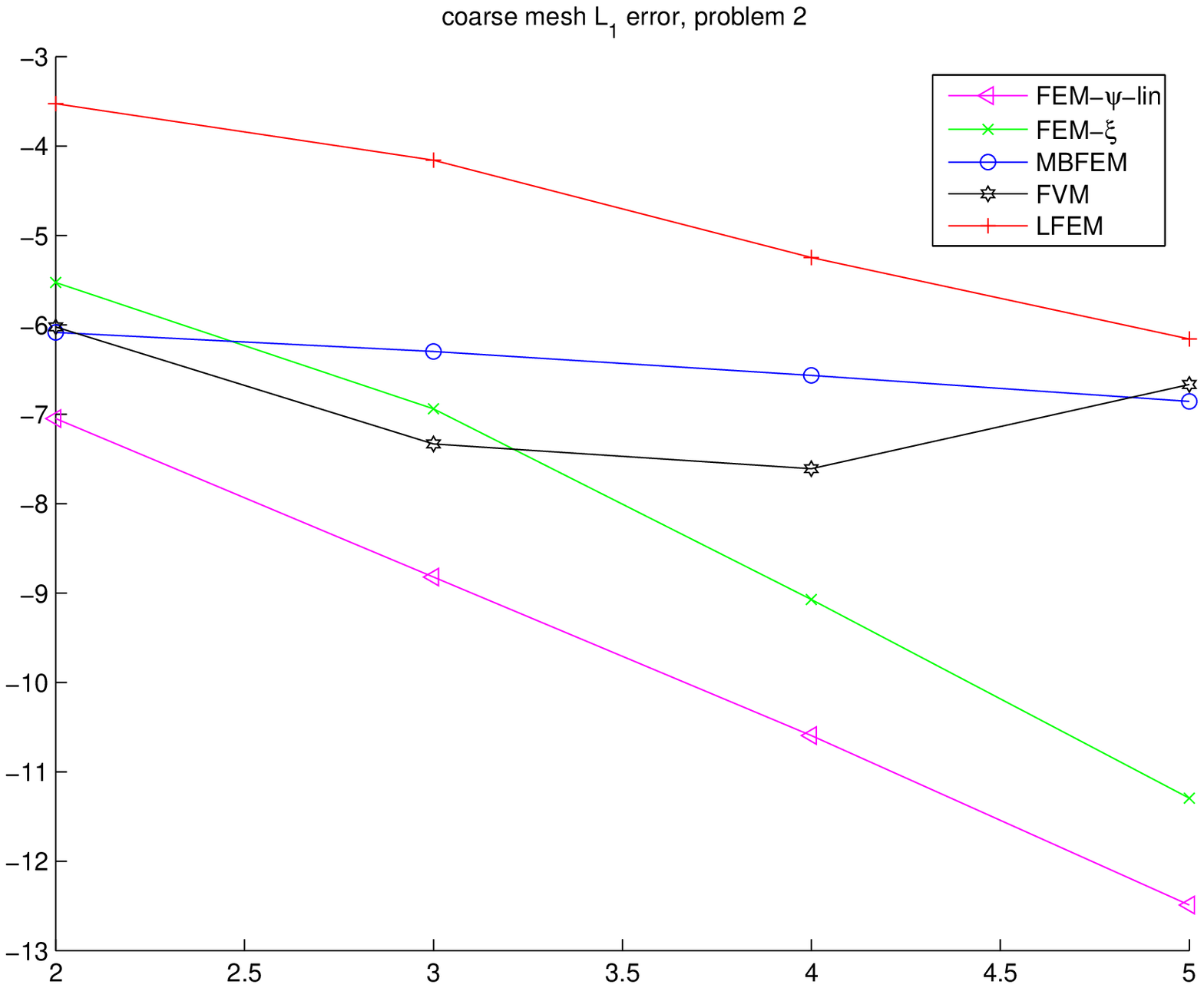}}\\
    \caption{Example \ref{exa:channel}, High Conductivity Channel }
\end{center}
\end{figure}

\begin{example}
\label{exa:channel}High Conductivity channel
\end{example}
In this example $a$ is random and characterized by a fine and long
ranged high conductivity channel. We choose $a(x)=100$, if $x$ is in
the channel, and $a(x)=O(1)$, if $x$ is not in the channel.  The
weak aspect ratio for triangles the coarse mesh in the metric
 induced by $F$ is  $\eta^{\star}_{\min}=2.2630$. Table \ref{caswErxswr122}
 gives the relative error estimated on the nodes of the coarse mesh between the solution $u$ of the initial PDE
 \eref{ghjh52} and an approximation obtained from the up-scaled
 operator on the nodes of the coarse mesh.
 Table
 \ref{caswErxswr1} gives
 the relative error estimated on the nodes of the fine mesh between  $u$ and the  $\J_h$-interpolation of the previous approximations with respect
 to $F$ on a fine resolution. Figure \ref{condnumb2} gives the condition number of the
stiffness matrix associated to the up-scaled operator versus
$-\log_2 h$ (logarithm of the resolution). Figure \ref{l1err2} gives
the relative $L_1$-distance between $u$ and its approximation on the
coarse mesh in log scale versus $-\log_2 h$.

\begin{figure}[!]
  \begin{center}
    \subfigure[$a$.]
    {\includegraphics[width=0.35\textwidth,height= 0.3\textwidth]{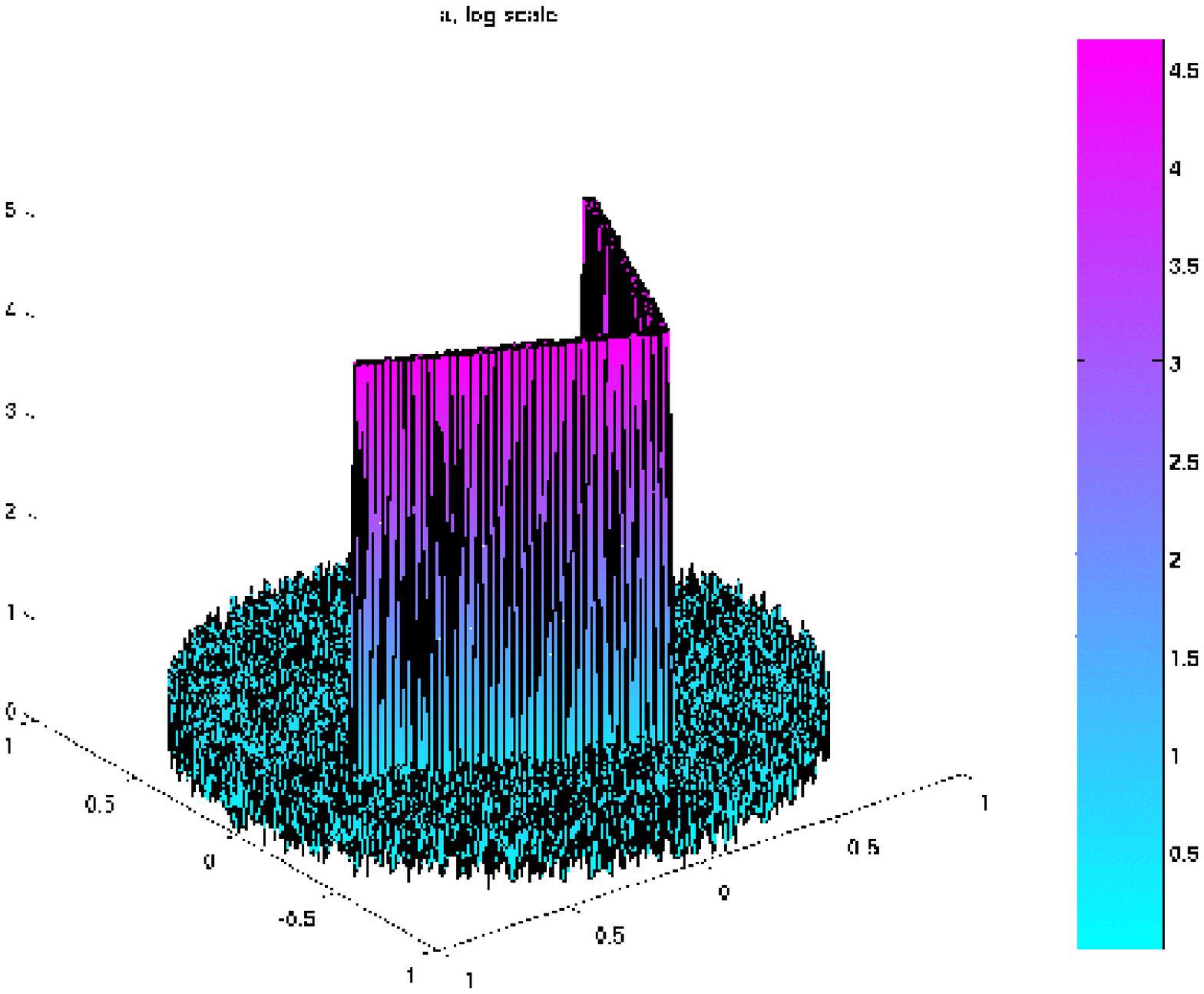}}
    \goodgap
    \subfigure[$\T^F$]
    {\includegraphics[width=0.35\textwidth,height= 0.3\textwidth]{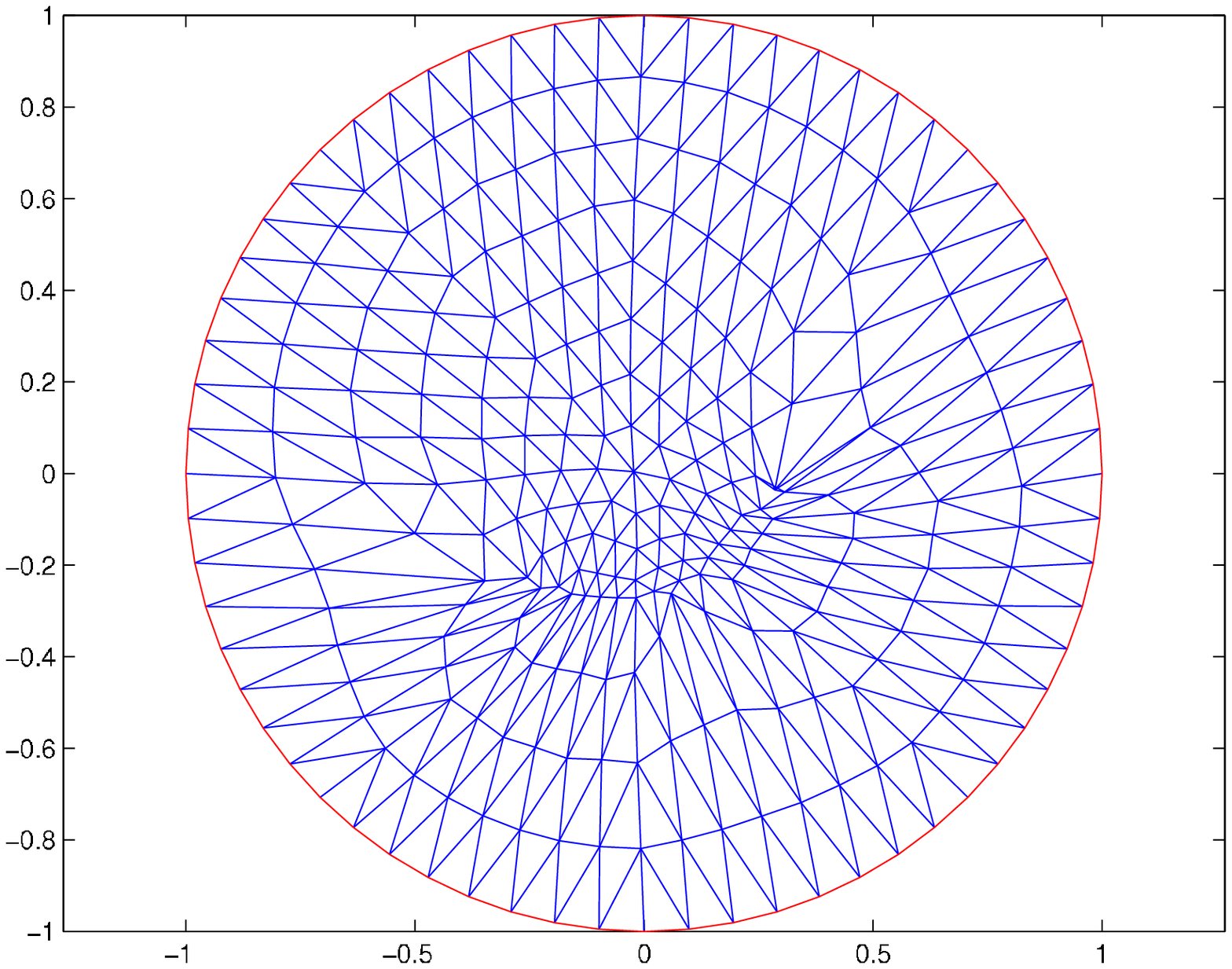}}\\
    \caption{Example \ref{exa:channel}, High Conductivity Channel. \label{figexa:channel}}
\end{center}
\end{figure}

 Observe in figure \ref{figexa:channel}
that the effect of the new metric on the mesh is to bring close
together nodes linked by a path of low electrical resistance.
\begin{Remark}\label{hgshsh512}
Let us recall that the
 natural distance associated to the Laplace operator on a fractal
 space is also the so called resistance metric \cite{MR1317710},
 \cite{MR1962353}, \cite{BCK05}. It is thus natural to find that a
 similar (not equivalent)
 notion
 of distance allows the numerical homogenization PDEs with arbitrary
 coefficients. More precisely the analogue of the resistance metric
 here are the harmonic mappings. The analysis of these mappings
allows to bypass boundary layers effects in homogenization in
periodic media \cite{AlBr04}, it allows to obtain quantitative
estimates on the heat kernel of periodic operators \cite{Ow00a} or
to analyze PDEs characterized by an infinite number of non separated
scales \cite{BeOw00b}, \cite{Ow04}.
\end{Remark}

\begin{table}[!]
\begin{center}
\begin{tabular}{|c||c|c|c|c|c|}
\hline Coarse Mesh Error& FEM\_$\psi$& FEM\_$\xi$& MBFEM& FVM&
LFEM\tabularnewline \hline \hline
\begin{tabular}{c}
$L^{1}$\tabularnewline \hline $L^{2}$\tabularnewline \hline
$L^{\infty}$\tabularnewline \hline $H^{1}$\tabularnewline
\end{tabular}&
\begin{tabular}{c}
0.0022\tabularnewline \hline 0.0025\tabularnewline \hline
0.0120\tabularnewline \hline 0.0120\tabularnewline
\end{tabular}&
\begin{tabular}{c}
0.0081\tabularnewline \hline 0.0096\tabularnewline \hline
0.0227\tabularnewline \hline 0.0384\tabularnewline
\end{tabular}&
\begin{tabular}{c}
0.0127\tabularnewline \hline 0.0179\tabularnewline \hline
0.0549\tabularnewline \hline 0.0919\tabularnewline
\end{tabular}&
\begin{tabular}{c}
0.0062\tabularnewline \hline 0.0081\tabularnewline \hline
0.0174\tabularnewline \hline 0.0265\tabularnewline
\end{tabular}&
\begin{tabular}{c}
0.0519\tabularnewline \hline 0.0606\tabularnewline \hline
0.1223\tabularnewline \hline 0.1514\tabularnewline
\end{tabular}\tabularnewline
\hline
\end{tabular}
\caption{\label{caswErxswr122}Example \ref{exa:channel}, High
Conductivity Channel.}
\end{center}
\end{table}

\begin{table}[!]
\begin{center}
\begin{tabular}{|c||c|c|c|c|c|}
\hline Fine mesh Error& FEM\_$\psi$& FEM\_$\xi$& MBFEM& FVM&
LFEM\tabularnewline \hline \hline
\begin{tabular}{c}
$L^{1}$\tabularnewline \hline $L^{2}$\tabularnewline \hline
$L^{\infty}$\tabularnewline \hline $H^{1}$\tabularnewline
\end{tabular}&
\begin{tabular}{c}
0.0070\tabularnewline \hline 0.0069\tabularnewline \hline
0.0133\tabularnewline \hline 0.0760\tabularnewline
\end{tabular}&
\begin{tabular}{c}
0.0155\tabularnewline \hline 0.0153\tabularnewline \hline
0.0227\tabularnewline \hline 0.1032\tabularnewline
\end{tabular}&
\begin{tabular}{c}
0.0164\tabularnewline \hline 0.0202\tabularnewline \hline
0.0573\tabularnewline \hline 0.1838\tabularnewline
\end{tabular}&
\begin{tabular}{c}
0.0121\tabularnewline \hline 0.0123\tabularnewline \hline
0.0214\tabularnewline \hline 0.0820\tabularnewline
\end{tabular}&
\begin{tabular}{c}
0.0612\tabularnewline \hline 0.0743\tabularnewline \hline
0.1226\tabularnewline \hline 0.2142\tabularnewline
\end{tabular}\tabularnewline
\hline
\end{tabular}
\caption{\label{caswErxswr1}Example \ref{exa:channel}, High
Conductivity Channel.}
\end{center}
\end{table}

\clearpage

\begin{figure}[!]
  \begin{center}
    \subfigure[$a$.]
    {\includegraphics[width=0.35\textwidth,height= 0.3\textwidth]{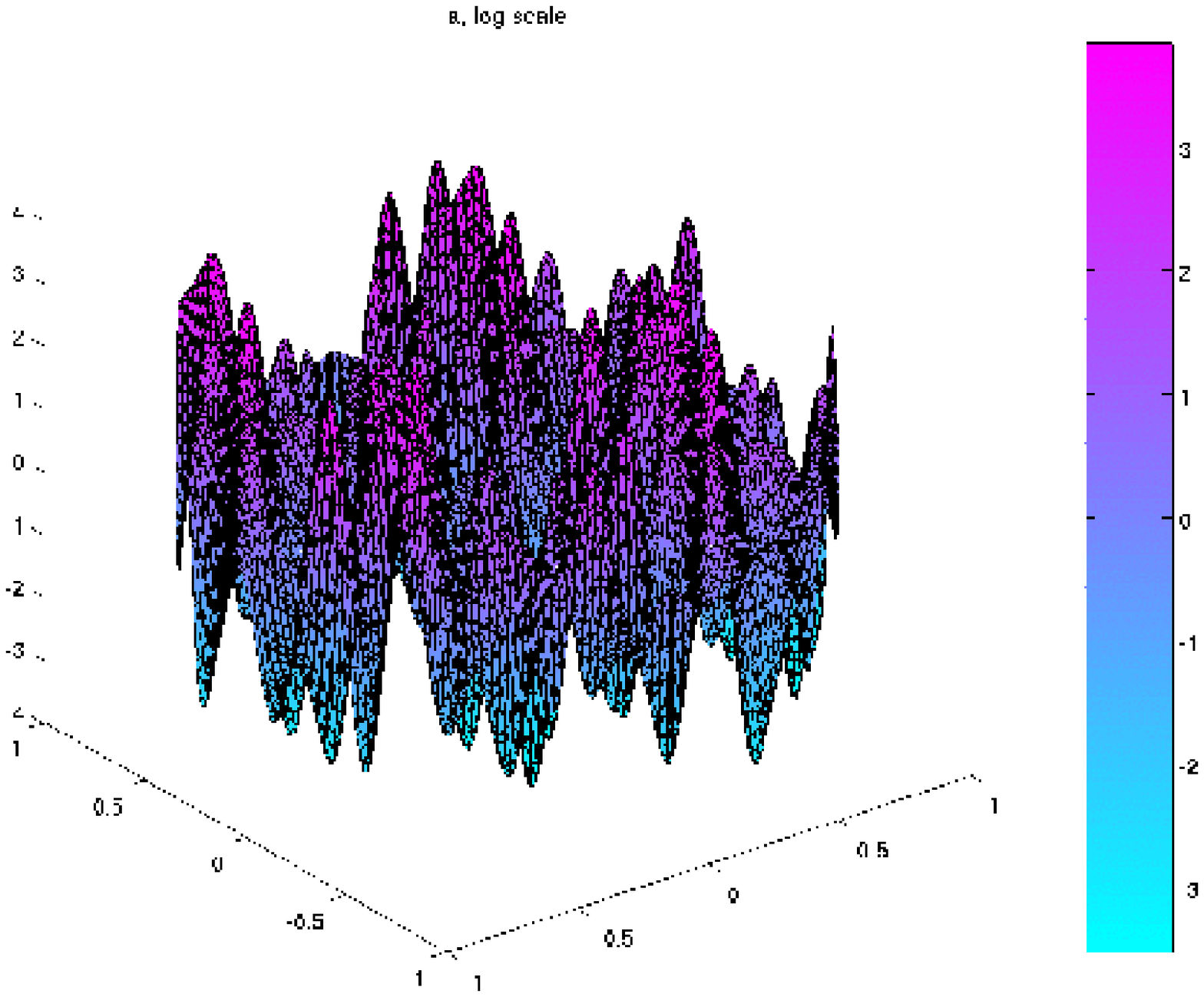}}
    \goodgap
    \subfigure[$\T^F$]
    {\includegraphics[width=0.35\textwidth,height= 0.3\textwidth]{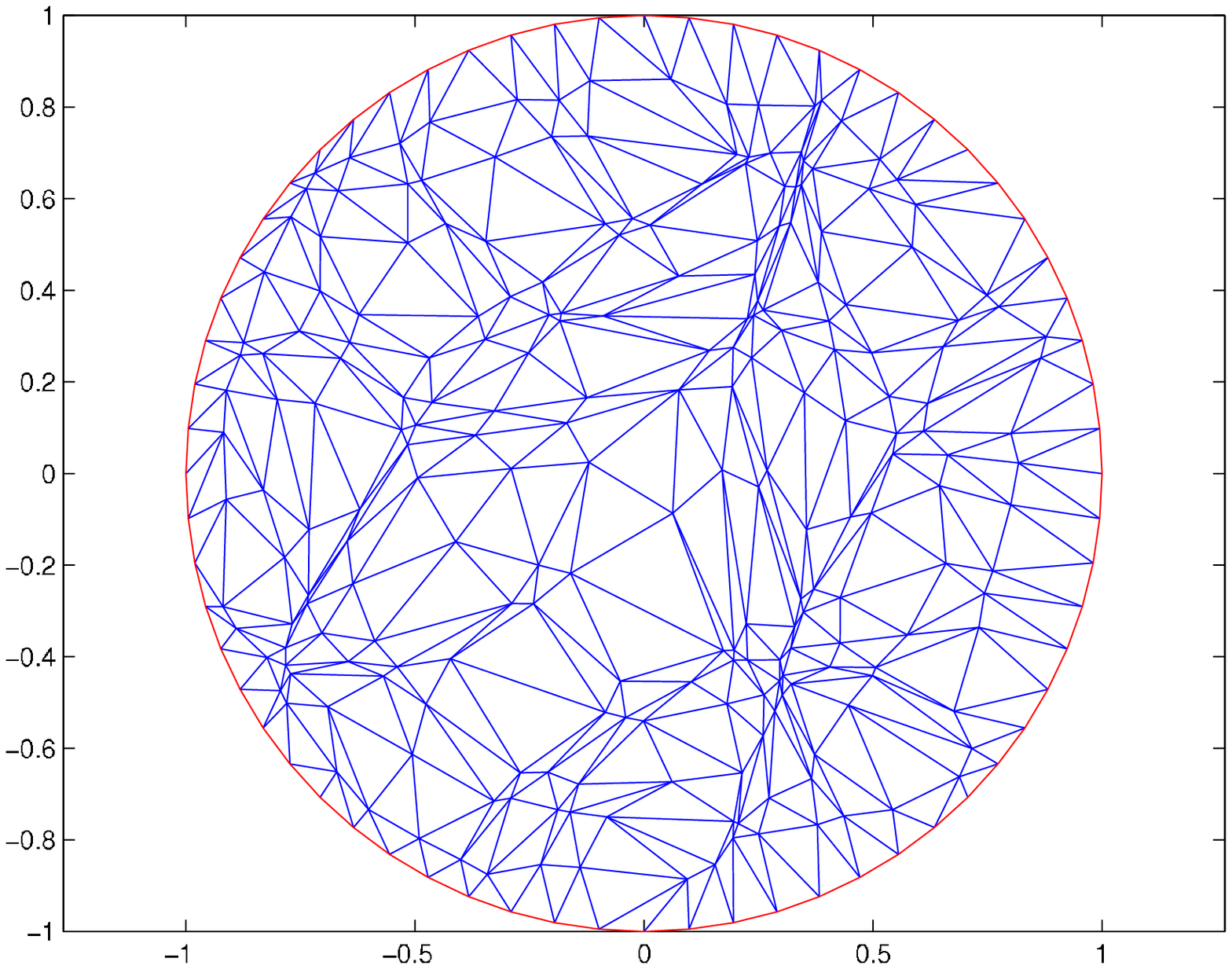}}\\
    \caption{\label{cap:ExampleFourier}Example \ref{exa:Fourier}, Random Fourier Modes .}
\end{center}
\end{figure}

\begin{figure}[!]
  \begin{center}
    \subfigure[$L^1$ error\label{Coarse-mesh-error-history-L1}.]
    {\includegraphics[width=0.35\textwidth,height= 0.3\textwidth]{circleg_4_errL1c}}
    \goodgap
    \subfigure[$L^2$ error\label{Coarse-mesh-error-history-L2}.]
    {\includegraphics[width=0.35\textwidth,height= 0.3\textwidth]{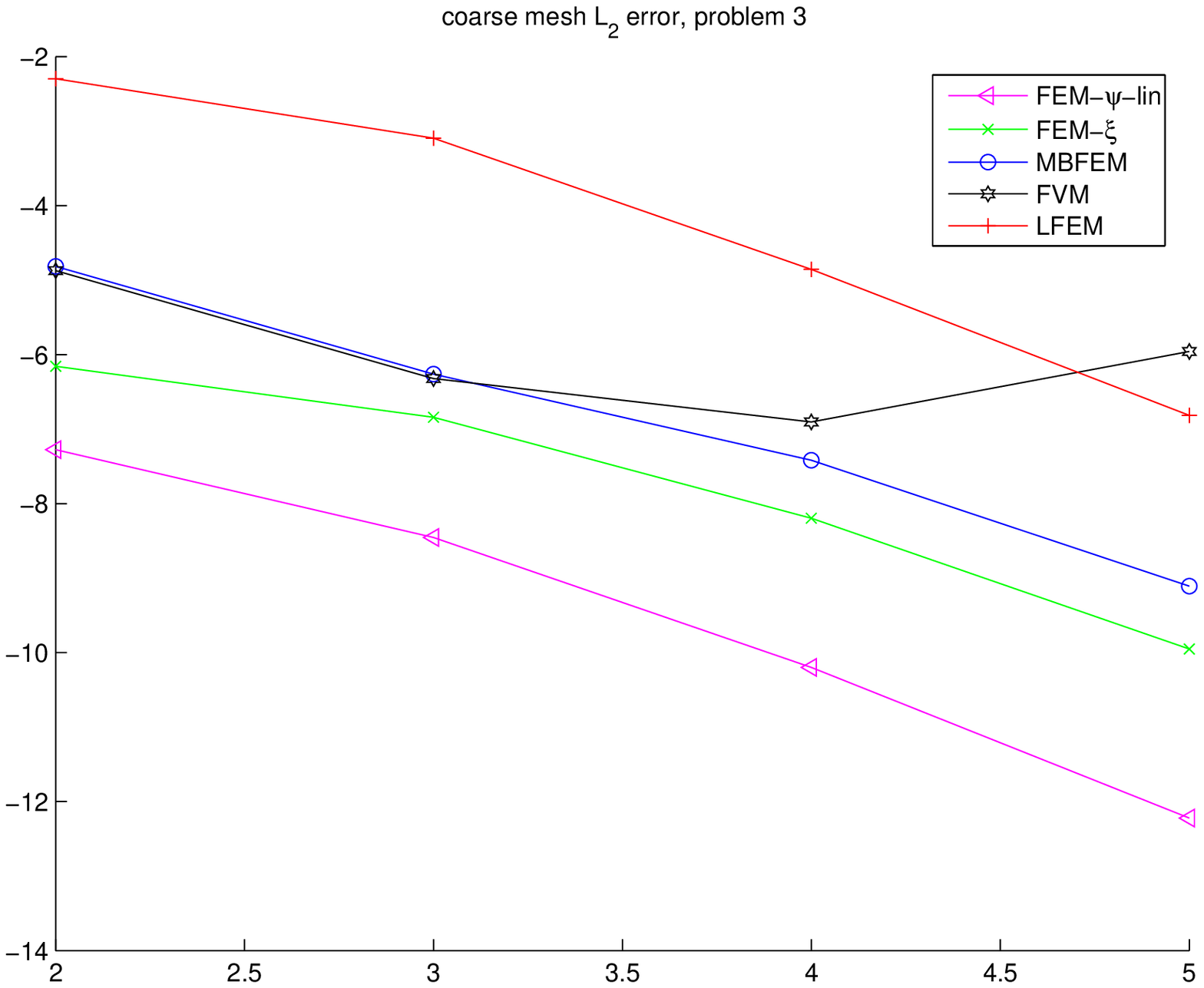}}\\
    \subfigure[$L^{\infty}$ error\label{Coarse-mesh-error-history-Li}.]
    {\includegraphics[width=0.35\textwidth,height= 0.3\textwidth]{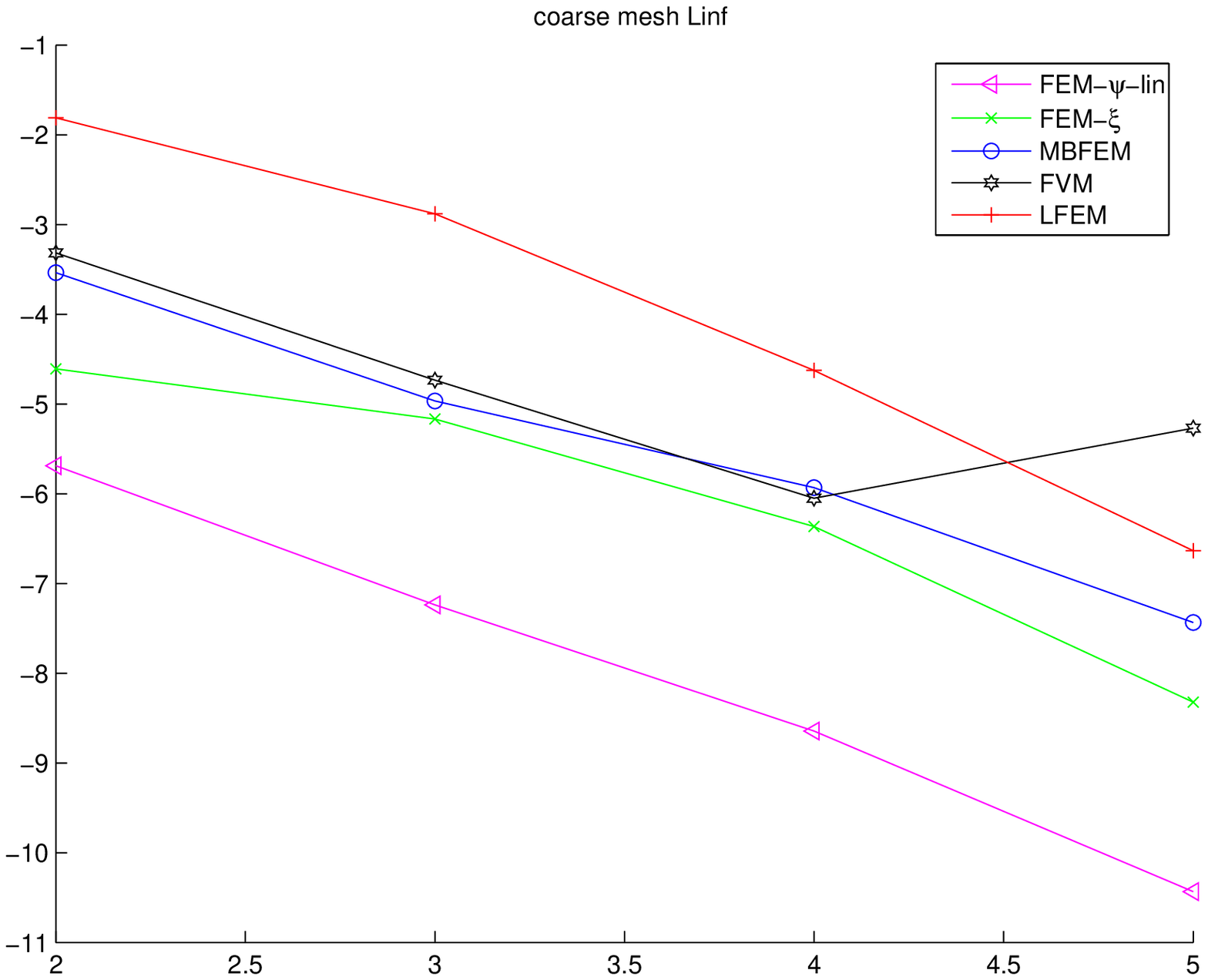}}
    \goodgap
    \subfigure[$H^1$ error\label{Coarse-mesh-error-history-H1}.]
    {\includegraphics[width=0.35\textwidth,height= 0.3\textwidth]{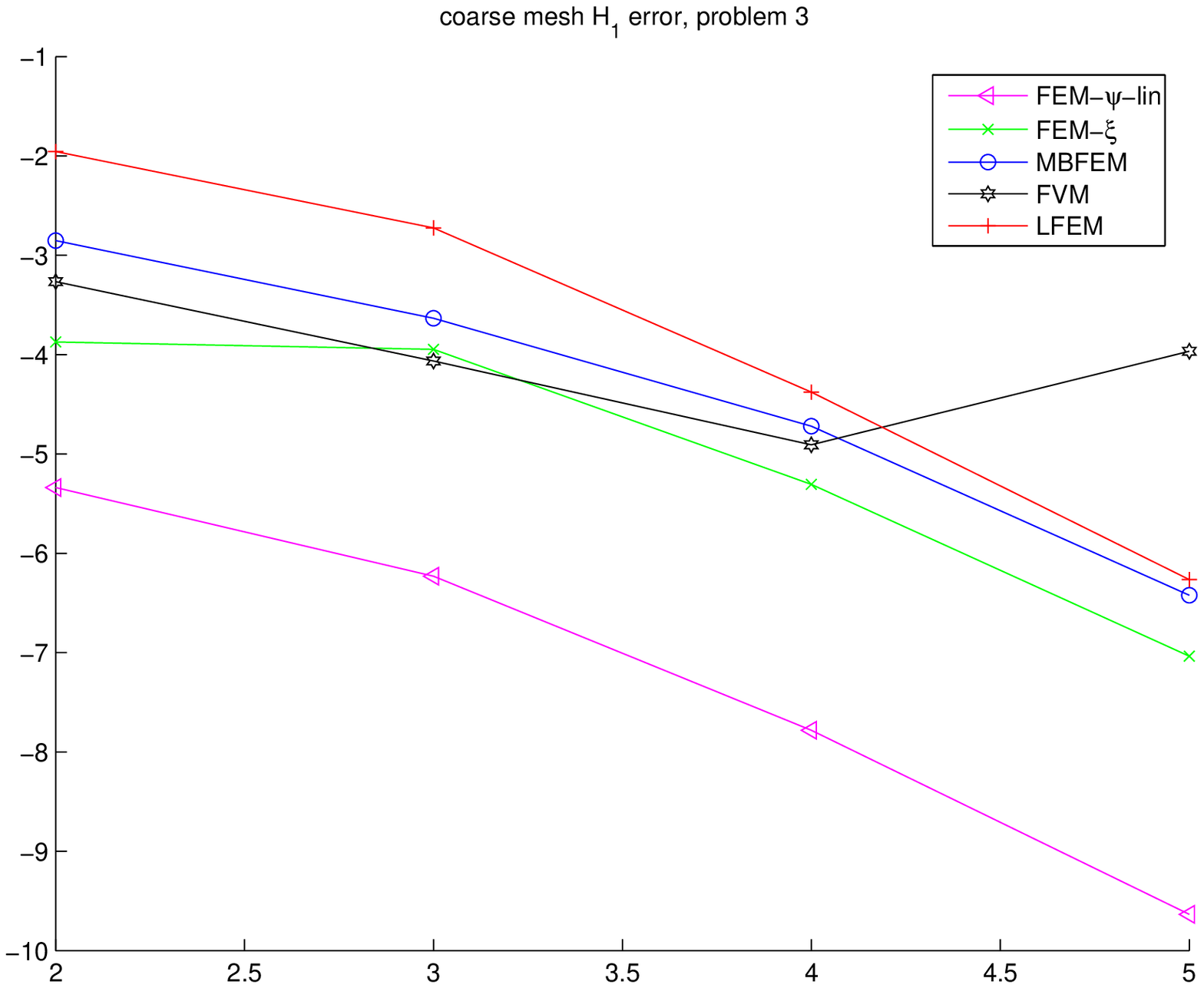}}\\[-10pt]
    \caption{\label{cap:Coarse-mesh-error-history}Coarse mesh error
      ($\log_{2}$) $L^{1}$, $L^{2}$, $L^{\infty}$ and $H^{1}$ errors
      v.s coarse mesh refinement, Example \ref{exa:Fourier}, Random Fourier Modes.}
\end{center}
\end{figure}

\begin{figure}[!]
  \begin{center}
    \subfigure[$L^1$ error\label{Fine-mesh-error-history-L1}.]
    {\includegraphics[width=0.35\textwidth,height= 0.3\textwidth]{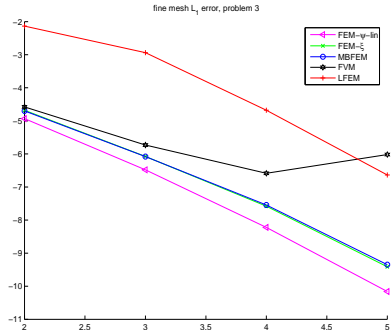}}
    \goodgap
    \subfigure[$L^2$ error\label{Fine-mesh-error-history-L2}.]
    {\includegraphics[width=0.35\textwidth,height= 0.3\textwidth]{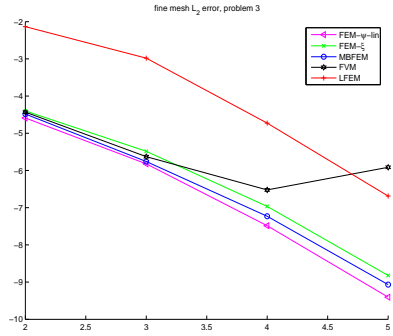}}\\
    \subfigure[$L^{\infty}$ error\label{Fine-mesh-error-history-Li}.]
    {\includegraphics[width=0.35\textwidth,height= 0.3\textwidth]{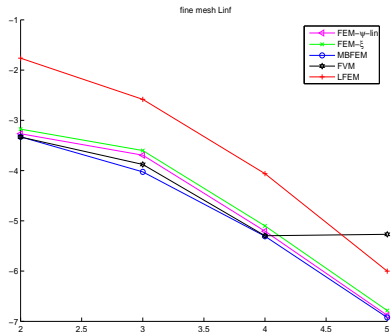}}
    \goodgap
    \subfigure[$H^1$ error\label{Fine-mesh-error-history-H1}.]
    {\includegraphics[width=0.35\textwidth,height= 0.3\textwidth]{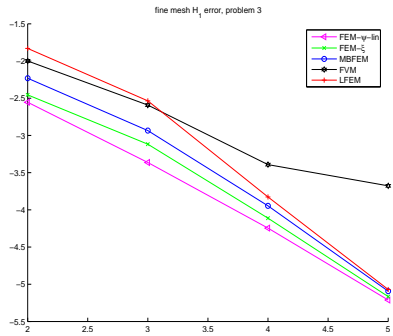}}\\[-10pt]
    \caption{\label{cap:Fine-mesh-error-history}Fine mesh approximation error
      ($\log_{2}$) $L^{1}$, $L^{2}$, $L^{\infty}$ and $H^{1}$ errors v.s coarse mesh refinement,
      Example \ref{exa:Fourier}, Random Fourier Modes.}
\end{center}
\end{figure}

\begin{example}
\label{exa:Fourier}Random Fourier modes.
\end{example}
In this case, $a(x)=e^{h(x)}$, with\[ h(x)=\sum_{|k|\leq
R}(a_{k}\sin(2\pi k\cdot x)+b_{k}\cos(2\pi k\cdot x))\] where
$a_{k}$ and $b_{k}$ are independent identically distributed random
variables on $[-0.3,0.3]$ and $R=6$. This is an other example where
scales are not separated. The weak aspect ratio of the triangles in
the metric induced by $F$ is $\eta^{\star}_{\min}=3.4997$. Observe
the deformation induced by the new metric (figure
\ref{cap:ExampleFourier}). Observe that distances between $u$ and
the interpolation of the coarse mesh approximations to the fine mesh
are larger (tables \ref{caswEswarrswaor1} and
\ref{caswErrsswwaor1}), this is due to the fact that those errors
depend on the aspect ratio $\eta^*_{\max}$ (which is not the case
for the coarse mesh errors) of course one could improve the
compression by adapting the mesh to the new metric but this has not
been our point of view here. We have preferred to show raw data
obtained with a given coarse mesh. The figures
\ref{cap:Coarse-mesh-error-history} and
\ref{cap:Fine-mesh-error-history} give the $L^{1}$, $L^{2}$,
$L^{\infty}$ and $H^{1}$ relative error ($\log_{2}$ basis versus
$\log_{2}$ basis of the resolution). The $x$-axis corresponds to the
refinement of coarse mesh, the $y$-axis is the error. The tables
\ref{cap:Coarse-mesh-solution-rate} and
\ref{cap:Fine-mesh-approximation-rate} give the convergence rate in
different norms (the parameter $\alpha$ in the error of the order of
$h^\alpha$).

\begin{figure}[httb]
  \begin{center}
    \subfigure[Condition Number.\label{condnumb3}]
    {\includegraphics[width=0.35\textwidth,height= 0.3\textwidth]{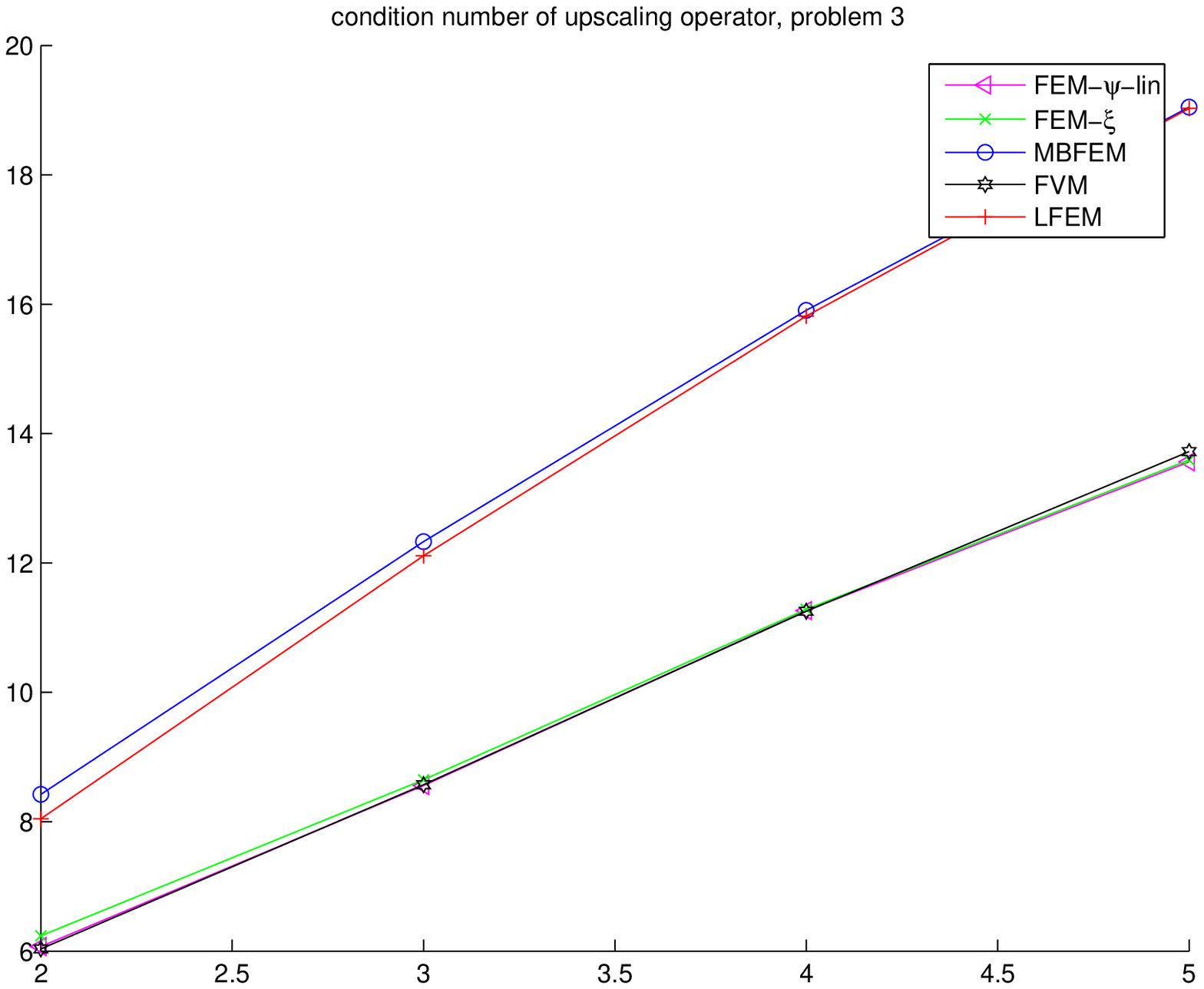}}
    \goodgap
    \subfigure[Coarse Mesh $L^1$ error. \label{l1err3}]
    {\includegraphics[width=0.35\textwidth,height= 0.3\textwidth]{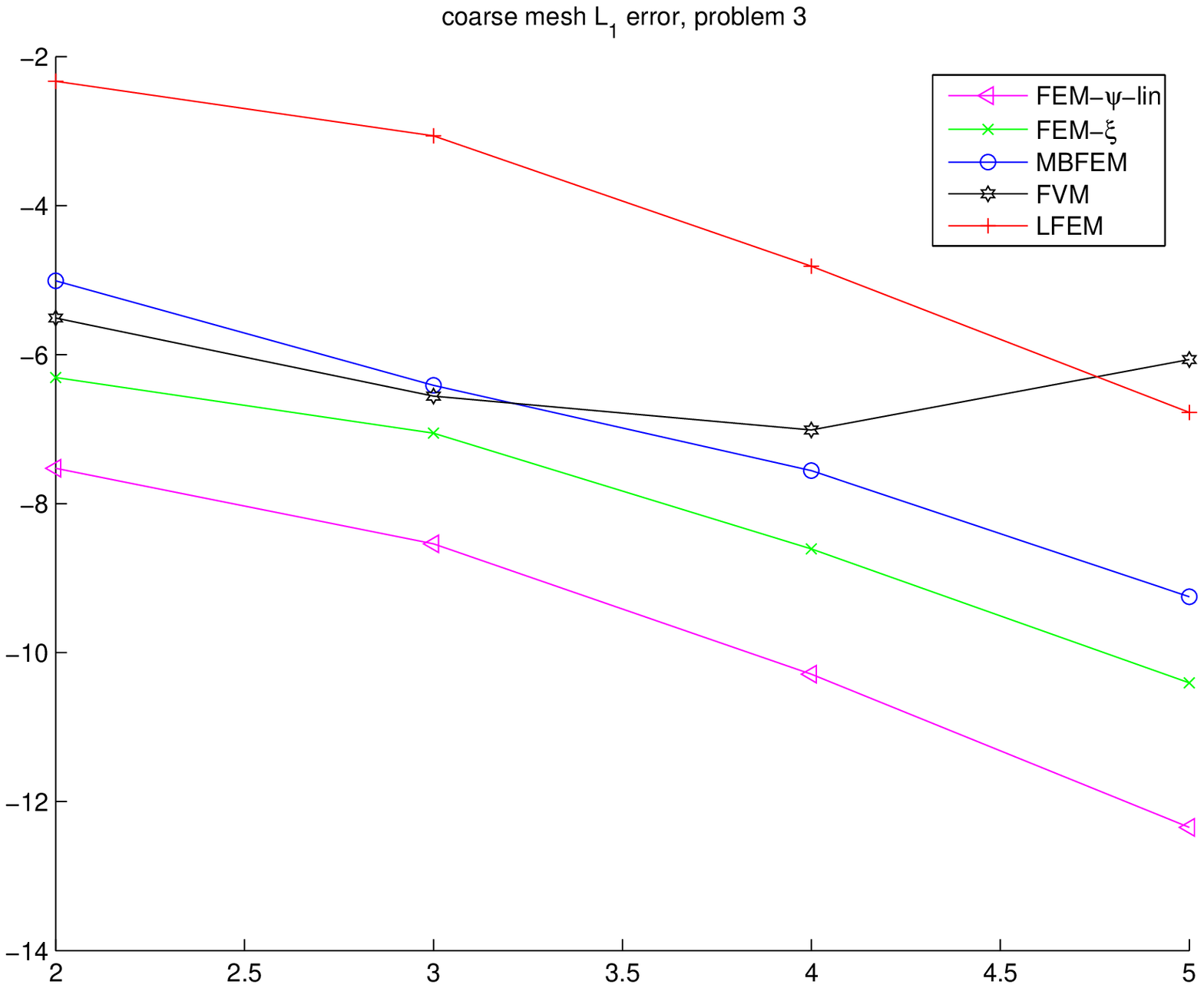}}\\
    \caption{\label{cap:ExamplesaFourier}Example \ref{exa:Fourier}, Random Fourier Modes }
\end{center}
\end{figure}

\begin{table}[!]
\begin{center}
\begin{tabular}{|c||c|c|c|c|c|}
\hline Coarse Mesh Error& FEM\_$\psi$& FEM\_$\xi$& MBFEM& FVM&
LFEM\tabularnewline \hline \hline
\begin{tabular}{c}
$L^{1}$\tabularnewline \hline $L^{2}$\tabularnewline \hline
$L^{\infty}$\tabularnewline \hline $H^{1}$\tabularnewline
\end{tabular}&
\begin{tabular}{c}
0.0027\tabularnewline \hline 0.0028\tabularnewline \hline
0.0066\tabularnewline \hline 0.0133\tabularnewline
\end{tabular}&
\begin{tabular}{c}
0.0075\tabularnewline \hline 0.0087\tabularnewline \hline
0.0278\tabularnewline \hline 0.0648\tabularnewline
\end{tabular}&
\begin{tabular}{c}
0.0117\tabularnewline \hline 0.0130\tabularnewline \hline
0.0320\tabularnewline \hline 0.0805\tabularnewline
\end{tabular}&
\begin{tabular}{c}
0.0106\tabularnewline \hline 0.0125\tabularnewline \hline
0.0376\tabularnewline \hline 0.0597\tabularnewline
\end{tabular}&
\begin{tabular}{c}
0.1197\tabularnewline \hline 0.1169\tabularnewline \hline
0.1358\tabularnewline \hline 0.1514\tabularnewline
\end{tabular}\tabularnewline
\hline
\end{tabular}
\caption{\label{caswEswarrswaor1}Example \ref{exa:Fourier}, Random
Fourier Modes }
\end{center}
\end{table}

\begin{table}[!]
\begin{center}
\begin{tabular}{|c||c|c|c|c|c|}
\hline Fine mesh Error& FEM\_$\psi$& FEM\_$\xi$& MBFEM& FVM&
LFEM\tabularnewline \hline \hline
\begin{tabular}{c}
$L^{1}$\tabularnewline \hline $L^{2}$\tabularnewline \hline
$L^{\infty}$\tabularnewline \hline $H^{1}$\tabularnewline
\end{tabular}&
\begin{tabular}{c}
0.0112\tabularnewline \hline 0.0177\tabularnewline \hline
0.0773\tabularnewline \hline 0.0972\tabularnewline
\end{tabular}&
\begin{tabular}{c}
0.0148\tabularnewline \hline 0.0223\tabularnewline \hline
0.0824\tabularnewline \hline 0.1152\tabularnewline
\end{tabular}&
\begin{tabular}{c}
0.0148\tabularnewline \hline 0.0184\tabularnewline \hline
0.0614\tabularnewline \hline 0.1307\tabularnewline
\end{tabular}&
\begin{tabular}{c}
0.0188\tabularnewline \hline 0.0202\tabularnewline \hline
0.0680\tabularnewline \hline 0.1659\tabularnewline
\end{tabular}&
\begin{tabular}{c}
0.1304\tabularnewline \hline 0.1265\tabularnewline \hline
0.1669\tabularnewline \hline 0.1725\tabularnewline
\end{tabular}\tabularnewline
\hline
\end{tabular}
\caption{\label{caswErrsswwaor1}Example \ref{exa:Fourier}, Random
Fourier Modes.}
\end{center}
\end{table}

\begin{table}[!]
\begin{center}
\begin{tabular}{|c|c|c|c|c|}
\hline Method& $L^{1}$& $L^{2}$& $L^{\infty}$&
$H^{1}$\tabularnewline \hline \hline FEM\_$\psi$& 1.62& 1.66& 1.56&
1.44\tabularnewline \hline FEM\_$\xi$& 1.38& 1.27& 1.23&
1.18\tabularnewline \hline MBFEM& 1.38& 1.40& 1.27&
1.08\tabularnewline \hline FVM& 0.53& 1.14& 1.26&
1.03\tabularnewline \hline LFEM& 1.51& 1.53& 1.62&
1.46\tabularnewline \hline
\end{tabular}

\caption{\label{cap:Coarse-mesh-solution-rate}Coarse mesh
approximation convergence rate}
\end{center}
\end{table}

\begin{table}[!]
\begin{center}
\begin{tabular}{|c|c|c|c|c|}
\hline Method& $L^{1}$& $L^{2}$& $L^{\infty}$&
$H^{1}$\tabularnewline \hline \hline FEM\_$\psi$& 1.74& 1.61& 1.23&
0.89\tabularnewline \hline FEM\_$\xi$& 1.57& 1.47& 1.23&
0.91\tabularnewline \hline MBFEM& 1.54& 1.52& 1.21&
0.96\tabularnewline \hline FVM& 0.75& 1.16& 1.22&
0.58\tabularnewline \hline LFEM& 1.52& 1.54& 1.42&
1.10\tabularnewline \hline
\end{tabular}

\caption{\label{cap:Fine-mesh-approximation-rate}Fine mesh
approximation convergence rate}
\end{center}
\end{table}

\clearpage

\begin{figure}[httb]
  \begin{center}
    \subfigure[Condition Number.\label{hddghd652}]
    {\includegraphics[width=0.35\textwidth,height= 0.3\textwidth]{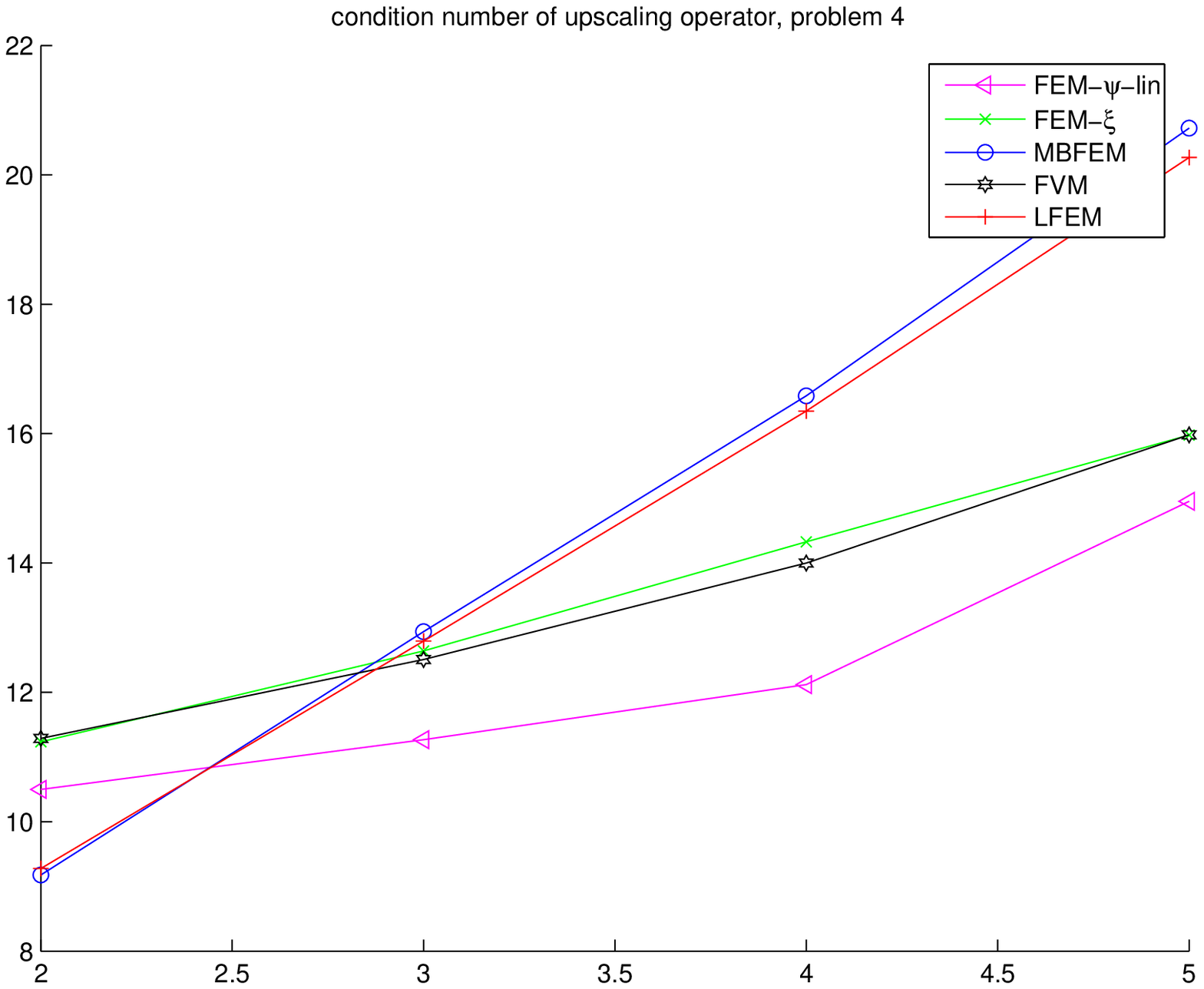}}
    \goodgap
    \subfigure[Coarse Mesh $L^1$ error.\label{hsyst6512}]
    {\includegraphics[width=0.35\textwidth,height= 0.3\textwidth]{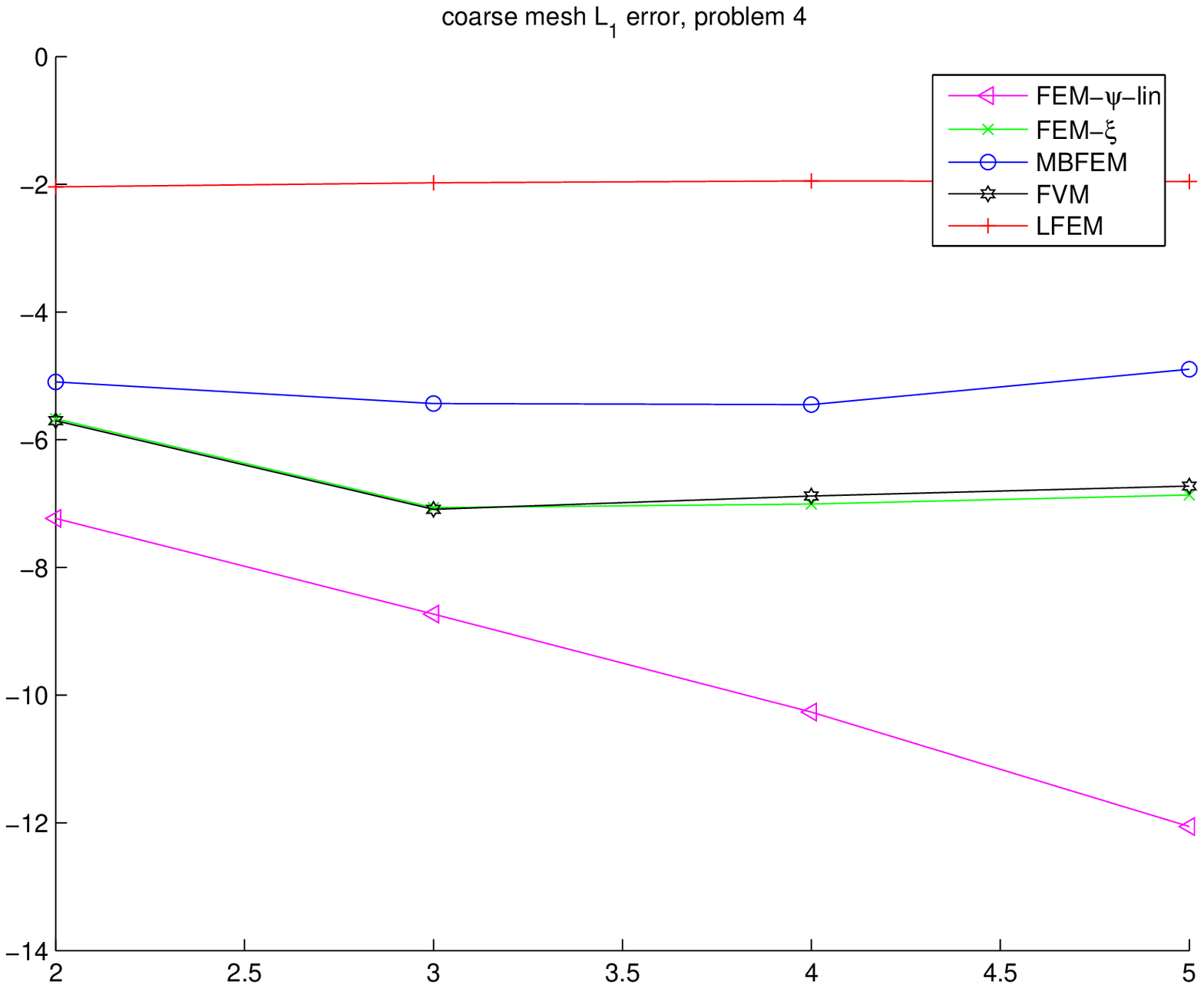}}\\
    \caption{\label{cap:Example nonergodics}Example \ref{exa:nonergodic}, Random Fractal}
\end{center}
\end{figure}

\begin{example}
\label{exa:nonergodic} Random fractal
\end{example}
In this case, $a$ is given by a product of discontinuous functions
oscillating randomly at different scales,
$a(x)=a_{1}(x)a_{2}(x)\cdots a_{n}(x)$, and $a_{i}(x)=c_{pq}$ for
$x\in[\frac{p}{2^{i}},\frac{p+1}{2^{i}})\times[\frac{q}{2^{i}},\frac{q+1}{2^{i}})$,
$c_{pq}$ is uniformly random in $[\frac{1}{\gamma},\gamma]$,  $n=5,$
and $\gamma=2$. The weak aspect ratio is
$\eta^{\star}_{\min}=2.4796$.  Table \ref{caswaqwarrswaor1}
 gives the relative error estimated on the nodes of the coarse mesh between the solution $u$ of the initial PDE
 \eref{ghjh52} and an approximation obtained from the up-scaled
 operator on the nodes of the coarse mesh. Table \ref{caswaqwarrsewoew3r1}
 gives
 the relative error estimated on the nodes of the fine mesh between  $u$ and the  $\J_h$-interpolation of the previous approximations with respect
 to $F$ on a fine resolution. Figure \ref{hddghd652} gives the condition number of the
stiffness matrix associated to the up-scaled operator versus
$-\log_2 h$ (logarithm of the resolution). Figure \ref{hsyst6512}
gives the relative $L_1$-distance between $u$ and its approximation
on the coarse mesh in log scale versus $-\log_2 h$.

\begin{figure}[!]
  \begin{center}
    \subfigure[ $a$.]
    {\includegraphics[width=0.35\textwidth,height= 0.3\textwidth]{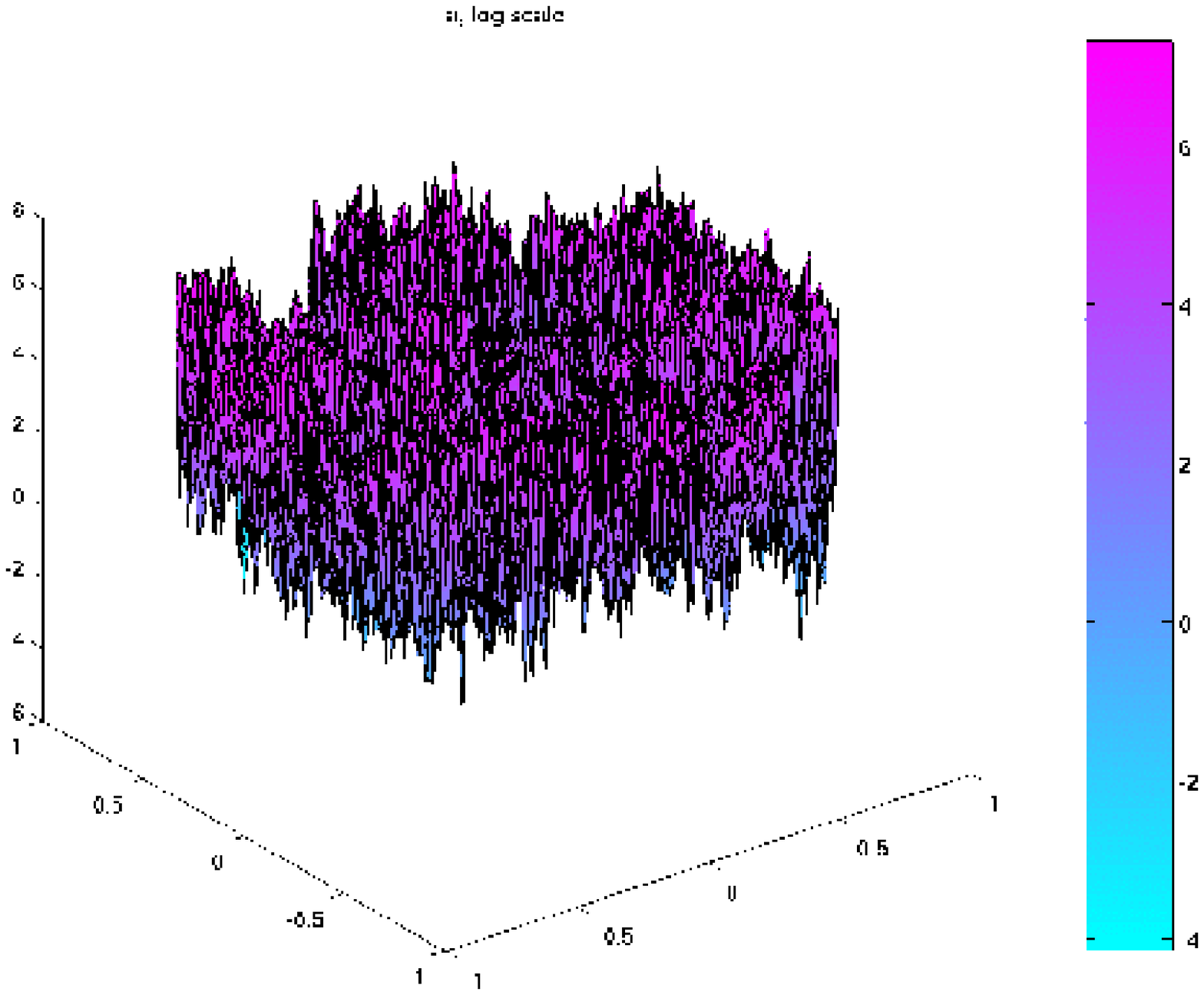}}
    \goodgap
    \subfigure[$\T^F$]
    {\includegraphics[width=0.35\textwidth,height= 0.3\textwidth]{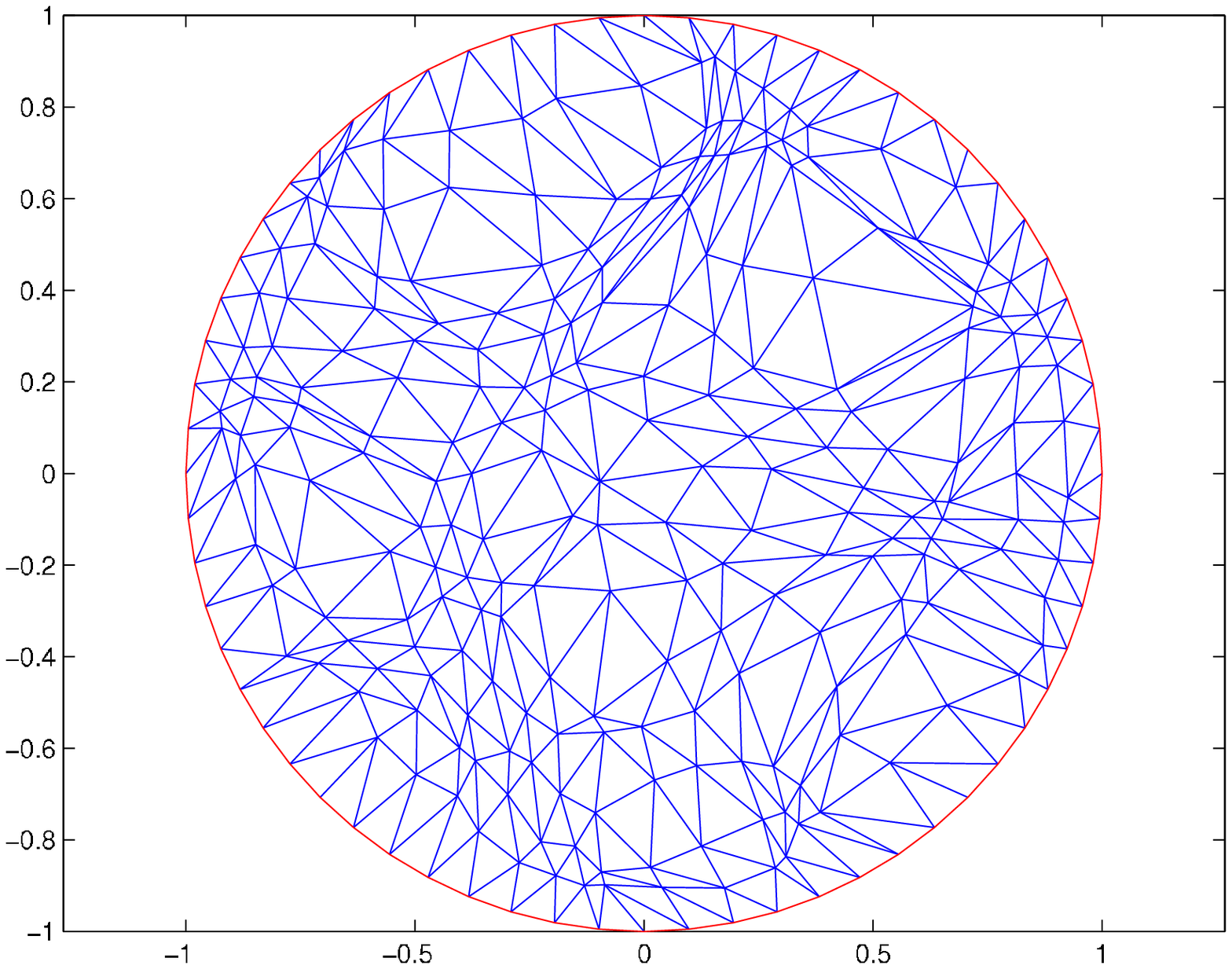}}\\
    \caption{\label{cap:Example nonergodic}Example \ref{exa:nonergodic}, Random Fractal.}
\end{center}
\end{figure}

\begin{table}[!]
\begin{center}
\begin{tabular}{|c||c|c|c|c|c|}
\hline Coarse Mesh Error& FEM\_$\psi$& FEM\_$\xi$& MBFEM& FVM&
LFEM\tabularnewline \hline \hline
\begin{tabular}{c}
$L^{1}$\tabularnewline \hline $L^{2}$\tabularnewline \hline
$L^{\infty}$\tabularnewline \hline $H^{1}$\tabularnewline
\end{tabular}&
\begin{tabular}{c}
0.0024\tabularnewline \hline 0.0025\tabularnewline \hline
0.0094\tabularnewline \hline 0.0161\tabularnewline
\end{tabular}&
\begin{tabular}{c}
0.0075\tabularnewline \hline 0.0085\tabularnewline \hline
0.0399\tabularnewline \hline 0.0718\tabularnewline
\end{tabular}&
\begin{tabular}{c}
0.0231\tabularnewline \hline 0.0241\tabularnewline \hline
0.0920\tabularnewline \hline 0.1553\tabularnewline
\end{tabular}&
\begin{tabular}{c}
0.0073\tabularnewline \hline 0.0100\tabularnewline \hline
0.0398\tabularnewline \hline 0.0493\tabularnewline
\end{tabular}&
\begin{tabular}{c}
0.0519\tabularnewline \hline 0.0606\tabularnewline \hline
0.1694\tabularnewline \hline 0.3107\tabularnewline
\end{tabular}\tabularnewline
\hline
\end{tabular}%
\caption{\label{caswaqwarrswaor1}Example \ref{exa:nonergodic},
Random Fractal.}
\end{center}
\end{table}

\begin{table}[!]
\begin{center}
\begin{tabular}{|c||c|c|c|c|c|}
\hline Fine mesh Error& FEM\_$\psi$& FEM\_$\xi$& MBFEM& FVM&
LFEM\tabularnewline \hline \hline
\begin{tabular}{c}
$L^{1}$\tabularnewline \hline $L^{2}$\tabularnewline \hline
$L^{\infty}$\tabularnewline \hline $H^{1}$\tabularnewline
\end{tabular}&
\begin{tabular}{c}
0.0108\tabularnewline \hline 0.0155\tabularnewline \hline
0.0662\tabularnewline \hline 0.1015\tabularnewline
\end{tabular}&
\begin{tabular}{c}
0.0147\tabularnewline \hline 0.0198\tabularnewline \hline
0.0802\tabularnewline \hline 0.1231\tabularnewline
\end{tabular}&
\begin{tabular}{c}
0.0245\tabularnewline \hline 0.0280\tabularnewline \hline
0.0919\tabularnewline \hline 0.1838\tabularnewline
\end{tabular}&
\begin{tabular}{c}
0.0142\tabularnewline \hline 0.0173\tabularnewline \hline
0.0720\tabularnewline \hline 0.1433\tabularnewline
\end{tabular}&
\begin{tabular}{c}
0.0765\tabularnewline \hline 0.0812\tabularnewline \hline
0.1694\tabularnewline \hline 0.2642\tabularnewline
\end{tabular}\tabularnewline
\hline
\end{tabular}
\caption{\label{caswaqwarrsewoew3r1}Example \ref{exa:nonergodic},
Random Fractal.}
\end{center}
\end{table}

\clearpage

\begin{figure}[httb]
  \begin{center}
    \subfigure[Condition Number.\label{hdhgd632}]
    {\includegraphics[width=0.35\textwidth,height= 0.3\textwidth]{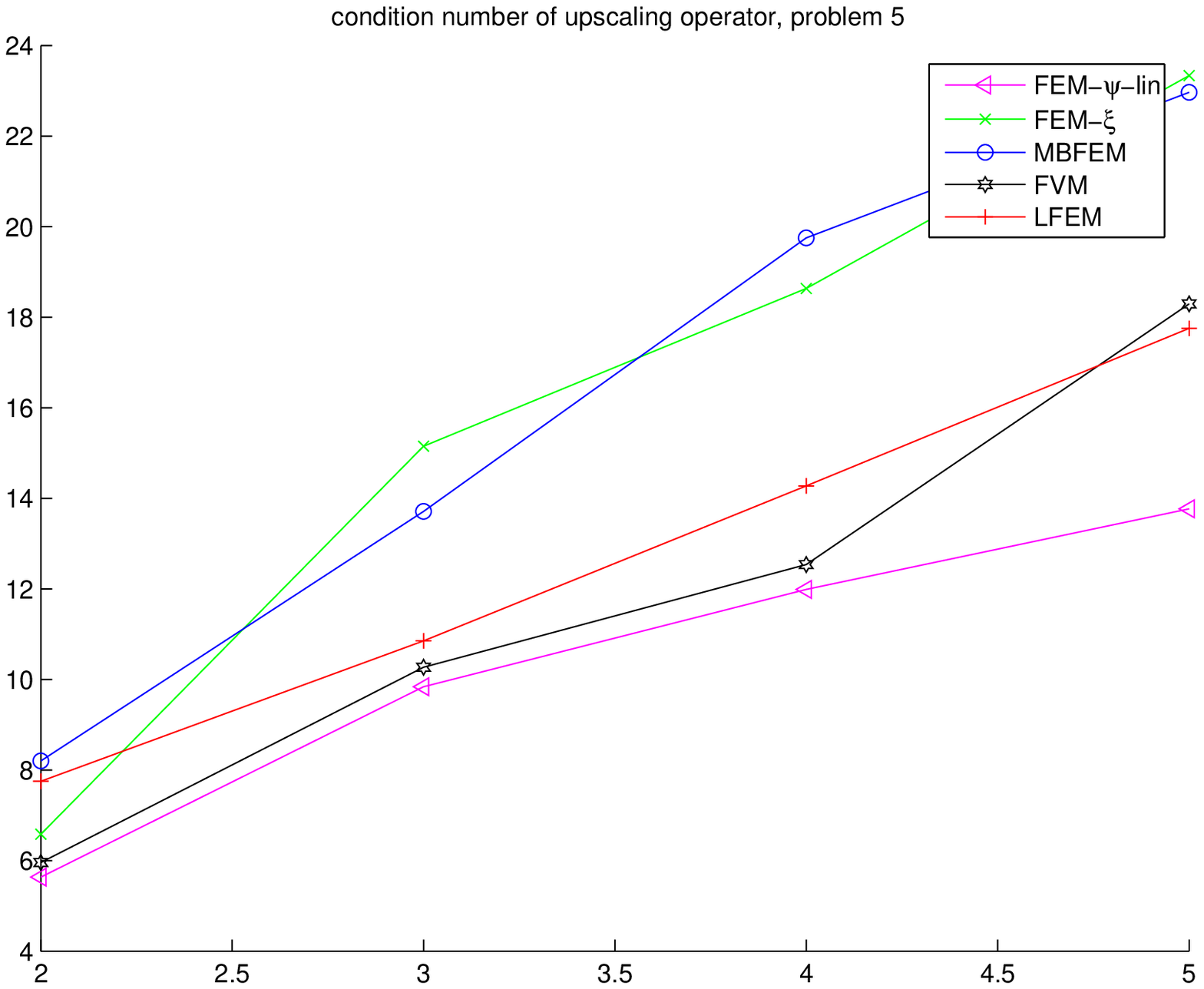}}
    \goodgap
    \subfigure[Coarse Mesh $L^1$ error.\label{hgs637hd}]
    {\includegraphics[width=0.35\textwidth,height= 0.3\textwidth]{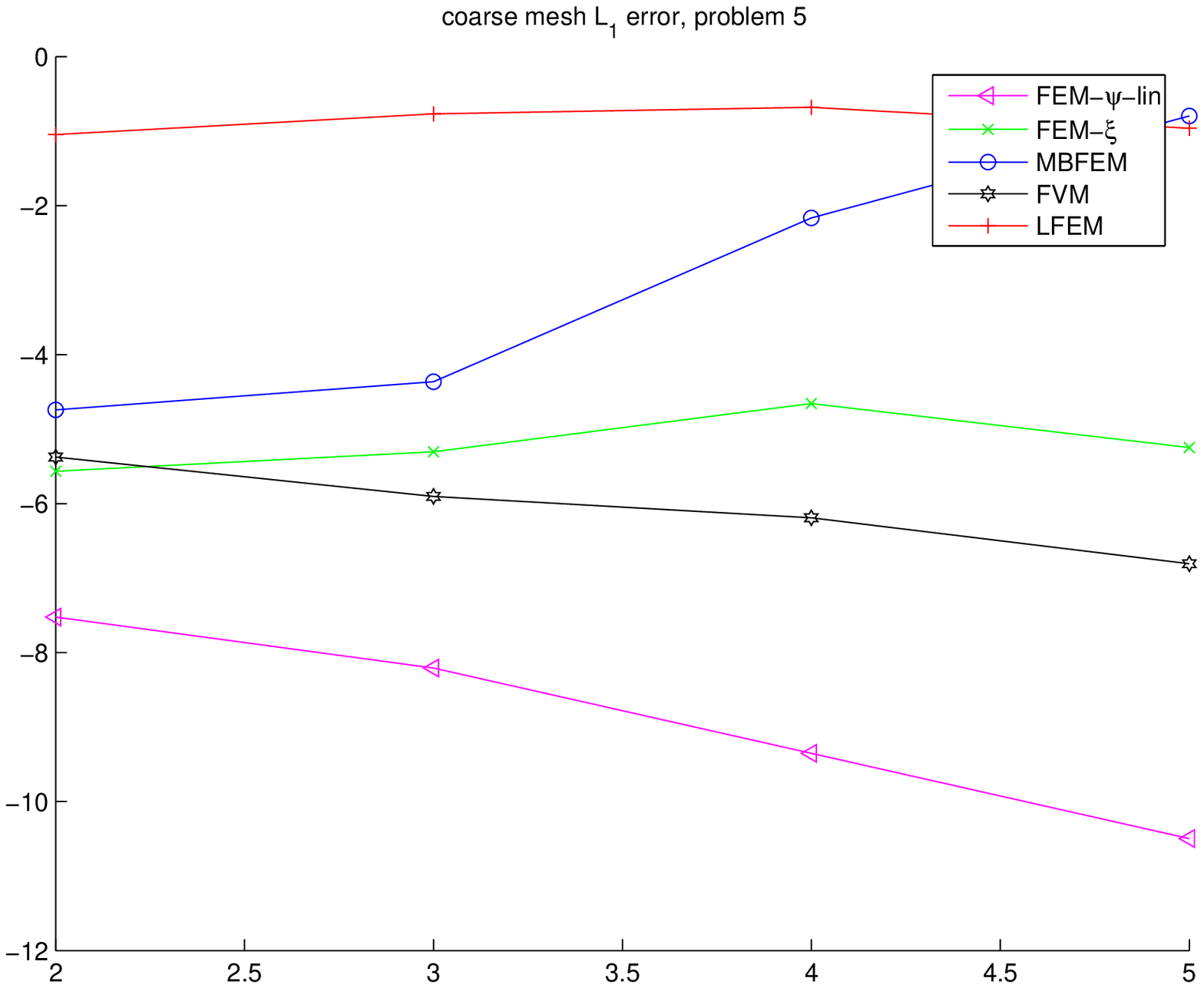}}
\caption{Example \ref{exa:percolation}, Percolation }
\label{cap:Example Percolations}
\end{center}
\end{figure}

\begin{example}
\label{exa:percolation}Percolation at criticality
\end{example}
In this case, the conductivity of each site is equal to $\gamma$ or
$1/\gamma$ with probability $1/2$. We have chosen $\gamma=4$ in this
example. Observe that some errors are larger for this challenging
case because a percolating medium generates flat triangles in the
new metric, indeed $\eta^{\star}_{\min}=22.3395$. Table
\ref{caswaqwarrseaswaor1}
 gives the relative error estimated on the nodes of the coarse mesh between the solution $u$ of the initial PDE
 \eref{ghjh52} and an approximation obtained from the up-scaled
 operator on the nodes of the coarse mesh. Table \ref{caswaqwarrsaawswaor1}
 gives
 the relative error estimated on the nodes of the fine mesh between  $u$ and the  $\J_h$-interpolation of the previous approximations with respect
 to $F$ on a fine resolution. Figure \ref{hdhgd632} gives the condition number of the
stiffness matrix associated to the up-scaled operator versus
$-\log_2 h$ (logarithm of the resolution). Figure \ref{hgs637hd}
gives the relative $L_1$-distance between $u$ and its approximation
on the coarse mesh in log scale versus $-\log_2 h$. Observe that the
methods based on a global change of metric do converge but when that
numerical homogenization is done by computing only local coarse
parameters (in averaging or finite elements techniques) then
convergence is not guaranteed without further assumptions on $a$
(see the curve of LFEM).

\begin{figure}[!]
  \begin{center}
    \subfigure[$a$.]
    {\includegraphics[width=0.4\textwidth,height= 0.35\textwidth]{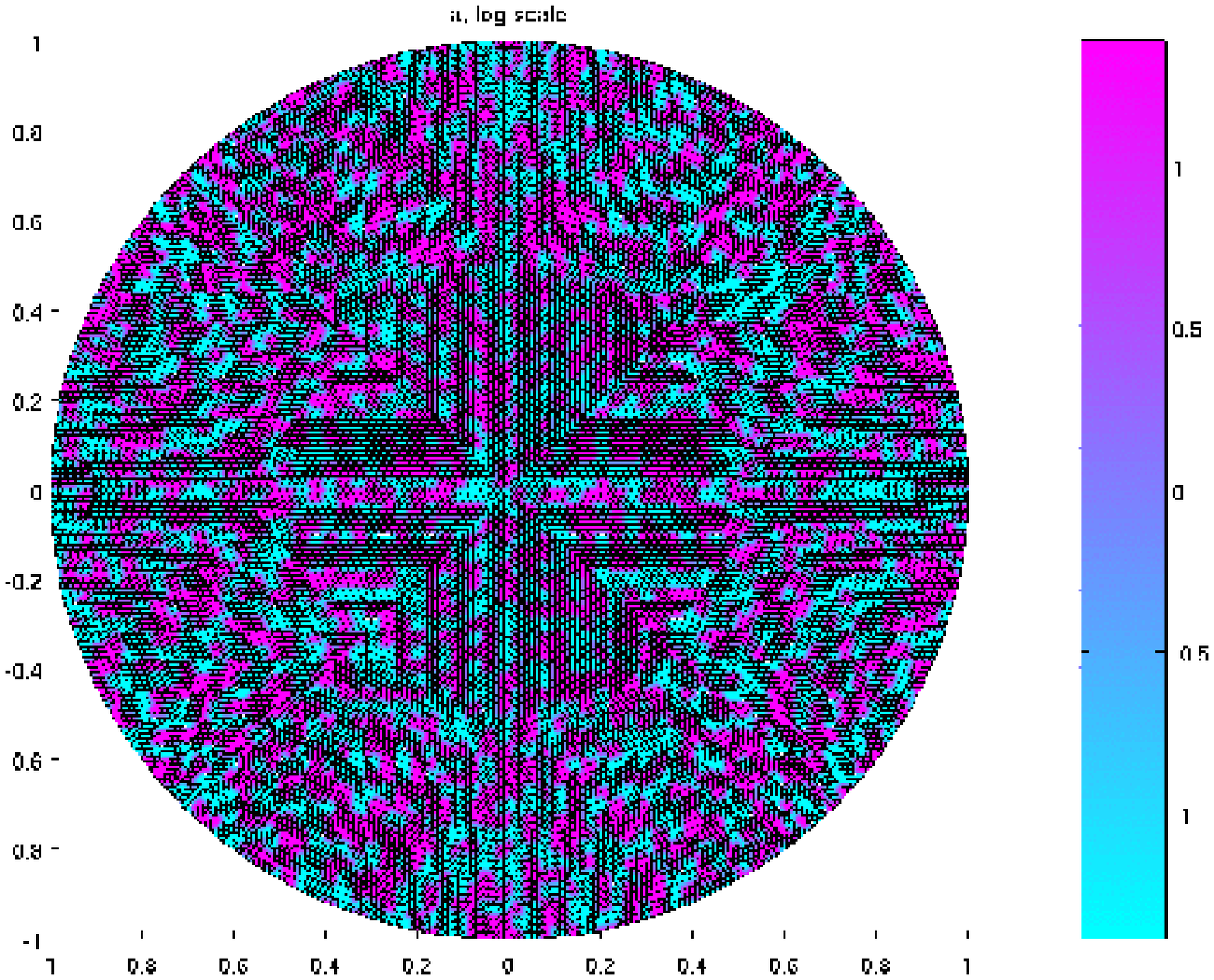}}
    \subfigure[$\T^F$.]
    {\includegraphics[width=0.4\textwidth,height= 0.35\textwidth]{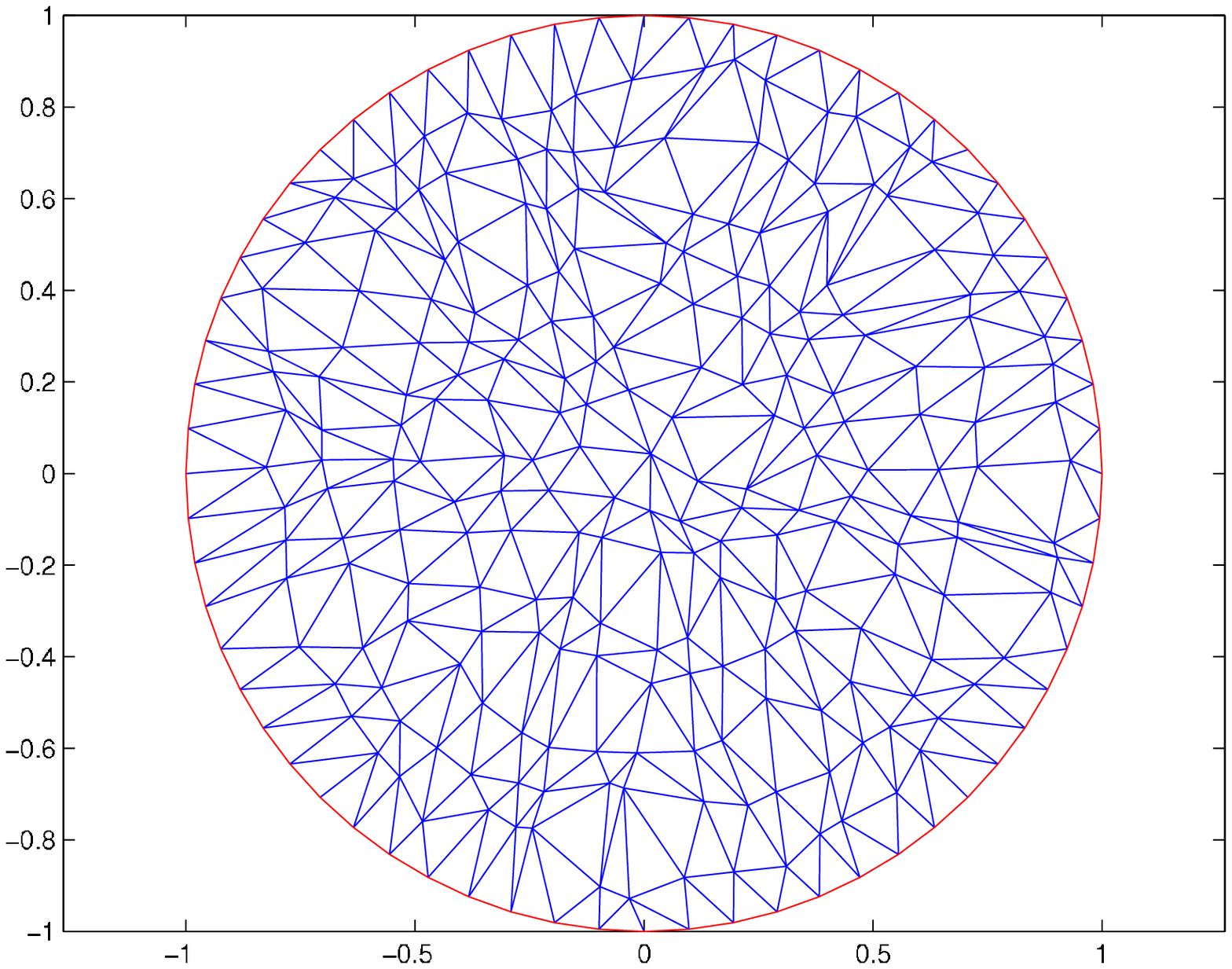}}
\caption{Example \ref{exa:percolation}, Percolation.}
\label{cap:Example Percolation}
\end{center}
\end{figure}

\begin{table}[!]
\begin{center}
\begin{tabular}{|c||c|c|c|c|c|}
\hline Coarse Mesh Error& FEM\_$\psi$& FEM\_$\xi$& MBFEM& FVM&
LFEM\tabularnewline \hline \hline
\begin{tabular}{c}
$L^{1}$\tabularnewline \hline $L^{2}$\tabularnewline \hline
$L^{\infty}$\tabularnewline \hline $H^{1}$\tabularnewline
\end{tabular}&
\begin{tabular}{c}
0.0034\tabularnewline \hline 0.0041\tabularnewline \hline
0.0163\tabularnewline \hline 0.0343\tabularnewline
\end{tabular}&
\begin{tabular}{c}
0.0253\tabularnewline \hline 0.0265\tabularnewline \hline
0.0813\tabularnewline \hline 0.0843\tabularnewline
\end{tabular}&
\begin{tabular}{c}
0.0485\tabularnewline \hline 0.0523\tabularnewline \hline
0.0643\tabularnewline \hline 0.1070\tabularnewline
\end{tabular}&
\begin{tabular}{c}
0.0167\tabularnewline \hline 0.0189\tabularnewline \hline
0.0499\tabularnewline \hline 0.0713\tabularnewline
\end{tabular}&
\begin{tabular}{c}
0.2848\tabularnewline \hline 0.2851\tabularnewline \hline
0.3018\tabularnewline \hline 0.3740\tabularnewline
\end{tabular}\tabularnewline
\hline
\end{tabular}
\caption{\label{caswaqwarrseaswaor1}Example \ref{exa:percolation},
Percolation.}
\end{center}
\end{table}

\begin{table}[!]
\begin{center}
\begin{tabular}{|c||c|c|c|c|c|}
\hline Fine mesh Error& FEM\_$\psi$& FEM\_$\xi$& MBFEM& FVM&
LFEM\tabularnewline \hline \hline
\begin{tabular}{c}
$L^{1}$\tabularnewline \hline $L^{2}$\tabularnewline \hline
$L^{\infty}$\tabularnewline \hline $H^{1}$\tabularnewline
\end{tabular}&
\begin{tabular}{c}
0.0115\tabularnewline \hline 0.0152\tabularnewline \hline
0.0500\tabularnewline \hline 0.1000\tabularnewline
\end{tabular}&
\begin{tabular}{c}
0.0265\tabularnewline \hline 0.0268\tabularnewline \hline
0.0527\tabularnewline \hline 0.1712\tabularnewline
\end{tabular}&
\begin{tabular}{c}
0.0585\tabularnewline \hline 0.0628\tabularnewline \hline
0.0940\tabularnewline \hline 0.1954\tabularnewline
\end{tabular}&
\begin{tabular}{c}
0.0216\tabularnewline \hline 0.0229\tabularnewline \hline
0.0497\tabularnewline \hline 0.1343\tabularnewline
\end{tabular}&
\begin{tabular}{c}
0.3024\tabularnewline \hline 0.3015\tabularnewline \hline
0.3135\tabularnewline \hline 0.3964\tabularnewline
\end{tabular}\tabularnewline
\hline
\end{tabular}
\caption{\label{caswaqwarrsaawswaor1}Example \ref{exa:percolation},
Percolation.}
\end{center}
\end{table}

\clearpage

\subsection{Numerical experiments with splines.}
We have seen that if $\sigma$ is stable then $u\circ {F^{-1}}$
belongs to $W^{2,p}(\Omega)$ with $p>2$. It is thus natural to
expect a better accuracy by using $C^1$-continuous elements in the
method described in subsection \ref{gfssg5} instead of piecewise
linear elements. This increase of accuracy has already been observed
in \cite{AlBr04} when $F$ is obtained as the solution of a local
cell problem. In our case (the harmonic coordinates are computed
globally) we also observe a sharp increase of the accuracy of finite
element method of subsection \ref{gfssg5} by using splines for the
elements $\varphi_i$.

We refer to \cite{key-2,key-3} for methods using $C^{1}$ finite
element.  One possibility is the weighted extened B-splines (WEB)
method developed by K. H$\ddot{o}$llig in \cite{key-1,key-2}, these
elements are $C^2$-continuous. They are obtained from tensor
products of one dimensional elements. The Dirichlet boundary
condition is satisfied using a smooth weight function $\omega$, such
that $\omega=0$ at the boundary. The condition number of the
stiffness matrix is bounded from above by $O(h^{-2})$ (we have the
same optimal bound on a Galerkin system with piecewise linear
elements).

We have considered two challenging multi-scale medium for our
numerical experiments:  random Fourier modes and percolation. For
the simplicity of the implementation a square domain has been
considered, and weighted spline basis are used instead of the WEB
spline basis. For a square domain $[-1,1]\times[-1,1]$, the weight
is $\omega=(1-x^2)(1-y^2)$. Two methods have been compared,
\begin{itemize}
\item The Galerkin scheme using the finite elements
$\psi_i=\varphi_{i}\circ F$, where $\varphi_i$ are the piecewise
linear nodal basis elements of subsection \ref{gfssg5}, noted
FEM\_$\psi_{lin}$
\item The Galerkin scheme using the  finite element
$\psi_i=\phi_{i}\circ F$, where $\phi_i$ are weighted cubic B-spline
basis elements, noted FEM\_$\psi_{sp}$
\end{itemize}

\begin{figure}[httb]
  \begin{center}
    \subfigure[Condition Number.\label{condnusmber3q}]
    {\includegraphics[width=0.3\textwidth,height= 0.3\textwidth]{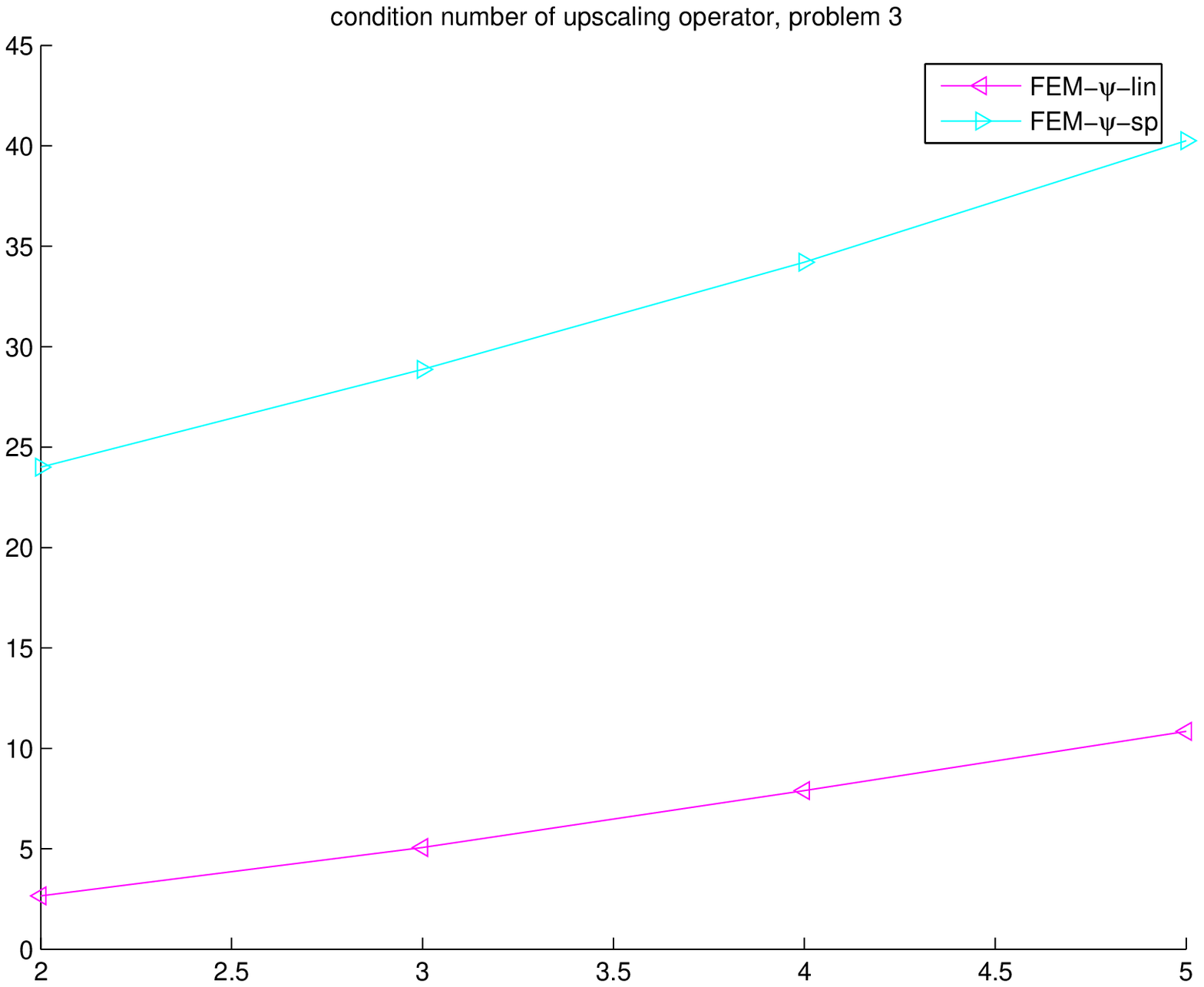}}
    \goodgap
    \subfigure[Coarse Mesh $L^1$ error. \label{l1cersr3q}]
    {\includegraphics[width=0.3\textwidth,height= 0.3\textwidth]{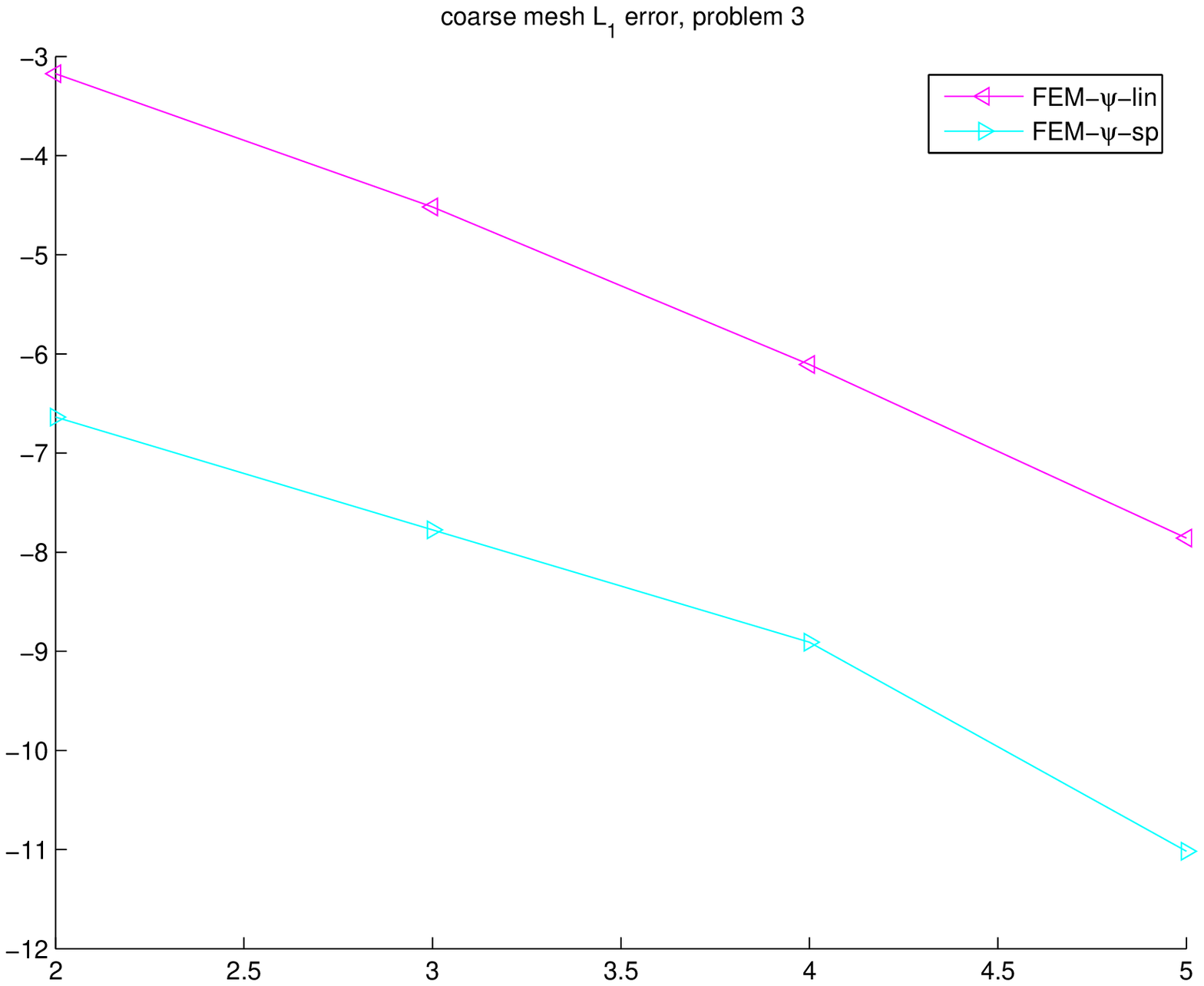}}
    \goodgap
    \subfigure[Fine Mesh $L^1$ error. \label{l1fesrr3q}]
    {\includegraphics[width=0.3\textwidth,height= 0.3\textwidth]{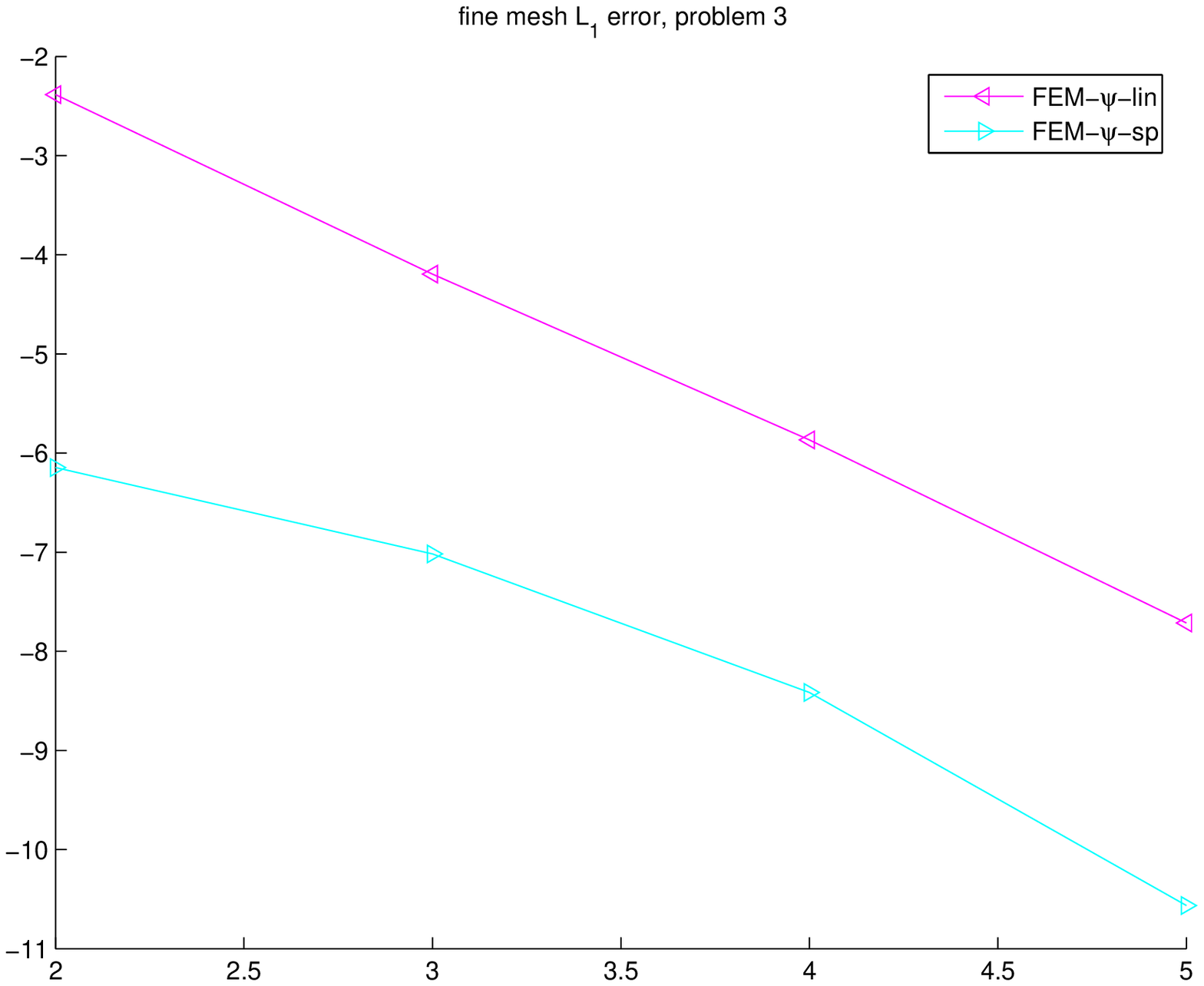}}\\

    \caption{\label{cap:ExamplesaFouwrier}Example \ref{exa:Fourier}. Random Fourier Modes.}
\end{center}
\end{figure}

\begin{figure}[httb]
  \begin{center}
    \subfigure[Condition Number.\label{condnumber7q}]
    {\includegraphics[width=0.3\textwidth,height= 0.3\textwidth]{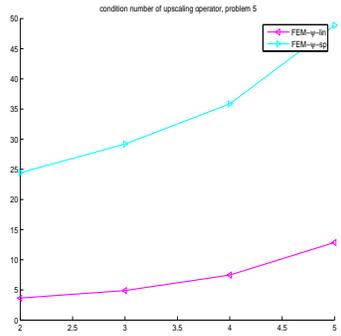}}
    \goodgap
    \subfigure[Coarse Mesh $L^1$ error. \label{l1cerr3q}]
    {\includegraphics[width=0.3\textwidth,height= 0.3\textwidth]{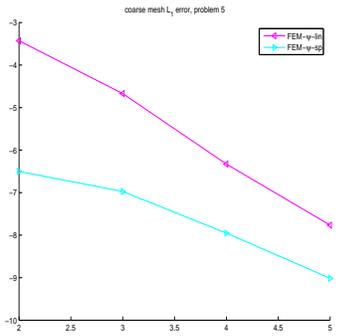}}
    \goodgap
    \subfigure[Fine Mesh $L^1$ error. \label{l1ferr3q}]
    {\includegraphics[width=0.35\textwidth,height= 0.3\textwidth]{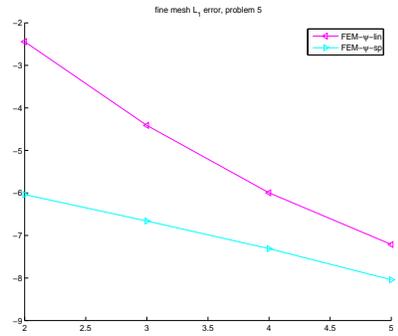}}\\

    \caption{\label{cap:ExamplesswaFourier}Example \ref{exa:percolation}. Percolation at criticality.}
\end{center}
\end{figure}

The error obtain with the method FEM\_$\psi_{sp}$ is much smaller
than the one obtained with the method FEM\_$\psi_{lin}$ as it is
shown in figures \ref{cap:ExamplesaFouwrier} and
\ref{cap:ExamplesswaFourier}.

\begin{table}[!]
\begin{center}
\begin{tabular}{|c||c|c|}
\hline Coarse Mesh Error& FEM\_$\psi_{lin}$& FEM\_$\psi_{sp}$
\tabularnewline \hline \hline
\begin{tabular}{c}
$L^{1}$\tabularnewline \hline $L^{2}$\tabularnewline \hline
$L^{\infty}$\tabularnewline \hline $H^{1}$\tabularnewline
\end{tabular}&
\begin{tabular}{c}
0.0437\tabularnewline \hline 0.0426\tabularnewline \hline
0.0614\tabularnewline \hline 0.0746\tabularnewline
\end{tabular}&
\begin{tabular}{c}
0.0046\tabularnewline \hline 0.0052\tabularnewline \hline
0.0096\tabularnewline \hline 0.0227\tabularnewline
\end{tabular}
\tabularnewline \hline
\end{tabular}
\caption{\label{caswEswarrswafdor1op}Example \ref{exa:Fourier},
Random Fourier Modes }
\end{center}
\end{table}

\begin{table}[!]
\begin{center}
\begin{tabular}{|c||c|c|}
\hline Fine mesh Error& FEM\_$\psi_{lin}$&
FEM\_$\psi_{sp}$\tabularnewline \hline \hline
\begin{tabular}{c}
$L^{1}$\tabularnewline \hline $L^{2}$\tabularnewline \hline
$L^{\infty}$\tabularnewline \hline $H^{1}$\tabularnewline
\end{tabular}&
\begin{tabular}{c}
0.0546\tabularnewline \hline 0.0529\tabularnewline \hline
0.0920\tabularnewline \hline 0.2109\tabularnewline
\end{tabular}&
\begin{tabular}{c}
0.0077\tabularnewline \hline 0.0096\tabularnewline \hline
0.0289\tabularnewline \hline 0.0547\tabularnewline
\end{tabular}
\tabularnewline \hline
\end{tabular}
\caption{\label{caswErrsswwavor1op}Example \ref{exa:Fourier}, Random
Fourier Modes.}
\end{center}
\end{table}

\begin{table}[!]
\begin{center}
\begin{tabular}{|c||c|c|}
\hline Coarse Mesh Error& FEM\_$\psi_{lin}$& FEM\_$\psi_{sp}$
\tabularnewline \hline \hline
\begin{tabular}{c}
$L^{1}$\tabularnewline \hline $L^{2}$\tabularnewline \hline
$L^{\infty}$\tabularnewline \hline $H^{1}$\tabularnewline
\end{tabular}&
\begin{tabular}{c}
0.0393\tabularnewline \hline 0.0379\tabularnewline \hline
0.0622\tabularnewline \hline 0.0731\tabularnewline
\end{tabular}&
\begin{tabular}{c}
0.0080\tabularnewline \hline 0.0098\tabularnewline \hline
0.0309\tabularnewline \hline 0.0404\tabularnewline
\end{tabular}
\tabularnewline \hline
\end{tabular}
\caption{\label{caswEswarrswaonr1op}Example \ref{exa:percolation},
Percolation }
\end{center}
\end{table}

\begin{table}[!]
\begin{center}
\begin{tabular}{|c||c|c|}
\hline Fine mesh Error& FEM\_$\psi_{lin}$&
FEM\_$\psi_{sp}$\tabularnewline \hline \hline
\begin{tabular}{c}
$L^{1}$\tabularnewline \hline $L^{2}$\tabularnewline \hline
$L^{\infty}$\tabularnewline \hline $H^{1}$\tabularnewline
\end{tabular}&
\begin{tabular}{c}
0.0470\tabularnewline \hline 0.0464\tabularnewline \hline
0.1174\tabularnewline \hline 0.2030\tabularnewline
\end{tabular}&
\begin{tabular}{c}
0.0099\tabularnewline \hline 0.0130\tabularnewline \hline
0.0554\tabularnewline \hline 0.0838\tabularnewline
\end{tabular}
\tabularnewline \hline
\end{tabular}
\caption{\label{caswErrsswwaor1opb}Example \ref{exa:percolation},
Percolation.}
\end{center}
\end{table}








%








\clearpage

\paragraph*{Acknowledgments.}
Part of the work of the first author has been supported by CNRS. The
authors would like to thank Jean-Michel Roquejoffre for indicating
us the correct references on nonlinear PDEs, Mathieu Desbrun for
enlightening discussions on discrete exterior calculus
\cite{MR2047000} (a powerful tool that has put into evidence the
intrinsic way to define discrete differential operators on irregular
triangulations), Tom Hou and Jerry Marsden for stimulating
discussions on multi-scale computation, Clothilde Melot and
St\'{e}phane Jaffard for stimulating discussions on multi-fractal
analysis. Thanks are also due to Lexing Ying and Laurent Demanet for
useful comments on the manuscript and G. Ben Arous for indicating us
reference \cite{MR1473568}. Many thanks are also due to Stefan
M\"{u}ller (MPI, Leipzig) for valuable suggestions and for
indicating us the Hierarchical Matrices methods. We would like also
to thank G. Allaire, F. Murat and S.R.S. Varadhan for stimulating
discussions at the CIRM workshop on random homogenization and P.
Schr\"{o}der for stimulating discussions on splines based methods.
We also thank an anonymous referee for detailed comments and
suggestions.

\end{document}